\documentclass{amsart}

\usepackage{amsmath,amsfonts, latexsym,graphicx, footmisc, amssymb, mathtools, enumerate, amscd}

\usepackage{hyperref}
\hypersetup{colorlinks=true, urlcolor=blue, citecolor=blue, linkcolor=blue}

\newcommand{\Adot}{{\mathbf A^\bullet}}
\newcommand{\Bdot}{{\mathbf B^\bullet}}
\newcommand{\Cdot}{{\mathbf C^\bullet}}
\newcommand{\Edot}{{\mathbf E^\bullet}}
\newcommand{\Fdot}{{\mathbf F^\bullet}}
\newcommand{\Gdot}{{\mathbf G^\bullet}}
\newcommand{\Idot}{{\mathbf I^\bullet}}
\newcommand{\Kdot}{{\mathbf K^\bullet}}
\newcommand{\Mdot}{{\mathbf M^\bullet}}
\newcommand{\Pdot}{{\mathbf P^\bullet}}
\newcommand{\Qdot}{{\mathbf Q^\bullet}}
\newcommand{\C}{{\mathbb C}}
\newcommand{\R}{{\mathbb R}}
\newcommand{\D}{{\mathbb D}}
\newcommand{\vdual}{{\mathcal D}}

\newcommand{\hyp}{{\mathbb H}}

\newcommand{\dm}{{\operatorname{dim}}}

\newcommand{\pcoh}{{{}^{\mu}\hskip -0.02in H^0}}

\newcommand{\map}{{{f\circ\pi_1+g\circ\pi_2}}}
\newcommand{\lotimes}{{\ {\otimes^L}\ }}
\newcommand{\lboxtimes}{{\ {\boxtimes^L}\ }}
\newcommand{\piten}{{{\pi_1^*\Adot\lotimes\pi_2^*\Bdot}}}

\title{Notes on Perverse Sheaves and Vanishing Cycles}

\author{David B. Massey}

\address{David B. Massey, Dept. of Mathematics, Northeastern University, Boston,
MA, 02115, USA} 

\email{d.massey@northeastern.edu}

\begin{document}

\maketitle

\vskip .2in

\noindent\S0. {\bf Introduction to Version 10-26}  

\vskip .2in

These notes are my continuing effort to provide a sort of working mathematician's
guide to the derived category and perverse sheaves.  They began in 1991 as handwritten notes for David Mond
and his students, and then I decided to type them for possible inclusion as an appendix in a paper or book. 
Now, however, these notes represent more of a journal of my understanding of this machinery;
whenever I believe that I have understood a new, significant chunk of the subject, I put it here -- so
that the next time that I want to understand or use that chunk, I will only have to look in one
place.  Moreover, the place that I have to look will be written in a manner that I can follow easily
(and that I can carry around easily).  The version number indicates the month and year when I last
added to this effort. I should mention that the August 2025 update contains a major revision of the section on nearby and vanishing cycles.

\vskip .1in

The only results of my own that appear here are the Sebastiani-Thom Isomorphism and the fact that Verdier dualizing commutes with the shifted nearby and vanishing cycles up to natural isomorphisms. For the most part, I have merely attempted to pull together some results from a number of sources.  Primary sources
are  \cite{BBD},\cite{G-M3}, \cite{H}, \cite{I}, \cite{K-S}, \cite{Mac1}, \cite{M-V}, \cite{Nee}, and \cite{V}. There are no
proofs given in these notes, though many of the results follow easily from earlier statements.  In addition, there are two more-recent books of Dimca \cite{Di} and Sch\"urmann \cite{Sch} which are very nice; the book of Dimca is relatively expository, while the book of Sch\"urmann is highly technical. 

While the results described here may seem very formal, in fact, the
treatment here is fairly informal.  If one wishes to avoid the formality of the derived category
altogether and, yet, still understand perverse sheaves, there are the  AMS notes of
MacPherson \cite{Mac2} which describe intersection homology and perverse sheaves
via Eilenberg-Steenrod type axioms, and MacPherson's more-recent treatment of perverse sheaves on regular cell complexes.

\vskip .1in

I have made some attempt to note results that conflict with the
statements that appear in other places.  I note these not to emphasize the
mistakes in those papers, but rather to let the reader know that I am aware of
the discrepancy and believe that the statement that I give is the correct one.  

Of course, this is not to say that there may not be mistakes in these notes -- in fact, I find
typographical mistakes constantly.  There may also be mistakes of others that I have copied, or
mistakes resulting from my own lack of understanding.  For all of these mistakes, I apologize.  

\vskip 0.2in There is another issue for which I definitely need to apologize. I initially learned about the derived category from Goresky and MacPherson's ``Intersection Homology II'', \cite{G-M3}. The convention of $j$ being a closed inclusion, while $i$ is an open inclusion -- which is reversed from what most people use -- is theirs. This is how I am used to thinking and writing, and it is far too late for me to get motivated to change these notes. I apologize for the confusion this may cause.

\vskip .2in

Many find it hard to believe that all the
machinery in these notes is necessary -- or even very useful -- for investigating problems in
the topology of singularities.  So, in the rest of the introduction, I give my
own initial motivation for learning this material.

\vskip .2in

Suppose that $f : (\mathbb C^{n+1}, \mathbf 0) \rightarrow (\mathbb C,0)$ is a
polynomial with a critical point at the origin.  We wish to discuss the case
where this critical point is non-isolated; so, let $\Sigma f$ denote the critical locus
of $f$ and let $s$ denote $\text{dim}_\mathbf 0\Sigma f$.  

The Milnor fibration for $f$ at the origin exists even for non-isolated
singularities.
 The Milnor fiber of $f$ at the origin has possibly non-trivial cohomology 
 only in dimensions between $n-s$ and $n$ (inclusive).
 
 If $s \geqslant 1$, the origin is not an isolated point in $\Sigma f$  and
so, at points $\mathbf p \in \Sigma f$ arbitrarily close to $\mathbf 0$, we may talk
about the Milnor fiber of $f$ {\bf at $\mathbf p$}.  Thus, we have a collection of
local data at each point of $\Sigma f$, and sheaf theory encodes how all this
local data patches together.

This is fine.  So why does one need a derived category of
complexes of sheaves, instead of just plain old normal everyday sheaves?  The
problem is:  at each point of $\Sigma f$, we wish to have the information about
the cohomology of the Milnor fiber at that point.  This means that {\bf after we
look at the stalk at a point $\mathbf p \in \Sigma f$}, we still wish to have
cohomology groups in all dimensions at our disposal.  It does not take long to
realize that what you need is a complex of sheaves.  But, really, one frequently
only cares about this complex of sheaves up to cohomology.  Very loosely
speaking, this is what the derived category gives you.

In this example, the two complexes of sheaves that one associates with the Milnor
fiber data correspond to the cohomology and reduced cohomology of the Milnor
fiber -- they are the complexes of sheaves of nearby and vanishing cycles,
respectively.  The earlier statement that the Milnor fiber  has possibly
non-trivial cohomology only in dimensions between $n-s$ and $n$ follows from certain 
vanishing conditions on these two complexes of sheaves, a vanishing conditions
which go by the name ``perverse".  The general results on perverse
sheaves are what give so much topological power to the machinery of the derived category.

\vskip .2in

Currently, this paper is organized as follows:

\vskip .1in

\S1. {\bf Constructible Complexes} -- This section contains general results on
bound\-ed, constructible complexes of sheaves and the derived category.

\vskip .1in

\S2. {\bf Perverse Sheaves} -- This section contains the definition and basic
results on perverse sheaves. Here, we also give the axiomatic characterization of
the intersection cohomology complex.  Finally in this section, we also give some
results on the category of perverse sheaves.  This categorical information is
augmented by that in section 5.

\vskip .1in

\S3. {\bf Nearby and Vanishing Cycles} -- In this section, we define and
examine the complexes of sheaves of nearby and vanishing cycles of an
analytic function.  These complexes contain hypercohomological information on the
Milnor fiber of the function under consideration.

\vskip .1in

\S4. {\bf Some Quick Applications} -- In this section, we give three easy
examples of results on Milnor fibers which follow from the machinery described in
the previous three sections.  

\vskip .1in

\S5. {\bf Truncation and Perverse Cohomology} -- This section contains an
informal discussion on $t$-structures. This enables us to describe truncation
functors and the perverse cohomology of a complex.  It also sheds some light on
our earlier discussion of the categorical structure of perverse sheaves.

\vskip .3in

\noindent\S1. {\bf Constructible Complexes}  

\vskip .1in

Much of this section is  lifted from Goresky and MacPherson's
``Intersection  Homology II" \cite{G-M3}.

\vskip .1 in

In these notes, we are primarily interested in sheaves on complex analytic
spaces, and we make an effort to state most results in this context.  However,
as one frequently wishes to do such things as intersect with a closed ball,
one really needs to consider at least the real semi-analytic case (that is,
spaces locally defined by finitely many real analytic inequalities).  In fact,
one can treat the subanalytic case.  Generally, when we leave the analytic
category we shall do so without comment, assuming the natural generalizations
of any needed results.  However, the precise statements in the subanalytic case
can be found in \cite{G-M2}, \cite{G-M3}, and \cite{K-S}.

\vskip .2in

Let $R$ be a regular Noetherian ring with finite Krull dimension (e.g., $\mathbb Z,
\mathbb Q, \text{or}\ \mathbb C$).  A complex $(\Adot, d^\bullet)$ (usually denoted simply by $\Adot$ if the differentials are
clear or arbitrary)
$$\begin{CD}
\cdots
\rightarrow \mathbf A^{-1} @>d^{-1}>> \mathbf A^0  @>d^{0}>> \mathbf A^1
@>d^{1}>> \mathbf A^2 @>d^{2}>> \cdots\end{CD}$$ of sheaves of $R$-modules on a
complex analytic space, $X$, is {\it bounded} if $\mathbf A^p = 0$ for $|p|$ large.

\smallskip

The cohomology sheaves $\mathbf H^p(\Adot)$ arise by taking the
(sheaf-theoretic) cohomology of the complex.  The stalk of $\mathbf H^p(\Adot)$
at a point $x$ is written $\mathbf H^p(\Adot)_x$ and is isomorphic to what one
gets by first taking stalks and then taking cohomology, i.e., $H^p(\Adot_x)$. The {\it support of $\Adot$}, $\operatorname{supp}\Adot$, is the closure of the set of points where $\Adot$ has non-zero stalk cohomology, i.e.,
$$
\operatorname{supp}\Adot \ = \ \overline{\left\{x\in X \ | \ H^*(\Adot)_x\neq 0\right\}}.
$$

\vskip .1in

The complex $\Adot$ is {\it constructible} with respect to a complex
analytic stratification, $\mathcal S = \{S_\alpha\}$, of $X$ provided that, for all
$\alpha$ and $i$, the cohomology sheaves $\mathbf H^i({\Adot})_{|_{S_\alpha}}$ are locally constant and have finitely-generated stalks;
we write  \hbox{$\Adot \in \mathbf D_{{}_{\mathcal S}}(X)$.}  If $\Adot \in \mathbf  D_{{}_{\mathcal S}}(X)$ and $\Adot$ is bounded, we
write $\Adot \in \mathbf D^b_{{}_{\mathcal S}}(X)$. 

If $\Adot \in \mathbf D^b_{{}_{\mathcal S}}(X)$ for some stratification
(and, hence, for any refinement of $\mathcal S$) we say that $\Adot$ is a
bounded, constructible complex and write $\Adot \in \mathbf D^b_c(X)$. 
(Note, however, that $\mathbf D^b_c(X)$ actually denotes the {\it derived} category and, while the
objects of this category are, in fact, the bounded, constructible complexes, the
morphisms are not merely maps between complexes.  We shall return to this.)

When it is important to indicate the base ring in the notation, we write
$\mathbf D_{{}_{\mathcal S}}(R_{{}_X})$, $\mathbf D^b_{{}_{\mathcal S}}(R_{{}_X})$, and $\mathbf D^b_c(R_{{}_X})$.

\vskip .1in

A single sheaf $\mathbf A$ on $X$ is considered a complex, $\Adot$, on $X$ by letting $\mathbf
A^0 = \mathbf A$ and $\mathbf A^i = 0$ for $i \neq 0$; thus, $\mathbf R^\bullet_{{}_X}$ 
denotes the constant sheaf on $X$.

\vskip .1in

The shifted complex $\Adot[n]$ is defined by $(\Adot[n])^k =
\mathbf A^{n + k}$ and differential $d^k_{[n]} = (-1)^nd^{k+n}$.

\vskip .1in

A map of complexes is a graded collection of sheaf maps $\phi^\bullet : \Adot \rightarrow \Bdot$ which commute
with the differentials. The shifted sheaf map $\phi^\bullet_{[n]}: \Adot[n] \rightarrow \Bdot[n]$ is defined by
$\phi^k_{[n]} := \phi^{k+n}$ (note the lack of a $(-1)^n$). A map of complexes is a {\it quasi-isomorphism} provided that
the induced maps
$$\mathbf H^p(\phi^\bullet) :
\mathbf H^p(\Adot)
\rightarrow  \mathbf H^p(\Bdot)$$ are isomorphisms for all $p$.  We use
the term ``quasi-isomorphic" to mean the equivalence relation generated by
``existence of a quasi-isomorphism"; this is sometimes referred to as
``generalized" quasi-isomorphic.

If  $\phi^\bullet : \Adot \rightarrow \mathbf I^\bullet$ is a quasi-isomorphism
and each $\mathbf I^p$ is injective, then $\mathbf I^\bullet$ is called an {\it
injective resolution} of $\Adot$.  Injective resolutions always exist
(in our setting), and are unique up to chain homotopy.  However, it is sometimes
important to associate one particular resolution to a complex, so it is
important that there is a {\it canonical injective resolution} which can be
associated to any complex (we shall not describe the canonical resolution here).

\vskip .1in

If $\Adot$ is a complex on $X$, then the {\it hypercohomology module},
$\mathbb H^p(X ; \Adot)$, is defined to be the $p$-th cohomology of the
global section functor applied to the canonical injective resolution of $\Adot$.

Note that if $\mathbf A$ is a single sheaf on $X$ and we form $\Adot$,
then $\mathbb H^p(X ; \Adot) = H^p(X ; \mathbf A) =$ ordinary sheaf
cohomology. In particular, $\mathbb H^p(X ; \mathbf R^\bullet_{{}_X}) = H^p(X ; R)$.  

\vskip .1in

Note also that if $\Adot$ and $\Bdot$ are quasi-isomorphic,
then  $\mathbb H^*(X ; \Adot) \cong  \mathbb H^*(X ; \Bdot)$.

\vskip .1in

If  $Y$ is a subspace of $X$ and $\Adot \in \mathbf D^b_c(X)$, then one
usually writes $\mathbb H^*(Y ; \Adot)$ in place of  $\mathbb H^*(Y ; {\Adot}_{|_Y})$.

\vskip .1in

The usual Mayer-Vietoris sequence is valid for hypercohomology; that is, if $U$
and $V$ form an open cover of $X$ and $\Adot \in \mathbf D^b_c(X)$, then
there is an exact sequence
$$
\dots\rightarrow \mathbb H^i(X; \Adot)
\rightarrow \mathbb H^i(U; \Adot)\oplus\mathbb H^i(V; \Adot)
\rightarrow \mathbb H^i(U\cap V; \Adot)
\rightarrow \mathbb H^{i+1}(X; \Adot)\rightarrow\dots .
$$

\vskip .1in

Of course, hypercohomology is not a homotopy invariant.  However, it is true
that:  if $\mathcal S$ is a real analytic Whitney stratification of $X$, $\Adot\in \mathbf D^b_{{}_\mathcal S}(X)$, and
$r: X\rightarrow [0,1)$ is a proper real analytic map such that, for all $S\in\mathcal
S$, $r_{|_S}$ has no critical values in $(0,1)$, then the inclusion $r^{-1}(0)\hookrightarrow X$ induces an isomorphism 
$$
\mathbb H^i(X; \Adot) \cong \mathbb H^i(r^{-1}(0); \Adot).
$$

\vskip .1in

If $R$ is an integral domain, we may talk about the rank of a finitely-generated
$R$-module.  In this case, if $\Adot \in \mathbf D^b_c(X)$, then the Euler
characteristic, $\chi$, of the stalk cohomology is defined as the
alternating sum of the ranks of the cohomology modules, i.e., $\chi(\Adot)_x =
\sum (-1)^i\text{ rank }\mathbf H^i(\Adot)_x$.  If the hypercohomology modules
are finitely-generated -- for instance, if $\Adot \in \mathbf D^b_c(X)$ and $X$
is compact -- then the Euler characteristic 
$\chi\bigl(\mathbb H^*(X ; \Adot)\bigr)$ is defined analogously.

\vskip .1in

If $\mathbf H^i(\Adot) = 0$ for all but, possibly, one value of $i$ - say, $i = p$, then
 $\Adot$ is quasi-isomorphic to the complex that has $\mathbf
H^p(\Adot)$ in degree $p$ and zero elsewhere.  We reserve the
term {\it local system} for a locally constant single sheaf or a complex
which is concentrated in degree zero and is locally constant.  If $M$ is
the stalk of a local system $\mathcal L$ on a path-connected space $X$, then $\mathcal L$ is determined  up to isomorphism by
a monodromy representation $\pi_1(X, \mathbf x) \rightarrow \operatorname{Aut}(M)$, where $x$ is a fixed point in $X$.

\vskip .1in

For any $\Adot \in \mathbf D^b_c(X)$, there is an $E_2$ cohomological spectral
sequence: $$E_2^{p,q} = H^p(X ; \mathbf H^q(\Adot)) \Rightarrow \mathbb
H^{{}^{p+q}}(X ; \Adot) .$$

\vskip .1in

If $\mathbf H^i(\Adot) = 0$ for all $i<d$, then $\mathbb
H^i(X ; \Adot)=0$ for all $i<d$ and $\mathbb H^d(X; \Adot)\cong H^0(X; \mathbf H^d(\Adot))$.

\vskip 0.1in

If $\Adot \in \mathbf D^b_c(X)$, $x \in X$, and $(X, x)$ is locally
embedded in some $\mathbb C^n$, then for all $\epsilon > 0$ small, the restriction
map $\mathbb H^q(B^\circ_\epsilon(x) ; \Adot) \rightarrow
\mathbf H^q(\Adot)_\mathbf x$ is an isomorphism $\bigl(\text{here},
B^\circ_\epsilon(x) = \bigl\{\mathbf z \in \bigl.\mathbb C^n \bigr|\ |
\mathbf z - \mathbf x | < \epsilon\bigr\}\ \bigr)$.  If, in addition, $R$ is a principal ideal
domain, the Euler characteristic $\chi\bigl(\mathbb H^*(B^\circ_\epsilon(x) - x ; \Adot)\bigr)$ is defined and 
$$\chi\bigl(\mathbb H^*(B^\circ_\epsilon(x) - x ; \Adot)\bigr) = \chi\bigl(\mathbb H^*(S_{\epsilon^\prime}(x) ;
\Adot)\bigr) = 0,$$
 where $0 < \epsilon^\prime < \epsilon$ and 
$S_{\epsilon^\prime}(x)$ denotes the sphere of radius $\epsilon^\prime$ centered at $x$.

\vskip .1in

We now wish to say a little about the morphisms in the derived category $\mathbf
D^b_c(X)$.  The derived category is obtained by formally inverting the
quasi-isomorphisms so that they become isomorphisms in $\mathbf D^b_c(X)$.  Thus,
$\Adot$ and $\Bdot$ are isomorphic in $\mathbf D^b_c(X)$
provided that there exists a complex $\mathbf C^\bullet$ and quasi-isomorphisms
$\Adot \leftarrow \mathbf C^\bullet \rightarrow \Bdot$;  $\Adot$ and $\Bdot$ are then said to be {\it incarnations} of the
same isomorphism class in $\mathbf D^b_c(X)$.  

More generally, a morphism in
$\mathbf D^b_c(X)$ from $\Adot$ to $\Bdot$ is an equivalence
class of diagrams of maps of complexes $\Adot \leftarrow \mathbf C^\bullet
\rightarrow \Bdot$ where $\Adot \leftarrow \mathbf C^\bullet$
is a quasi-isomorphism.  Two such diagrams,
$$\begin{CD}
\Adot @<f_1<< \mathbf C^\bullet_1
@>g_1>> \Bdot , \hskip .4in \Adot @<f_2<< \mathbf C^\bullet_2
@>g_2>> \Bdot
\end{CD}
$$
are equivalent provided that there exists a third such diagram 
$$\begin{CD}\Adot @<f<< \mathbf C^\bullet
@>g>> \Bdot\end{CD}$$ and a diagram

\vbox{$$\mathbf C^\bullet_1$$
$${}^{f_1}\swarrow\hskip .2in\uparrow\hskip .2in\searrow {}^{g_1}$$

$$\begin{CD}\Adot @<\ \ f\ << \mathbf C^\bullet @>\ g\ \ >> \Bdot\end{CD}$$

$${}_{f_2}\nwarrow\hskip .2in\downarrow\hskip .2in\nearrow {}_{g_2}$$
$$\mathbf C^\bullet_2$$}

\noindent which commutes up to (chain) homotopy.

Composition of morphisms in $\mathbf D^b_c(X)$ is not difficult to describe.  If we have two
representatives of morphisms, from $\Adot$ to $\Bdot$ and from 
$\Bdot$ to $\mathbf D^\bullet$, respectively,
$$\begin{CD}
\Adot @<f_1<< \mathbf C^\bullet_1
@>g_1>> \Bdot , \hskip .4in \Bdot @<f_2<< \mathbf C^\bullet_2
@>g_2>> \mathbf D^\bullet\end{CD}
$$
then we consider the pull-back $\mathbf C^\bullet_1\times_{{}_{\Bdot}}\mathbf C^\bullet_2$ (in the category of chain
complexes) and the projections $\pi_1$ and $\pi_2$ to $\mathbf C^\bullet_1$ and $\mathbf C^\bullet_2$,
respectively.  As $f_2$ is a quasi-isomorphism, so is $\pi_1$, and the composed morphism from $\Adot$ to $\mathbf D^\bullet$ is represented by 
$$\begin{CD}\Adot @<f_1\circ\pi_1<< \mathbf C^\bullet_1\times_{{}_{\Bdot}}\mathbf C^\bullet_2
 @>g_2\circ\pi_2>> \mathbf D^\bullet\end{CD}.$$

\vskip .1in

If we restrict ourselves to considering only injective complexes, by
associating to any complex its canonical injective resolution, then morphisms in
the derived category become easy to describe -- they are chain-homotopy classes
of maps between the injective complexes.  

\vskip .1in

The moral is: in $\mathbf D^b_c(X)$, we essentially only care about complexes up to
quasi-isomorphism.  Note, however, that the objects of $\mathbf D^b_c(X)$ are {\bf
not} equivalence classes -- this is one reason why it is important that to each
complex we can associate a {\bf canonical} injective resolution.  It allows us
to talk about certain functors in $\mathbf D^b_c(X)$ being {\bf naturally}
isomorphic.  When we write $\Adot \cong \Bdot$, we mean in
$\mathbf D^b_c(X)$.  As we shall discuss later, $\mathbf D^b_c(X)$ is an additive category,
but is not Abelian.

\vskip .1in

\noindent{\bf Warning}:  While morphisms of complexes which induce isomorphisms on cohomology sheaves become isomorphisms in the
derived category, there are morphisms of complexes which induce the zero map on cohomology sheaves but are not zero in the 
derived category.  The easiest example of such a morphism is given by the following.

Let $X$ be a space consisting of two complex lines $L_1$ and $L_2$ which intersect in a single point $\mathbf p$. For $i=1,2$,
let $\widetilde{\mathbb C}_{L_i}$ denote the $\mathbb C$-constant sheaf on $L_i$ extended by zero to all of $X$. There is a canonical
map,
$\alpha$, from the sheaf
$\mathbb C_X$ to the direct sum of sheaves
$\widetilde{\mathbb C}_{L_1}\oplus\widetilde{\mathbb C}_{L_2}$, which on
$L_1-\mathbf p$ is
$\operatorname{id}\oplus 0$, on $L_2-\mathbf p$ is $0\oplus\operatorname{id}$, and is the diagonal map on the stalk at $\mathbf p$.
Consider the complex, $\Adot$, which has  $\mathbb C_X$ in degree $0$, $\widetilde{\mathbb C}_{L_1}\oplus\widetilde{\mathbb
C}_{L_2}$ in degree $1$, zeroes elsewhere, and the coboundary map from degree $0$ to degree $1$ is $\alpha$.  This complex has
cohomology only in degree
$1$.  Nonetheless, the morphism of complexes from $\Adot$ to
$\mathbb C^\bullet_X$ which is the identity in degree $0$ and is zero elsewhere determines a non-zero morphism in the derived
category.

\vskip .3in

We now wish to describe {\it derived functors};  for this, we will need the derived category of an
arbitrary Abelian category $\mathcal C$.

\vskip .1in

Let $\mathcal C$ be an Abelian category.  Then, the derived category of bounded complexes in $\mathcal C$
is the category whose objects consist of bounded differential complexes of objects of $\mathcal C$, and
where the morphisms are obtained exactly as in the case of $\mathbf D^b_c(X)$ -- namely, by inverting
the quasi-isomorphisms as we did above.  Naturally, we denote this derived category by $\mathbf
D^b(\mathcal C)$.

We need some more general notions before we come back to complexes of sheaves.  If $\mathcal C$ is an
Abelian category, then we let $\mathbf K^b(\mathcal C)$ denote the category whose objects are again 
bounded differential complexes of objects of $\mathcal C$, but where the morphisms are chain-homotopy
classes of maps of differential complexes.  A {\it triangle} in $\mathbf K^b(\mathcal C)$ is a sequence
of morphisms \hbox{$\Adot\rightarrow\Bdot\rightarrow\mathbf C^\bullet\rightarrow\Adot[1]$}, which is usually written in the
more``triangular'' form

\vskip .1in

\vbox{
$\hskip 2.5in \Adot \longrightarrow \ \Bdot$ \vskip .03in
$\hskip 2.35in \ _{[1]} \ \nwarrow\hskip .15in \swarrow$ \vskip .03in $\hskip
2.75in \mathbf C^\bullet$
}

\vskip .1in

A triangle in $\mathbf K^b(\mathcal C)$ is called {\it distinguished} if it
is isomorphic in $\mathbf K^b(\mathcal C)$ to a diagram of maps of complexes

\vskip .2in

\vbox{
$\hskip 2.5in \widetilde{\mathbf A}^\bullet \xrightarrow{\ \phi\ }\ \widetilde{\mathbf B}^\bullet$ \vskip .03in $\hskip 2.35in
\ _{[1]}
\
\nwarrow\hskip .15in \swarrow$
\vskip .03in $\hskip 2.75in \mathbf M^\bullet$
}

\vskip .1in

\noindent where $\mathbf M^\bullet$ is the algebraic mapping cone of $\phi$ and
$\widetilde{\mathbf B}^\bullet \rightarrow \mathbf M^\bullet \rightarrow \widetilde{\mathbf A}^\bullet[1]$ are
the canonical maps.  (Recall that the algebraic mapping cone is defined by
$$\mathbf M^k := \widetilde{\mathbf A}^{k+1} \oplus\ \widetilde{\mathbf B}^k \longrightarrow 
\widetilde{\mathbf A}^{k+2} \oplus\ \widetilde{\mathbf B}^{k+1} =: \mathbf M^{k+1}$$ $$(a,b)
\longmapsto (-\partial a, \phi a + \delta b)$$ where $\partial$ and $\delta$ are
the differentials of $\widetilde{\mathbf A}^\bullet$ and  $\widetilde{\mathbf B}^\bullet$
respectively.) Note that if $\phi=0$, then we have an equality $\Mdot= \Adot[1]\oplus\Bdot$ (recall that the shifted complex
$\Adot[1]$ has as its differential the negated, shifted differential of $\Adot$).

\vskip .1in

Now we can define derived functors. Let $\mathcal C$ denote the Abelian category of sheaves of
$R$-modules on an analytic space $X$, and let $\mathcal C^\prime$ be another Abelian category.  Suppose
that $F$ is an additive, covariant functor from $\mathbf K^b(\mathcal C)$ to $\mathbf K^b(\mathcal C^\prime)$
such that $F\circ [1] =[1]\circ F$ and such that $F$ takes distinguished triangles to distinguished
triangles (such an $F$ is called a {\it functor of triangulated categories}). Suppose also that, for
all complexes of injective sheaves $\mathbf I^\bullet\in \mathbf K^b(\mathcal C)$ which are
quasi-isomorphic to $0$, $F(\mathbf I^\bullet)$ is also quasi-isomorphic to $0$.

Then, $F$ induces a morphism $RF$ -- the {\it right derived functor of $F$} -- from $\mathbf D^b(X)$ to
$\mathbf D^b(\mathcal C^\prime)$; for any $\Adot\in \mathbf D^b(X)$, let $\Adot\rightarrow\mathbf I^\bullet$ denote the canonical injective resolution of $\Adot$, and define $RF(\Adot) := F(\mathbf I^\bullet)$.  The action of $RF$ on the
morphisms is the obvious associated one.

\vskip .1in

A morphism $F: \mathbf K^b(\mathcal C)\rightarrow\mathbf K^b(\mathcal C^\prime)$ as described above is
frequently obtained by starting with a left-exact functor $T: \mathcal C\rightarrow\mathcal C^\prime$ and
then extending $T$ in a term-wise fashion to be a functor from 
$\mathbf K^b(\mathcal C)$ to $\mathbf K^b(\mathcal C^\prime)$.  In this case, we naturally write $RT$ for the
derived functor.

This is the process which is applied to:

\vskip .1in

\noindent $\Gamma(X; \cdot)$ (global sections);

\vskip .1in

\noindent $\Gamma_c(X; \cdot)$ (global sections with compact support);

\vskip .1in

\noindent $f_*$ (direct image);

\vskip .1in

\noindent  $f_!$ (direct image with proper supports); and

\vskip .1in

\noindent $f^*$ (pull-back or inverse image),

\vskip .1in

\noindent where $f : X \rightarrow Y$ is a continuous map (actually, in these
notes, we would need an analytically constructible map; e.g., an analytic map).

\vskip.1in

If the functor $T$ is an exact functor from sheaves to sheaves, then $RT(\Adot) \cong T(\Adot)$; in this
case, we normally suppress the $R$. Hence, if $f : X \rightarrow Y$, $\Adot\in \mathbf D^b(X)$, and $\Bdot\in
\mathbf D^b(Y)$, we write:

\hbox{}\hskip .25in $f^*\Bdot$; 

\hbox{}\hskip .25in $f_!\Adot$, if $f$ is the inclusion of a locally closed subset and,
hence, $f_!$ is extension by zero;

\hbox{}\hskip .25in $f_*\Adot$, if $f$ is the inclusion of a closed
subspace.

\vskip .2in

Note that hypercohomology is just the cohomology of the derived global section functor, i.e., $\mathbb H^*(X; \cdot) = H^*\circ
R\Gamma(X; \cdot)$.  The cohomology of the derived functor of global sections with compact support is the {\it compactly
supported hypercohomology} and is denoted $\mathbb H^*_c(X ; \Adot)$.

\vskip .2in

If $f:X\rightarrow Y$ is the inclusion of a subset and $\Bdot\in \mathbf D^b(Y)$, then the {\it restriction of
$\Bdot$} to $X$ is defined to be
$f^*(\Bdot)$, and is usually denoted by ${\Bdot}_{|_{X}}$.

\vskip .4in

It is important to note that, in general, for an inclusion map $f:X\rightarrow Y$, $Rf_*$ and $Rf_!$ need not take constructible complexes to constructible complexes. However, in the algebraic setting, the derived functors of inclusions of Zariski open or closed subsets do take constructible complexes to constructible complexes. In the analytic setting, the derived functors of inclusions of Zariski closed subsets take constructible complexes to constructible complexes, but the derived functors of inclusions of Zariski open subsets need not take constructible complexes to constructible complexes. However, suppose $\Adot\in \mathbf D^b_c(X)$, $Z$ is analytically Zariski closed in $X$, and we have the inclusion $i:X-Z\rightarrow X$; then, $Ri_*i^*\Adot$ and $Ri_!i^*\Adot$ are in $\mathbf D^b_c(X)$.

\bigskip

If $f : X \rightarrow Y$ is continuous and $\Adot \in \mathbf D^b_c(X)$, there is a
natural map 
$$Rf_!\Adot \rightarrow Rf_*\Adot .$$

\vskip .1in

For $f : X \rightarrow Y$ continuous, there are canonical isomorphisms
$$
R\Gamma(X; \Adot)\cong R\Gamma(Y; Rf_*\Adot)\text{\ and\ }R\Gamma_c(X; \Adot)\cong R\Gamma_c(Y; Rf_!\Adot)
$$
which lead to canonical isomorphisms
$$\mathbb H^*(X ; \Adot) \cong \mathbb H^*(Y ; Rf_*\Adot)
 \text{ \ and \ }
 \mathbb H^*_c(X ; \Adot) \cong  \mathbb H^*_c(Y ; Rf_!\Adot)$$  
 for all $\Adot$ in $\mathbf D^b_c(X)$.

\vskip .2in

If $f : X \rightarrow Y$ is continuous, $\Adot \in \mathbf D^b_c(X)$, and $\Bdot \in \mathbf D^b_c(Y)$, there are natural maps induced by  restriction of sections 
$$\Bdot\rightarrow Rf_*f^*\Bdot \text{ \ and \ } 
f^*Rf_*\Adot \rightarrow \Adot.$$
If $f$ is the inclusion of a locally closed subset, then $f^*Rf_*\Adot \rightarrow \Adot$ is an isomorphism; more precisely, there is a natural isomorphism of functors $f^*Rf_*\cong\operatorname{id}$.

\vskip .2in

If $\{S_\alpha\}$ is a stratification of $X$, $\Adot \in  \mathbf
D^b_{{}_{\{S_\alpha \times \mathbb C^k\}}}(X \times \mathbb C^k)$, and $\pi : X \times
\mathbb C^k \rightarrow X$ is the projection, then restriction of sections induces
a quasi-isomorphism $\pi^*R\pi_*\Adot \rightarrow \Adot$.

It follows easily that if $j : X \hookrightarrow X \times \mathbb C^k$ is the zero section,
then \hbox{$\pi^*j^*\Adot \cong \Adot$.}  This says exactly what one
expects: the complex $\Adot$ has a product  structure in the $\mathbb C^k$
directions.

An important consequence of this is the following:  let $\mathcal S = \{S_\alpha\}$ be a
Whitney stratification of $X$ and let $\Adot \in  \mathbf
D^b_{{}_{\mathcal S}}(X)$.  Let $x \in S_\alpha \subseteq X$.  As $S_\alpha$ is
a Whitney stratum, $X$ has a product structure along $S_\alpha$ near $x$.  By the
above, $\Adot$ itself also has a product structure along $S_\alpha$. 
Hence, by taking a normal slice, many problems concerning the complex $\Adot$ can be reduced to considering a zero-dimensional stratum.

\vskip .2in

Let $\Adot, \Bdot\in \mathbf D^b_c(X)$.  Define $\Adot\otimes\Bdot$ to be the single complex which is associated to the double complex $\mathbf A^p\otimes\mathbf
B^q$.  The left derived functor $\Adot\otimes^L*$ is defined by 
$$\Adot\otimes^L\Bdot = \Adot\otimes\mathbf J^\bullet ,$$
where $\mathbf J^\bullet$ is a flat resolution of $\Bdot$, i.e., the stalks of 
$\mathbf J^\bullet$ are flat $R$-modules and there exists a
quasi-isomorphism $\mathbf J^\bullet\rightarrow\Bdot$.

For all $\Adot, \Bdot\in \mathbf D^b_c(X)$, there is an isomorphism 
$\Adot\otimes^L\Bdot \cong\Bdot\otimes^L\Adot$.

\vskip .1in

For any map $f:X\rightarrow Y$ and any $\Adot, \Bdot\in \mathbf D^b_c(Y)$,
$$f^*(\Adot\otimes^L\Bdot)\cong f^*\Adot\otimes^L f^*\Bdot .$$

\vskip .2in

Fix a complex $\Bdot$ on $X$. There are two covariant functors which we wish to
consider: the functor $\mathbf {Hom}^\bullet(\Bdot , *)$ from the category of
complexes of sheaves on $X$ to complexes of sheaves on $X$ and the functor 
$Hom^\bullet(\Bdot , *)$ from the category of complexes of sheaves on $X$ to the
category of complexes of $R$-modules.  These functors are given by 
$$\left(\mathbf
{Hom}^\bullet(\Bdot, \Adot)\right)^n = \prod_{p \in \mathbb Z} \mathbf
{Hom}(\mathbf B^p, \mathbf A^{n+p})$$
and 
$$\left(Hom^\bullet(\Bdot, \Adot)\right)^n = \prod_{p \in \mathbb Z}
Hom(\mathbf B^p, \mathbf A^{n+p})$$ 
with differential given by
$$
[\partial^nf]^p = \partial^{n+p}f^p + (-1)^{n+1}f^{p+1}\partial^p
$$ 
(there is an indexing error in \cite{I}, 12.4]).  The associated derived functors are, respectively,  $R\mathbf
{Hom}^\bullet(\Bdot , *)$ and $RHom^\bullet(\Bdot, *)$.  

\smallskip

Note that if $\mathbf R^\bullet_{{}_X}$ is the constant sheaf on $X$, then there is a natural isomorphism of functors $R\mathbf
{Hom}^\bullet(\mathbf R^\bullet_{{}_X} , *)\cong \operatorname{id}$.

\smallskip

For all $k$, 
$$R\mathbf {Hom}^\bullet(\Bdot, \Adot[k]) = R\mathbf {Hom}^\bullet(\Bdot[-k], \Adot)=\big(R\mathbf {Hom}^\bullet(\Bdot, \Adot)\big)[k].$$

\smallskip

If $\Pdot\rightarrow\Bdot$ is a projective resolution of $\Bdot$, then, in $\mathbf D^b_c(X)$, $R\mathbf{Hom}^\bullet(\Bdot ,
\Adot)$ is isomorphic to
$\mathbf{Hom}^\bullet(\Pdot, \Adot)$. 

The functor $RHom^\bullet(\Bdot, *)$ is naturally isomorphic to the derived global
sections functor applied to $R\mathbf {Hom}^\bullet(\Bdot, *)$, i.e., for any 
$\Adot \in  
\mathbf D^b_c(X)$,
$$
RHom^\bullet(\Bdot, \Adot)\cong 
R\Gamma\left(X;R\mathbf{Hom}^\bullet(\Bdot, \Adot)\right) .
$$
$H^0(RHom^\bullet(\Bdot, \Adot))$ is naturally isomorphic as an
$R$-module to the derived
category homomorphisms from $\Bdot$ to $\Adot$, i.e., $$
H^0(RHom^\bullet(\Bdot, \Adot))\cong 
Hom_{{}_{\mathbf D^b_c(X)}}(\Bdot, \Adot) .
$$
If $\Bdot$ and $\Adot$  have locally constant cohomology sheaves on $X$ 
then, for all $x\in X$, $R\mathbf {Hom}^\bullet(\Bdot, \Adot)_x$ is 
naturally isomorphic to
$RHom^\bullet(\Bdot_x, \Adot_x)$.

\vskip .1in

For all $\Adot, \Bdot, \mathbf C^\bullet\in \mathbf D^b_c(X)$, there is a natural 
isomorphism
$$
R\mathbf{Hom}^\bullet(\Adot\otimes^L\Bdot, \mathbf C^\bullet) 
\cong R\mathbf{Hom}^\bullet(\Adot, R\mathbf{Hom}^\bullet(\Bdot, \mathbf C^\bullet)).
$$
Moreover, if $\mathbf C^\bullet$ has locally constant cohomology sheaves, then there is an isomorphism 
$$
R\mathbf{Hom}^\bullet(\Adot, \Bdot\otimes^L\mathbf C^\bullet) 
\cong R\mathbf{Hom}^\bullet(\Adot, \Bdot)\otimes^L\mathbf C^\bullet .
$$

\vskip .2in

For all $j$, we define $\mathbf{Ext}^j(\Bdot, \Adot) := 
\mathbf H^j(R\mathbf {Hom}^\bullet(\Bdot, \Adot))$  and  
define $$Ext^j(\Bdot, \Adot) := H^j(RHom^\bullet(\Bdot, \Adot)).$$  It is immediate that we have 
isomorphisms of $R$-modules 
$$
Ext^j(\Bdot, \Adot) = 
H^0(RHom^\bullet(\Bdot, \Adot[j])) \cong 
Hom_{{}_{\mathbf D^b_c(X)}}(\Bdot, \Adot[j]) .
$$

\vskip .1in

If $X = point$, $\Bdot\in \mathbf D^b_c(X)$, and the base ring is a PID, then  
$\Bdot \cong \bigoplus_k\mathbf H^k(\Bdot)[-k]$ in 
$\mathbf D^b_c(X)$; if we also have $\Adot\in \mathbf D^b_c(X)$, then 
$$
\mathbf H^i(\Adot\otimes^L\Bdot)\cong\Big(\bigoplus_{p+q=i}\mathbf H^p(\Adot)\otimes\mathbf H^q(\Bdot)\Big)\ \oplus\ 
\big(\bigoplus_{r+s=i+1}\operatorname{Tor}(\mathbf H^r(\Adot), \mathbf H^s(\Bdot))\big).
$$
If, in addition, the cohomology modules of 
$\Adot$ are projective (hence, free), then
$$Hom_{{}_{\mathbf D^b_c(X)}}(\Adot, \Bdot) 
\cong \bigoplus_k Hom\left(\mathbf H^k(\Adot), \mathbf H^k(\Bdot)\right) .$$

\vskip .1in

If we have a map $f:X\rightarrow Y$, then the functors $f^*$ and $Rf_*$ are adjoints of
each other in the derived category of bounded complexes.  In fact, for all $\Adot$ on $X$ and $\Bdot$ on $Y$, there is a canonical isomorphism in $\mathbf D^b(Y)$
$$
R\mathbf{Hom}^\bullet(\Bdot, Rf_*\Adot)\cong 
Rf_*R\mathbf{Hom}^\bullet(f^*\Bdot, \Adot)
$$
and so
$$
Hom_{{}_{\mathbf D^b(Y)}}(\Bdot, Rf_*\Adot) \cong H^0\left(RHom^\bullet(\Bdot, Rf_*\Adot)\right)\cong \mathbb H^0\left(Y; R\mathbf{Hom}^\bullet(\Bdot, Rf_*\Adot)\right)\cong
$$
$$
\mathbb H^0\left(X; R\mathbf{Hom}^\bullet(f^*\Bdot, \Adot)\right)
\cong  H^0\left(RHom^\bullet(f^*\Bdot, \Adot)\right)\cong Hom_{{}_{\mathbf D^b(X)}}(f^*\Bdot,
\Adot).
$$

\vskip .3in

We wish now to describe an analogous adjoint for $Rf_!$

\vskip .1in

Let $\mathbf I^\bullet$ be a complex of injective sheaves on $Y$.  Then, $f^!(\mathbf
I^\bullet)$ is defined to be the sheaf associated to the presheaf given by
$$
\Gamma(U; f^!\mathbf I^\bullet) = Hom^\bullet(f_!\mathbf K^\bullet_{{}_U}, \mathbf I^\bullet) ,
$$
for any open $U \subseteq X$, where $\mathbf K^\bullet_{{}_U}$ denotes the canonical
injective resolution of the constant sheaf $\mathbf R^\bullet_{{}_U}$.  For any $\Adot\in\mathbf D^b(X)$, define $f^!\Adot$ to be $f^!\mathbf I^\bullet$, where
$\mathbf I^\bullet$ is the canonical injective resolution of $\Adot$.  

Now that we
have this definition, we may state:

\vskip .2in

\noindent ({\bf Poincar\'e-Verdier Duality}) \hskip .2in If $f: X \rightarrow Y$ is continuous, $\Adot \in \mathbf
D^b(X)$, and $\Bdot \in \mathbf D^b(Y)$, then there is a
canonical isomorphism in $\mathbf D^b(Y)$: $$Rf_*R\mathbf {Hom}^\bullet(\Adot , f^!\Bdot) \cong R\mathbf {Hom}^\bullet(Rf_!\Adot ,
\Bdot)$$
and so
$$
Hom_{{}_{\mathbf D^b(X)}}(\Adot, f^!\Bdot) \cong  Hom_{{}_{\mathbf D^b(Y)}}(Rf_!\Adot,
\Bdot).
$$

\vskip .2in

In particular, replacing $\Adot$ with the constant sheaf $\mathbf R^\bullet_{{}_X}$, we obtain a natural isomorphism
$$
Rf_*f^!\Bdot\cong R\mathbf {Hom}^\bullet(Rf_!\mathbf R^\bullet_{{}_X} , \Bdot).
$$

It follows that, if $f$ is the inclusion of a locally closed subset, then
$$
f^!\Bdot\cong f^*R\mathbf {Hom}^\bullet(Rf_!\mathbf R^\bullet_{{}_X} , \Bdot).
$$

\medskip

If $\Bdot$ and $\mathbf C^\bullet$ are in $\mathbf D^b(Y)$ and $f: X \rightarrow Y$ is continuous, then we have an
isomorphism $$
f^!R\mathbf{Hom}^\bullet(\Bdot, \mathbf C^\bullet)\cong 
R\mathbf{Hom}^\bullet(f^*\Bdot, f^!\mathbf C^\bullet).$$

\bigskip

Let $f: X\rightarrow point$.  Then, the {\it dualizing complex}, $\mathbb D^\bullet_{{}_X}$,
is $f^!$ applied to the constant sheaf, i.e., $\mathbb D^\bullet_{{}_X} = f^!\mathbf
R^\bullet_{{}_{pt}}$.  For any complex $\Adot\in \mathbf D^b(X)$, the {\it
Verdier dual} (or, simply, {\it the dual}) of $\Adot$ is 
$R\mathbf{Hom}^\bullet(\Adot, \mathbb D^\bullet_{{}_X})$  and is denoted by 
$\mathcal D_{{}_X}\Adot$ (or just $\mathcal D\Adot$).  There is a canonical
isomorphism between $\mathbb D^\bullet_{{}_X}$ and the dual of the constant sheaf on $X$,
i.e., $\mathbb D^\bullet_{{}_X}\cong \mathcal D\mathbf R^\bullet_{{}_X}$. The Verdier dual yields a contravariant functor from $\mathbf D^b(X)$ to itself, and there is a natural morphism from $\Adot$ to $\mathcal D(\mathcal D \Adot)$, i.e., a natural transformation from $\operatorname{id}_{\mathbf D^b(X)}$ to $\mathcal D\mathcal D$.

Suppose that the base ring $R$ is a field or a Dedekind domain (e.g., a PID). Then:

1)  $\Adot\in \mathbf D^b(X)$ is constructible if and only if $\mathcal D\Adot$ is constructible. Thus $\mathcal D$ yields a contravariant functor from $\mathbf D^b_c(X)$ to itself.

\smallskip

2) On $\mathbf D^b_c(X)$, the natural transformation from $\operatorname{id}_{\mathbf D^b(X)}$ to $\mathcal D\mathcal D$ is a natural isomorphism.

\vskip .2in

Let $\Adot \in \mathbf D^b_c(X)$.  If the base ring $R$ is a field or a Dedekind domain (e.g., a PID), then $\mathcal D\Adot$
is well-defined up to quasi-isomorphism by:

\vskip .1in

\noindent for any open $U \subseteq X$, there is a natural split exact sequence:
$$0 \rightarrow Ext(\mathbb H^{q+1}_c(U ; \Adot) , R) \rightarrow \mathbb
H^{-q}(U ; \mathcal D\Adot) \rightarrow Hom(\mathbb H_c^q(U ; \Adot), R) \rightarrow 0  .$$

\vskip .2in

In particular, if $R$ is a field, then $\mathbb H^{-q}(U ; \mathcal D\Adot)
\cong 
 \mathbb H_c^q(U ; \Adot)$, and so 
 $$
 \mathbf H^q(\mathcal D\Adot)_x
\cong \mathbb H^q(B^\circ_\epsilon(x) ; \mathcal D\Adot) \cong 
\mathbb H^{-q}_c(B^\circ_\epsilon(x) ; \Adot).
$$
If, in
addition, $X$ is compact, $\mathbb H^{-q}(X ; \mathcal D\Adot) \cong 
 \mathbb H^q(X ; \Adot)$.
 
\vskip .2in

Dualizing is a local operation, i.e., if $i : U \hookrightarrow X$ is the
inclusion of an open subset and $\Adot \in \mathbf D^b_c(X)$, then
$i^*\mathcal D\Adot \cong \mathcal Di^*\Adot$.

  If $\mathcal L$ is a local system on a connected
real $m$-manifold, $N$, then $(\mathcal D\mathcal L^\bullet)[-m]$ is a local system; if, in addition, $N$ is  oriented, and $\mathcal L$ is actually
locally free with stalks $R^a$ and monodromy representation $\eta : \pi_1(N, \mathbf p) \rightarrow
\operatorname{Aut}(R^a)$, then $\mathcal D\mathcal L^\bullet[-m]$ is a local
system with stalks equal to the dual, $(R^a)^\vee:=Hom_R(R^a, R)\cong R^a$, of $R^a$ and monodromy \hbox{${}^\vee\eta : \pi_1(N, \mathbf p) \rightarrow
\operatorname{Aut}((R^a)^\vee)$}, where ${}^\vee\eta(\alpha)$ is the dual of $(\eta(\alpha))^{-1}$, i.e., $\big({}^\vee\eta(\alpha)\big)(f):=  f\circ(\eta(\alpha))^{-1} = f\circ(\eta(\alpha^{-1})) $.

 \vskip .2in
 
 If $\Adot \in \mathbf D^b_c(X)$, then $\mathcal D(\Adot[n]) = (\mathcal
D\Adot)[-n]$.

\vskip .2in

If $\pi : X \times \mathbb C^n \rightarrow X$ is projection, then
$\mathcal D(\pi^*\Adot)[-n] \cong \pi^*(\mathcal D\Adot)[n]$.

\vskip .2in

The dualizing complex, $\mathbb D_{{}_X}$, is quasi-isomorphic to the complex of sheaves of
singular chains on $X$ which is associated to the complex of
presheaves, $\mathbf C^\bullet$, given by $\Gamma(U ; \mathbf C^{-p}) := C_p(X , X - U ; R)$.

The cohomology sheaves of $\mathbb D^\bullet_{{}_X}$ are non-zero in negative
degrees only, with stalks $\mathbf H^{-p}(\mathbb D^\bullet_{{}_X})_x = H_p(X , X - x
; R)$.

If $X$ is an oriented, real $m$-manifold, then $\mathbb D^\bullet_{{}_X}[-m]$ is 
quasi-isomorphic to $\mathbf R^\bullet_{{}_X}$, and so  $\mathcal D\mathbf R^\bullet_{{}_X}\cong \mathbf R^\bullet_{{}_X}[m]$ and, for all $\Adot \in \mathbf D^b_c(X)$, 
$$\mathcal D\Adot \cong
R\mathbf{Hom}^\bullet(\Adot , \mathbf R^\bullet_{{}_X}[m]) = 
\left(R\mathbf{Hom}^\bullet(\Adot , \mathbf
R^\bullet_{{}_X})\right)[m] .$$
Note that if $X$ is an even-dimensional, oriented, real $m$-manifold, with $m=2n$, then $\mathcal D(\mathbf R^\bullet_{{}_X}[n])\cong \mathbf R^\bullet_{{}_X}[n]$, i.e., $\mathbf R^\bullet_{{}_X}[n]$ is self-dual.

\vskip .1in

$\mathbb D^\bullet_{{}_{V\times W}}$ is naturally isomorphic to 
$\pi_1^*\mathbb D^\bullet_{{}_V}\otimes^L\pi_2^*\mathbb D^\bullet_{{}_W}$, where $\pi_1$ and $\pi_2$ are
the projections onto $V$ and $W$, respectively.

\vskip .1in

Hypercohomology with coefficients in $\mathbb D^\bullet_{{}_X}$ is isomorphic to Borel-Moore homology, with a negation in the degree, i.e., $\mathbb H^{-i}(X ; \mathbb D^\bullet_{{}_X}) \cong\mathbb H_i^{BM}(X)$.

\vskip .1in

If $X$ is a real, smooth, oriented $m$-manifold and $R = \mathbb R$, then $\mathbb
D^\bullet_{{}_X}[-m]$ is naturally isomorphic to the complex of real
differential forms on $X$.

\vskip .1in

$\mathbb D^\bullet_{{}_X}$ is constructible with respect to any Whitney
stratification of $X$.  It follows that if $\mathcal S$ is a Whitney
stratification of $X$, then $\Adot \in
\mathbf D^b_{{}_{\mathcal S}}(X)$ if and only if $\mathcal D\Adot \in
\mathbf D^b_{{}_{\mathcal S}}(X)$.

\vskip .2in

 For all $\Adot, \Bdot\in\mathbf D^b_c(X)$, we have isomorphisms
$$
R\mathbf{Hom}^\bullet(\Adot, \Bdot) \cong
R\mathbf{Hom}^\bullet(\mathcal D\Bdot, \mathcal D\Adot) \cong
\mathcal D\left(\mathcal D\Bdot\otimes^L\Adot\right) .
$$

\vskip .2in

If $f: X \rightarrow Y$ is continuous, then we have natural isomorphisms
$$Rf_! \cong  \mathcal DRf_*\mathcal D \ \ \text{and}\ \ 
f^! \cong \mathcal Df^*\mathcal D .$$

\bigskip

Suppose that $k:Y\hookrightarrow X$ is the inclusion of a locally closed subset, and let 
$$
\operatorname{id}\xrightarrow{ \ \alpha \ } Rk_*k^*\hskip 0.2in \textnormal{ and }\hskip 0.2in   k_!k^!\xrightarrow{ \ \beta \ } \operatorname{id}
$$
be the standard natural transformations. Then there is a natural isomorphism 
$$\vdual \big(Rk_*k^*\big)\xrightarrow[\cong]{\gamma} k_!k^!\vdual$$ 
such that, for all $\Adot\in \mathbf D^b_c(X)$, $\vdual(\alpha_{\Adot})=\beta_{\vdual\Adot}\circ \gamma_{\Adot}$.

\vskip 0.2in

If $Y\subseteq X$ and $f : X -Y \hookrightarrow X$ is the inclusion, we define
$$\mathbb H^k(X , Y ; \Adot) := \mathbb H^k(X ; f_!f^!\Adot) .$$

\vskip .2in

Suppose that $X$ is compact, $U$ is an open subset of $X$, and $Y:=X-U$. Then,
$$
\mathbb H^k_c(U; \Adot)\cong \mathbb H^k(X, Y; \Adot).
$$

\vskip .2in

Let $x \in X$ and $\Adot \in \mathbf D^b_c(X)$. Let $j_x$ denote the inclusion of $\{x\}$ into $X$. Then, for all  sufficiently small $\epsilon >
0$,
$$H^q(j_x^!\Adot)\cong\mathbb H^q_c(B^\circ_\epsilon(x) ; \Adot) \cong
\mathbb H^q( B^\circ_\epsilon(x) , B^\circ_\epsilon(x)
- x ; \Adot) \cong \mathbb H^q(X, X-x; \Adot)$$
and so, if $R$ is a field,
$$
\mathbf H^{-q}(\mathcal D\Adot)_x\cong  
\mathbb H^q( B^\circ_\epsilon(x) , B^\circ_\epsilon(x)
- x ; \Adot) \cong \mathbb H^q(X, X-x; \Adot).
$$

\vskip .2in

Suppose that the base ring $R$ is a field, that $X$ is a compact space, and that $\Adot$ is a self-dual complex on  $X$, i.e., there exists an isomorphism $\alpha:\Adot\rightarrow\mathcal D\Adot$. Then, for all $k$, there are isomorphisms
$$
\mathbb H^{-k}(X; \Adot)\cong\mathbb H^{-k}(X;\mathcal D\Adot)\cong Hom(\mathbb H_c^k(X;\Adot), R)\cong Hom(\mathbb H^k(X;\Adot), R)
$$
and so, there is an induced bilinear map:
$$
\mathbb H^k(X; \Adot)\times\mathbb H^{-k}(X; \Adot)\rightarrow R.
$$ 
When $k=0$, the map above is generally referred to as the {\it intersection pairing} or {\it intersection form}. This agrees with the standard usage of the term for a compact, even-dimensional, oriented, real $2n$-manifold, where $\Adot=\mathbf R^\bullet_{{}_X}[n]$, for then $\mathbb H^0(X; \Adot)=\mathbb H^0(X; \mathbf R^\bullet_{{}_X}[n])\cong H^n(X; R)$.

\vskip .2in

If $f: X \rightarrow Y$ and $g: Y \rightarrow Z$, then there are natural isomorphisms

$$R(g\circ f)_* \cong Rg_* \circ Rf_* \hskip .4in R(g\circ f)_! \cong Rg_! \circ Rf_!$$
\centerline{and} 
$$(g \circ f)^* \cong f^* \circ g^* \hskip .4in  
(g \circ f)^! \cong f^! \circ g^! .$$

\vskip .2in

Suppose that $f: Y \hookrightarrow X$ is inclusion of a subset.  Then, if $Y$ is
open, there is a natural isomorphism $f^! \cong f^*$.  More generally, if $\Adot\in \mathbf D^b_c(X)$ and $Y\cap\operatorname{supp}\Adot$ is open in  $\operatorname{supp}\Adot$, then $f^!\Adot \cong f^*\Adot$.  If $Y$ is closed, then  there are natural isomorphisms $Rf_! \cong f_! \cong f_* \cong Rf_*$.

\medskip

Excision has the following form: if $Y \subseteq U \subseteq X$, where $U$ is
open in $X$ and $Y$ is closed in $X$, then $$\mathbb H^k(X, X-Y;\Adot)
\cong \mathbb H^k(U, U-Y;\Adot).$$
This isomorphism on hypercohomology actually follows from a natural isomorphism of functors. Let $j$ denote the inclusion of $Y$ into $U$, and let $i$ denote the inclusion of $U$ into $X$. Then, there is a natural isomorphism $(ij)_!(ij)^!\cong i_*j_!j^!i^*$.

\vskip .2in

Suppose that $f:X\rightarrow Y$ is continuous. One might hope that, as an analog to the natural map $Rf_!\rightarrow Rf_*$, there would be a natural map $f^*\rightarrow f^!$; however, the situation is more complicated than that. If $\Adot$ and $\Bdot$ are in $\mathbf D^b_c(Y)$, then there is a natural map
$$
f^!\Adot\ \otimes^L\  f^*\Bdot \ \rightarrow \ f^!\big(\Adot \otimes^L\Bdot\big).
$$
In particular, taking $\Adot$ to be the constant sheaf $\mathbf R^\bullet_Y$, there is a natural map
$$
f^!\mathbf R^\bullet_Y\ \otimes^L \ f^*\ \rightarrow \ f^!.
$$
The complex $f^!\mathbf R^\bullet_Y$ is referred to as the {\it relative dualizing complex}, and is denoted by $\omega_{X/Y}$ (when the map $f$ is clear).

\smallskip

The map $f^!\mathbf R^\bullet_Y\ \otimes^L \ f^*\ \rightarrow \ f^!$ is an isomorphism if $f$ is a {\it topological submersion} of some dimension $r$ (we use $r$ here because it is a {\bf real} dimension), i.e., if every point in $X$ possesses an open neighborhood $\mathcal U$ such that $\mathcal V:=f(\mathcal U)$ is open and such that the restriction $f:\mathcal U\rightarrow \mathcal V$ is homeomorphic to the projection $\mathcal V\times\mathbb R^r\rightarrow\mathcal V$. If $f:X\rightarrow Y$ is a topological submersion of dimension $r$, and $X$ and $Y$ are orientable manifolds, then $\omega_{X/Y}\cong \mathbf R^\bullet_X[r]$ and so, for all $\Bdot\in \mathbf D^b_c(Y)$,
$$
f^*\Bdot[r]\  \cong\  f^!\Bdot.
$$

\medskip

If $g : Y \hookrightarrow X$ is the inclusion of an orientable submanifold into another orientable manifold, and $r$ is the real codimension of $Y$ in $X$, and $\Fdot \in \mathbf D^b_c(X)$ has locally constant
cohomology on $X$, then $g^!\Fdot$ has locally constant cohomology on $Y$
and $$g^*\Fdot[-r] \cong g^!\Fdot.$$
(There is an error here in \cite{G-M3}; they have the negation of the correct shift.)

\medskip

If $f : X \rightarrow Y$ is continuous, $\Adot \in \mathbf D^b_c(X)$, and $\Bdot \in \mathbf D^b_c(Y)$, then dual to the canonical maps
$$\Bdot\rightarrow Rf_*f^*\Bdot \text{ \ and \ } 
f^*Rf_*\Adot \rightarrow \Adot$$
are the canonical maps
$$Rf_!f^!\Bdot \rightarrow \Bdot\text{ \ and \ }
\Adot\rightarrow f^!Rf_!\Adot,$$
and, if $f$ is the inclusion of a locally closed subset, $\Adot\rightarrow f^!Rf_!\Adot$ is an isomorphism.
\vskip .4in

If 
$$\CD Z            @>\hat f>>         W\\ 
@V\hat\pi VV              @VV\pi V\\
X            @>f>>        S 
\endCD$$

\vskip .2in

\noindent is a pull-back diagram (fiber square, Cartesian diagram), then for all
$\Adot\in \mathbf D^b_c(X)$, $R\hat f_!\hat\pi^*\Adot \cong
\pi^*Rf_!\Adot$ (there is an error  in \cite{G-M3}; they have
lower $*$'s, not lower $!$'s, but see below for when these agree) and, dually,
$R\hat f_*\hat\pi^!\Adot \cong \pi^!Rf_*\Adot$.  In
particular, if $f$ is proper (and, hence, $\hat f$ is proper) or $\pi$ is the inclusion of
an open subset (and, hence, so is $\hat\pi$, up to homeomorphism), then $R\hat
f_*\hat\pi^*\Adot \cong \pi^*Rf_*\Adot$; this is also true if $W = 
S\times \mathbb C^n$ and
$\pi : W \rightarrow S$ is projection (and, hence, up to homeomorphism, $\hat\pi$ is
projection from $X \times \mathbb C^n$ to $X$). 

Still looking at the pull-back diagram above, we find that the natural map $Rf_!f^!\rightarrow\operatorname{id}$  yields a natural map $R\hat f_!\hat\pi^*f^! \cong\pi^*Rf_!f^!\rightarrow \pi^*$, which in turn yields natural maps
$$
\hat\pi^*f^!\rightarrow \hat f^!R\hat f_!\hat\pi^*f^!\rightarrow \hat f^!\pi^*.
$$

\vskip .1in

If we have $\Adot\in \mathbf D^b_c(X)$ and $\Bdot\in
\mathbf D^b_c(W)$, then we let $\Adot\boxtimes^L_{{}_S}\Bdot:=\hat\pi^*\Adot\otimes^L \hat f^*\Bdot$,
assuming that the maps $\hat\pi$ and $\hat f$ are clear. If $S$ is a point, so that $Z\cong X\times W$, then we omit the $S$
in the notation and write simply $\Adot\boxtimes^L\Bdot$.

\vskip .2in

There is a K\"unneth formula, which we now state in its most general form, in terms of maps over a base space $S$. Suppose
that we have two maps $f_1:X_1\rightarrow Y_1$ and $f_2:X_2\rightarrow Y_2$ over $S$, i.e., there are commutative diagrams

\vskip .1in

\vbox{
$\begin{CD}\hskip 1in X_1 @>f_1>> \ Y_1 
\hskip 1.95in X_2 @>f_2>> \ Y_2\end{CD}$
\vskip .03in 
$\hskip 1.1in r_1\searrow\hskip .1in \swarrow t_1$\hskip .9in
and $\hskip 1in  r_2\searrow\hskip .1in \swarrow t_2$
\vskip .01in
\hskip 1.4in S\hskip 2.8in S\hskip 1in.
}

\vskip .1in

\noindent Then, there is an induced map $f=f_1\times_{{}_S}f_2:X_1\times_{{}_S}X_2\rightarrow Y_1\times_{{}_S} Y_2$. If
$\Adot\in \mathbf D^b_c(X_1)$ and $\Bdot\in
\mathbf D^b_c(X_2)$, there is the {\bf K\"unneth isomorphism}
$$
Rf_!(\Adot\boxtimes^L_{{}_S}\Bdot)\cong R{f_1}_! \Adot\ \boxtimes^L_{{}_S}\ R{f_2}_! \Bdot.
$$

Using the above notation, if $S$ is a point and $\Fdot\in \mathbf D^b_c(Y_1)$ and $\Gdot\in
\mathbf D^b_c(Y_2)$, then there is a natural isomorphism (the {\bf adjoint K\"unneth isomorphism})
$$
f^!(\Fdot\boxtimes^L\Gdot)\cong {f_1}^! \Fdot\ \boxtimes^L\ {f_2}^! \Gdot.
$$
If we let $q_1$ and $q_2$ denote the projections from $Y_1\times Y_2$ onto $Y_1$ and $Y_2$, respectively, then the adjoint
K\"unneth formula can be proved by using the following natural isomorphism twice
$$
\mathcal D\Fdot \boxtimes^L\Gdot\cong R\mathbf{Hom}^\bullet(q_1^*\Fdot, q_2^!\Gdot).
$$

\bigskip

Let $Z$ be a locally closed subset of an analytic space $X$. There are two derived functors, associated to $Z$, that we
wish to describe: the derived functors of restricting-extending to $Z$, and of taking the sections supported on $Z$. Let $k:Z\hookrightarrow X$ denote the inclusion.

If $\mathbf A$ is a (single) sheaf on $X$, then the restriction-extension of $\mathbf A$ to $Z$, $(\mathbf A)_Z$, is given by
$k_!k^*(\mathbf A)$. Thus, up to isomorphism, $(\mathbf A)_Z$ is characterized by $((\mathbf A)_Z)_{|_Z}\cong {\mathbf A}_{|_Z}$
and $((\mathbf A)_Z)_{|_{X-Z}}=0$. This functor is exact, and so we also denote the derived functor by $( )_Z$.

\vskip .1in

Now, we want to define the sheaf of sections of $\mathbf A$ supported by $Z$, $\Gamma_{{}_Z}(\mathbf A)$. If $\mathcal U$ is an open
subset of $X$ which contains $Z$ as a closed subset (of $\mathcal U$), then we define 
$$
\Gamma_{{}_Z}(\mathcal U; \mathbf A):= \operatorname{ker}\big\{\Gamma(\mathcal U; \mathbf A)\rightarrow \Gamma(\mathcal U-Z; \mathbf A)\big\}.
$$
Up to isomorphism, $\Gamma_{{}_Z}(\mathcal U; \mathbf A)$ is independent of the open set $\mathcal U$ (this uses that $\mathbf A$ is a
sheaf, not just a presheaf). The sheaf $\Gamma_{{}_Z}(\mathbf A)$ is defined by: for all open $\mathcal U\subseteq X$,
$\Gamma(\mathcal U; \Gamma_{{}_Z}(\mathbf A)):= \Gamma_{{}_{\mathcal U\cap Z}}(\mathcal U; \mathbf A)$.

One easily sees that $\operatorname{supp} \Gamma_{{}_Z}(\mathbf A)\subseteq \overline{Z}$. The functor $\Gamma_{{}_Z}()$ is left exact; of course, we denote the right derived functor by $R\Gamma_{{}_Z}()$. There is a canonical isomorphism 
$$k^!\cong k^* R\Gamma_{{}_Z}.$$
If $Z$ is closed, then
$R\Gamma_{{}_Z}\cong k_!k^!$. If $Z$ is open, then $R\Gamma_{{}_Z}\cong Rk_*k^*$.

\bigskip

\noindent{\bf Avoiding Injective Resolutions:}

\vskip .2in

To calculate right derived functors from the definition, one must use injective resolutions. However, this is inconvenient
in many proofs if some functor involved in the proof does not take injective complexes to injective complexes. There are (at
least) four ``devices'' which come to our aid, and enable one to prove many of the isomorphisms described earlier; these
devices are {\bf fine resolutions}, {\bf flabby resolutions}, {\bf $c$-soft resolutions}, and {\bf injective subcategories
with respect to a functor}.

If $T$ is a left-exact functor on the category of sheaves on $X$, then the right derived functor $RT$ is defined by applying
$T$ term-wise to the sheaves in a canonical injective resolution. The importance of saying that a certain subcategory of the
category of sheaves on $X$ is injective with respect to $T$ is that one may take a resolution in which the individual
sheaves are in the given subcategory, then apply $T$ term-wise, and end up with a complex which is canonically isomorphic to
that produced by $RT$.

\vskip .1in

Recall that a single sheaf $\mathbf A$ on $X$ is:

\vskip .1in

\noindent $\bullet$\ \    fine, if partitions of unity of $\mathbf A$ subordinate to any given locally finite open cover of $X$ exist;

\vskip .1in

\noindent  $\bullet$ \ \ flabby, if for every open subset $\mathcal U\subseteq X$, the restriction homomorphism $\Gamma(X; \mathbf
A)\rightarrow\Gamma(\mathcal U; \mathbf A)$ is a surjection;

\vskip .1in

\noindent  $\bullet$\ \ $c$-soft, if for every compact subset $\mathcal K\subseteq X$, the restriction homomorphism $\Gamma(X; \mathbf
A)\rightarrow\Gamma(\mathcal K; \mathbf A)$ is a surjection;

\vskip .1in

Injective sheaves are flabby, and flabby sheaves are $c$-soft. In addition, fine sheaves are $c$-soft.

\vskip .1in

 The subcategory of $c$-soft sheaves is injective with respect to the functors $\Gamma(X;*)$, $\Gamma_c(X;*)$,
and $f_!$. The subcategory of flabby sheaves is injective with respect to the functor 
$f_*$.

\vskip .1in

If $\mathbf A^\bullet$ is a bounded complex of sheaves, then a bounded $c$-soft resolution of $\mathbf A^\bullet$ is given by
$\mathbf A^\bullet\rightarrow\mathbf A^\bullet\otimes\mathbf S^\bullet$, where $\mathbf S^\bullet$ is a $c$-soft, bounded above,
resolution of the base ring (which always exists in our context).

\bigskip

\noindent{\bf Triangles:}

\medskip

$\mathbf D^b_c(X)$ is an additive category, but is not an Abelian category.  In place of short
exact sequences, one has distinguished triangles, just as we did in $\mathbf K^b(\mathcal C)$.  A triangle
of morphisms in
$\mathbf D^b_c(X)$

\vskip .1in

\vbox{
$\hskip 2.5in \Adot \longrightarrow \ \Bdot$ \vskip .03in
$\hskip 2.35in \ _{[1]} \ \nwarrow\hskip .15in \swarrow$ \vskip .03in $\hskip
2.75in \mathbf C^\bullet$
}

\vskip .1in

\noindent (the $[1]$ indicates a morphism shifted by one, i.e., a morphism $\mathbf
C^\bullet \rightarrow \Adot[1]$)  is called {\it distinguished} if it
is isomorphic in $\mathbf D^b_c(X)$ to a diagram of sheaf maps

\vskip .1in

\vbox{
$\hskip 2.5in {\widetilde\Adot}\xrightarrow{ \ \phi  \  } \ {\widetilde\Bdot}$ \vskip .03in $\hskip 2.35in \ _{[1]} \ \nwarrow\hskip .15in \swarrow$
\vskip .03in $\hskip 2.75in \mathbf M^\bullet$,
}

\vskip .1in

\noindent where $\mathbf M^\bullet$ is the algebraic mapping cone of $\phi$ and
$\widetilde\Bdot \rightarrow \mathbf M^\bullet \rightarrow \widetilde\Adot[1]$ are
the canonical maps.  The ``in-line" notation
for a triangle is \hbox{$\Adot\rightarrow\Bdot\rightarrow\mathbf C^\bullet\rightarrow\Adot[1]$} or 
$$\begin{CD}\Adot\rightarrow\Bdot\rightarrow\mathbf C^\bullet @>[1] >>.\end{CD}$$

\vskip .1in

Any short exact sequence of complexes becomes a distinguished triangle in 
$\mathbf D^b_c(X)$.  Any edge of a distinguished triangle determines the triangle
up to (non-canonical) isomorphism in $\mathbf D^b_c(X)$; more specifically, we can ``turn'' the distinguished triangle:
$$\begin{CD}\Adot@>\alpha>>\Bdot@>\beta>>\mathbf C^\bullet@>\gamma>>\Adot[1]\end{CD}$$ is a distinguished triangle if and only if
$$\begin{CD}\Bdot@>\beta>>\Cdot@>\gamma>>\mathbf A^\bullet [1]@>-\alpha[1]>>\Bdot[1]\end{CD}$$ is a distinguished triangle.

\bigskip

Given two distinguished triangles and maps $u$ and $v$ which make the left-hand square of the following diagram commute

\vbox{
$$\Adot\rightarrow\Bdot\rightarrow\Cdot\rightarrow\Adot[1]\ $$
\vskip -.15in
$$\downarrow u\hskip .2in \downarrow v\hskip .5in\downarrow u[1]$$
\vskip -.1in
$$\widetilde\Adot\rightarrow\widetilde\Bdot\rightarrow\widetilde\Cdot\rightarrow\widetilde\Adot[1],$$}

\smallskip

\noindent there exists a (not necessarily unique) $w:\Cdot\rightarrow\widetilde\Cdot$ such that

$$\Adot\rightarrow\Bdot\rightarrow\Cdot\rightarrow\Adot[1]\ $$
\vskip -.2in
$$\downarrow u\hskip .17in \downarrow v\hskip .17in\downarrow w\hskip .17in\downarrow u[1]$$
\vskip -.1in
$$\widetilde\Adot\rightarrow\widetilde\Bdot\rightarrow\widetilde\Cdot\rightarrow\widetilde\Adot[1]\ $$
also commutes. We say that the original commutative square embeds in a morphism of distinguished triangles.

\vskip .1in

We will now give the {\it octahedral lemma}, which allows one to realize an isomorphism between mapping cones of two
composed maps.

Suppose that we have two distinguished triangles 
$$\begin{CD}\Adot@>f>>\Bdot@>g>>\Cdot@>h>>\Adot[1]\end{CD}$$ and $$\begin{CD}\Bdot@>\beta>>\Edot@>\gamma>>\Fdot@>\delta>>\Bdot[1].\end{CD}$$ 

\bigskip

\noindent Then, there exists a complex $\Mdot$ and two distinguished triangles
$$
\begin{CD}\Adot@>\beta\circ f>>\Edot@>\tau>>\Mdot@>\omega>>\Adot[1]\end{CD}
$$
and
$$
\begin{CD}\Cdot@>\sigma>>\Mdot@>\nu>>\Fdot@>g[1]\circ\delta>>\Cdot[1].\end{CD}
$$
\bigskip

\noindent such that the following diagram commutes

\pagebreak

$$\begin{CD}
\Adot@>\ \ f\ \ >>\Bdot@>\ \ g\ \ >>\Cdot@>\ \ h\ \ >>\Adot[1]\end{CD}
$$
$$\begin{CD}
\operatorname{id}\downarrow\hskip .5in\beta\downarrow\hskip .5in\sigma\downarrow\hskip
.6in\operatorname{id}\downarrow\hskip .2in\hbox{}\end{CD}
$$
$$\begin{CD}
\Adot@>\ \beta\circ f\ >>\Edot@>\ \ \tau\ \ >>\Mdot@>\ \ \omega\ \ >>\Adot[1]\end{CD}
$$
$$\begin{CD}
f\downarrow\hskip .5in\operatorname{id}\downarrow\hskip .5in\nu\downarrow\hskip
.45in f[1]\downarrow\hskip .2in\hbox{}\end{CD}
$$
$$
\begin{CD}\Bdot@>\ \ \ \beta\ \ >>\Edot@>\ \ \gamma\ \ >>\Fdot@>\ \ \delta\ \ >>\Bdot[1]\end{CD}
$$
$$\begin{CD}
g\downarrow\hskip .5in\tau\downarrow\hskip .45in\operatorname{id}\downarrow\hskip
.4in g[1]\downarrow\hskip .15in\hbox{}\end{CD}
$$
$$\begin{CD}
\hbox{}\hskip .05in \Cdot@>\ \ \sigma\ \ >>\Mdot@>\ \ \nu\ \ >>\Fdot@>g[1]\circ\delta>>\Cdot[1].\end{CD}
$$

\bigskip

\noindent It is somewhat difficult to draw this in its octahedral form (and worse to type it); moreover, it is no easier to read the
relations from the octahedron. However, the interested reader can give it a try: the octahedron is formed by gluing together two
pyramids along their square bases. One pyramid has $\Bdot$ at its top vertex, with $\Adot$, $\Cdot$, $\Edot$, and $\Fdot$ at the
vertices of its base, and has the original two distinguished triangles as opposite faces. The other pyramid has $\Mdot$ at its
top vertex, with $\Adot$, $\Cdot$, $\Edot$, and $\Fdot$ at the vertices of its base, and has the other two distinguished
triangles (whose existence is asserted in the lemma) as opposite faces. The two pyramids are joined together by matching the
vertices of the two bases, forming an octahedron in which the faces are alternately distinguished and commuting.

\vskip .2in

A distinguished triangle determines long exact sequences on cohomology and
hypercohomology:

$$\cdots \rightarrow \mathbf H^p(\Adot) \rightarrow \mathbf H^p(\Bdot)
 \rightarrow \mathbf H^p(\mathbf C^\bullet) \rightarrow \mathbf H^{p+1}(\Adot) \rightarrow \cdots$$ $$\cdots \rightarrow \mathbb H^p(X;\Adot)
\rightarrow \mathbb H^p(X;\Bdot)
 \rightarrow \mathbb H^p(X;\mathbf C^\bullet) \rightarrow \mathbb H^{p+1}(X;\Adot) \rightarrow \cdots .$$

\vskip .1in

If $f : X \rightarrow Y$ and $\Fdot \in \mathbf D_c^b(X)$, then the functors $Rf_* , Rf_! ,
f^* , \ f^!,\ \text{and}\ \Fdot\otimes^L*$ all take distinguished triangles to
distinguished triangles (with all arrows in the same direction and the shift in the same place).

As for $R\mathbf{Hom}^\bullet$, if $\Fdot \in \mathbf D_c^b(X)$ and 
 \hbox{$\Adot\rightarrow\Bdot\rightarrow\mathbf C^\bullet\rightarrow\Adot[1]$} 
is a distinguished triangle in $\mathbf D_c^b(X)$, then we have distinguished triangles

\vskip .2in

{\moveright .5in\vbox{
$R\mathbf{Hom}^\bullet(\Fdot, \Adot) \longrightarrow
R\mathbf{Hom}^\bullet(\Fdot,\Bdot)  \hskip .23in 
R\mathbf{Hom}^\bullet(\Adot, \Fdot) \longleftarrow
R\mathbf{Hom}^\bullet(\Bdot, \Fdot)$ 
\vskip .03in 
$\hbox{}\hskip .33in\ _{[1]}
\ \nwarrow\hskip .42in \swarrow$\hskip 1in and  $\hskip .5in \ _{[1]} \ \searrow\hskip
.4in \nearrow$ 
\vskip .03in 
$\hskip .45in R\mathbf{Hom}^\bullet(\Fdot, \mathbf C^\bullet) \hskip 1.84in 
R\mathbf{Hom}^\bullet(\mathbf C^\bullet, \Fdot)$.
}}

\vskip .2in

By applying the right-hand triangle above to the special case where $\Fdot =
\mathbb D_X^\bullet$, we find that the dualizing functor $\mathcal D$ also takes
distinguished triangles to  distinguished triangles, but with a reversal of arrows, i.e., if we have a distinguished triangle \hbox{$\Adot\rightarrow\Bdot\rightarrow\mathbf C^\bullet\rightarrow\Adot[1]$} in $\mathbf D_c^b(X)$, 
then, by dualizing, we have distinguished triangles

\vskip .2in

\vbox{
$\hskip 1.15in \mathcal D\Adot \longleftarrow \ \mathcal D\Bdot 
\hskip 2.2in \mathcal D\mathbf C^\bullet \longrightarrow \ \mathcal D\Bdot$
\vskip .03in $\hskip 1.07in\ _{[1]} \ \searrow\hskip .23in \nearrow$\hskip 1.2in
or  $\hskip 1in \ _{[1]} \ \nwarrow\hskip .22in \swarrow$ \vskip .03in $\hskip
1.45in \mathcal D\mathbf C^\bullet \hskip 2.84in \mathcal D\Adot$.
}

\vskip.2in

There are (at least) six distinguished triangles associated to the functors  $()_Z$ and $R\Gamma_{{}_Z}$. Let $\Fdot$ be in $\mathbf D_c^b(X)$, $\mathcal U_1$ and $\mathcal U_2$ be open subsets of $X$, $Z_1$ and $Z_2$ be closed subsets of
$X$, $Z$ be locally closed in $X$, and $Z^\prime$ be closed in $Z$. Then, we have the following distinguished triangles:

\vskip .1in

$$
\begin{CD}R\Gamma_{{}_{\mathcal U_1\cup \mathcal U_2}}(\Fdot)\rightarrow R\Gamma_{{}_{\mathcal U_1}}(\Fdot)\oplus R\Gamma_{{}_{\mathcal U_2}}(\Fdot)\rightarrow R\Gamma_{{}_{\mathcal U_1\cap \mathcal U_2}}(\Fdot)@>[1]>>\end{CD}
$$

\vskip .1in

$$
\begin{CD}R\Gamma_{{}_{Z_1\cap  Z_2}}(\Fdot)\rightarrow R\Gamma_{{}_{ Z_1}}(\Fdot)\oplus R\Gamma_{{}_{ Z_2}}(\Fdot)\rightarrow R\Gamma_{{}_{Z_1\cup Z_2}}(\Fdot)@>[1]>>\end{CD}
$$

\vskip .1in

$$
\begin{CD}(\Fdot)_{{}_{\mathcal U_1 \cap\mathcal U_2}}\rightarrow (\Fdot)_{{}_{\mathcal U_1}}\oplus (\Fdot)_{{}_{\mathcal U_2}}\rightarrow (\Fdot)_{{}_{\mathcal U_1 \cup\mathcal U_2}}@>[1]>>\end{CD}
$$

\vskip .1in

$$
\begin{CD}(\Fdot)_{{}_{Z_1 \cup Z_2}}\rightarrow (\Fdot)_{{}_{Z_1}}\oplus (\Fdot)_{{}_{Z_2}}\rightarrow (\Fdot)_{{}_{Z_1 \cap Z_2}}@>[1]>>\end{CD}
$$

\vskip .1in

$$
\begin{CD}R\Gamma_{{}_{Z^\prime}}(\Fdot)\rightarrow R\Gamma_{{}_{Z}}(\Fdot)\rightarrow R\Gamma_{{}_{Z-Z^\prime}}(\Fdot)@>[1]>>\end{CD}
$$

\vskip .1in

$$
\begin{CD}(\Fdot)_{{}_{Z-Z^\prime}}\rightarrow (\Fdot)_{{}_{Z}}\rightarrow (\Fdot)_{{}_{Z^\prime}}@>[1]>>.\end{CD}
$$

\vskip .3in

If  $j : Y \hookrightarrow X$ is the inclusion of a closed subspace and $i : U
\hookrightarrow X$ the inclusion of the open complement, then for all $\Adot \in \mathbf D^b_c(X)$, the last two triangles above give us distinguished triangles

\vskip .1in

\vbox{
$\hskip .65in Ri_!i^!\Adot \longrightarrow \ \Adot  \hskip
2.2in Rj_!j^!\Adot \longrightarrow \  \Adot$ 
\vskip .03in
$\hskip .59in\ _{[1]} \ \nwarrow\hskip .35in \swarrow$\hskip 1.2in and  $\hskip
1in \ _{[1]} \ \nwarrow\hskip .22in \swarrow$ 
\vskip .03in 
$\hskip .95in Rj_*j^*\Adot \hskip 2.75in Ri_*i^*\Adot$,
}

\vskip .1in

\noindent where the second triangle can be obtained from the first by dualizing.
(Note that $Ri_! = i_!$, $Rj_* = j_* = j_! = Rj_!$, and $i^! = i^*$.)  The
associated long exact sequences on hypercohomology are those for the pairs $\mathbb
H^*(X, Y;\Adot)$ and $\mathbb H^*(X,U;\Adot)$, respectively. It is worth noting that the natural map  $Rj_*j^*[-1]\Adot\rightarrow Ri_!i^!\Adot$ induces the zero map on stalk cohomology, but would usually {\bf not} be the zero map in the derived category.

From either of these triangles, and the fact that $j^*i_!=0$, one easily obtains a natural isomorphism 
$$j^!Ri_![1]\cong j^*Ri_*.$$

Applying $j^*$ to the triangle on the right above, we obtain the distinguished triangle 

\medskip

\vbox{
\centerline{$j^!\Adot \longrightarrow \  j^*\Adot$}
\vskip .03in
 \centerline{$\hskip -0.3in\ _{[1]} \ \nwarrow\hskip .3in \swarrow$}
\vskip .05in 
\centerline{$j^*Ri_*i^*\Adot$,}
}

\medskip

\noindent which tells us how $j^!\Adot$ and $j^*\Adot$ are related.

\vskip .2in

As in our earlier discussion of the octahedral lemma, all of the morphisms of the last two paragraphs fit into the {\it
fundamental octahedron of the pair $(X,Y)$}. The four distinguished triangles making up the fundamental octahedron are the top
pair 
$$
Ri_!i^!\Adot\rightarrow\Adot\rightarrow Rj_*j^*\Adot\rightarrow Ri_!i^!\Adot[1]
$$
and
$$
\Adot\rightarrow Ri_*i^*\Adot\rightarrow Rj_!j^!\Adot[1]\rightarrow \Adot[1]
$$
and the bottom pair
$$
Ri_!i^!\Adot\rightarrow Ri_*i^*\Adot\rightarrow \Mdot\rightarrow Ri_!i^!\Adot[1]
$$
and
$$
Rj_*j^*\Adot\rightarrow \Mdot \rightarrow Rj_!j^!\Adot[1]\rightarrow Rj_*j^*\Adot[1],
$$

\medskip

\noindent where $\Mdot\cong Rj_!j^!Ri_!i^!\Adot[1]\cong Rj_*j^*Ri_*i^*\Adot$.

\bigskip

Consider again the situation in which we have a pull-back diagram

$$\CD Z            @>\hat f>>         W\\ 
@V\hat\pi VV              @VV\pi V\\
X            @>f>>        S,
\endCD$$

\vskip .2in

\noindent where we saw  that the natural map $Rf_!f^!\rightarrow\operatorname{id}$  yields a natural map $R\hat f_!\hat\pi^*f^! \cong\pi^*Rf_!f^!\rightarrow \pi^*$, which in turn yields natural maps
$$
\hat\pi^*f^!\rightarrow \hat f^!R\hat f_!\hat\pi^*f^!\rightarrow \hat f^!\pi^*.
$$
If $f$ is the inclusion of an open or closed subset, and $l:S-X\rightarrow S$ is the inclusion of the complement, then, for all $\Adot\in \mathbf D^b_c(S)$, we have the distinguished triangle
$$
Rf_!f^!\Adot\rightarrow\Adot\rightarrow Rl_*l^*\Adot\rightarrow Rf_!f^!\Adot[1],
$$
and so, have the induced triangle
$$
\hat f^!\pi^*Rf_!f^!\Adot\rightarrow\hat f^!\pi^*\Adot\rightarrow \hat f^!\pi^*Rl_*l^*\Adot\rightarrow \hat f^!\pi^*Rf_!f^!\Adot[1],
$$
which, up to isomorphism, is
$$
\hat\pi^*f^!\Adot\rightarrow\hat f^!\pi^*\Adot\rightarrow \hat f^!\pi^*Rl_*l^*\Adot\rightarrow \hat\pi^*f^!\Adot[1];
$$
it follows that, if $\pi^*Rl_*l^*\Adot\cong \pi^!Rl_*l^*\Adot$, then $ \hat f^!\pi^*Rl_*l^*\Adot=0$, and so $\hat\pi^*f^!\Adot\rightarrow\hat f^!\pi^*\Adot$ is an isomorphism. This would be the case, for instance, if $f$ and $\pi$ are inclusions of closed subsets, and $S=X\cup W$.

\vskip .4in

\noindent\S2. {\bf Perverse Sheaves}  

\bigskip

Let $X$ be a complex analytic
space, and for each $\mathbf x \in X$, let $j_\mathbf x :
\mathbf x \hookrightarrow X$ denote the inclusion. 

\smallskip

If $\Fdot \in \mathbf D_c^b(X)$, then the
{\it support of $\mathbf H^i(\Fdot)$} is the closure in $X$ of 
$$\{\mathbf x\in X |\ \mathbf
H^i(\Fdot)_\mathbf x \neq 0\} = \{\mathbf x\in X |\ \mathbf
H^i(j^*_\mathbf x \Fdot) \neq 0\} ;$$
we denote this by $\operatorname{supp^i}\Fdot$.

The {\it i-th cosupport of $\Fdot$} is the closure in $X$ of
$$\{\mathbf x\in X |\ \mathbf H^i(j^!_\mathbf x \Fdot) \neq 0\} \ = \ \{\mathbf x\in X |\ \mathbb H^i( B^\circ_\epsilon(x) , 
B^\circ_\epsilon(x) - x ;\ \Fdot) \neq 0\};$$ we denote this by $\operatorname{cosupp^i}\Fdot$.

\vskip .3in

If the base ring, $R$, is a field, then $\operatorname{cosupp^i}\Fdot = \operatorname{supp^{-i}}\mathcal D\Fdot$.

\vskip .3in

\noindent {\bf Definition}:  Let $X$ be a complex analytic
space (not necessarily pure-dimensional). Then, $\Pdot \in \mathbf D_c^b(X)$
is {\it  perverse} provided that for all $i$:

\vskip .1in

\noindent (support) \hskip .4in $\operatorname{dim}(\operatorname{supp^{-i}}\Pdot) \leqslant
i$;

\vskip .1in

\noindent (cosupport) \hskip .28in $\operatorname{dim}(\operatorname{cosupp^i}\Pdot)
\leqslant i$,

\vskip .1in

\noindent  where we set the dimension of the empty
set to be $-\infty$.

\vskip .1in

This definition is equivalent to:  let $\mathcal S:=\{S_\alpha\}$ be any stratification of $X$
with respect to which $\Pdot$ is constructible, and let $s_\alpha : S_\alpha
\hookrightarrow X$ denote the inclusion. Then,

\vskip .1in

\noindent (support) \hskip .4in $\mathbf H^k(s_\alpha^*\Pdot) = 0$ \hskip
.25in for $k > -\text{dim}_{{}_\mathbb C}S_\alpha$;

\vskip .1in

\noindent (cosupport) \hskip .28in $\mathbf H^k(s_\alpha^!\Pdot) = 0$ \hskip
.25in for $k < -\text{dim}_{{}_\mathbb C}S_\alpha$.

\vskip .1in

(There is a missing minus sign in \cite{G-M2}, 6.A.5.)

\vskip 0.1in

We will give another equivalent characterization of the support and cosupport conditions, which also uses a stratification $\mathcal S$ with respect to which $\Pdot$ is constructible. We define the {\bf upper filtration} of $X$ by: for all $m\geq0$, let $U^m:=\bigcup_{S_\alpha\in\mathcal S, \dim S_\alpha\geq m}S_\alpha$. We define the {\bf lower filtration} of $X$ by: for all $m\geq 0$, let $L^m:=\bigcup_{S_\alpha\in\mathcal S, \dim S_\alpha\leq m}S_\alpha$. For all $m\geq 0$, let $u_m: U^m\hookrightarrow X$ and $\ell_m: L^m\hookrightarrow X$ denote the inclusions.

Then, characterizations of the support and cosupport conditions are:

\vskip 0.1in

\noindent (support) for all $m\geq 0$, for all $k>-m$, $\mathbf H^k(u_m^*\Pdot)=0$;

\vskip 0.1in

\noindent (cosupport) for all $m\geq0$, for all $k<-m$, $\mathbf H^k(\ell_m^!\Pdot)=0$.

\vskip .3in

From the definition, it is clear that being  perverse is a
local property.

\vskip .2in

If the base ring $R$ is, in fact, a field, then the support and cosupport conditions can
be written in the following form, which is symmetric with respect to dualizing:

\vskip .1in

\noindent (support) \hskip .4in $\operatorname{dim} (\operatorname{supp^{-i}}\Pdot)
\leqslant i$;

\vskip .1in

\noindent (cosupport) \hskip .28in $\operatorname{dim} (\operatorname{supp^{-i}}\mathcal D\Pdot) \leqslant i$.

\vskip .2in

Suppose that $\Pdot$ is  perverse on $X$, $(X,x)$ is
locally embedded in $\mathbb C^n$, $S$ is a stratum of a Whitney
stratification with respect to which $\Pdot$ is constructible, and $x
\in S$.  Let $M$ be a normal slice of $X$ at $x$; that is, let $M$ be a smooth
submanifold of $\mathbb C^n$ of dimension $n - \text{dim }S$ (where $\text{dim } S$ is the complex dimension) which
transversely intersects $S$ at $x$.  Then, for some open neighborhood $U$ of
$x$ in $X$, $\Pdot_{|_{X\cap M\cap U}}[-\text{dim }S]$ is 
 perverse on $X\cap M\cap U$.

\vskip .1in

Let $\Pdot \in \mathbf D^b_c(X)$ be perverse. One can use the above normal slicing proposition to
prove that  $\mathbf H^i(\Pdot) = 0$ for all $i < -\text{dim }X$. For all $x\in X$, let $j_x$ denote the inclusion of $\{x\}$ into $X$. Then this last statement, together with the support condition, tells us that $H^k(j_x^*\Pdot)=0$ except possibly for $-\text{dim }X\leq k\leq 0$. The ``dual'' statement is also true: $H^k(j_x^!\Pdot)=0$ except possibly for $0\leq k\leq \text{dim }X$.

\vskip 0.2in

There is a global version of the statements in the previous paragraph; this global version sometimes goes by the name of {\it Artin Vanishing}: Let $\Pdot$ be a perverse sheaf on an affine analytic space $X$. Then, $\mathbb H^k(X; \Pdot)=0$, except possibly for $-\text{dim }X\leq k\leq 0$, and $\mathbb H_c^k(X; \Pdot)=0$, except possibly for $0\leq k\leq \text{dim }X$.

\vskip .2in

Suppose $\Pdot \in \mathbf D^b_c(X^n)$.  Let $\Sigma = \text{supp }
\Pdot$.  Then, 
$\Pdot$ is  perverse on $X$ if and only if 
$\Pdot_{|_\Sigma}$ is  perverse on $\Sigma$.  Another way of saying this is: if $m: \Sigma \rightarrow X$ is the inclusion of a closed
subspace, then $\mathbf Q^\bullet$ is  perverse on $\Sigma$ if and only if
$m_*\mathbf Q^\bullet\cong m_!\mathbf Q^\bullet$ is  perverse on $X$. It follows that if $\Pdot$ is  perverse, then $\mathbf
H^i(\Pdot) = 0$ except possibly for $-\text{dim }\Sigma \leqslant i \leqslant 0$. Furthermore, again letting $j_x$ denote the inclusion of $\{x\}$ into $X$, $H^i(j_x^!\Pdot)=0$ except possibly for $0\leq i\leq \text{dim }\Sigma$.

\vskip .2in

If $X$ is (not necessarily purely) $n$-dimensional, then it is trivial to show that the shifted constant sheaf $\mathbf R_{{}_X}^\bullet[n]$ satisfies the support condition.  If $X$ is a curve, then the constant sheaf $\mathbf R_{{}_X}^\bullet[1]$ is perverse. If $X$ is a surface, then $\mathbf R_{{}_X}^\bullet[2]$ is perverse if and only if, for all $p\in X$, for all sufficiently small open balls $B$ (in some Riemannian metric), centered at $p$, $B\cap X-\{p\}$ is connected; this would, for instance, be the case if $X$ were a surface which is locally irreducible everywhere.     If $X$ is a purely $n$-dimensional local complete intersection, then $\mathbf R_{{}_X}^\bullet[n]$ is perverse.  

\vskip .2in

A converse to the normal slicing proposition is: if $\pi : X \times \mathbb C^s \rightarrow X$ is projection and $\Pdot$ is  perverse on $X$, then
$\pi^*\Pdot[s]$ is  perverse.

\vskip .3in

The other basic example of a perverse sheaf that we wish to give is
that of intersection cohomology with local coefficients (with the perverse 
indexing, i.e., the only non-zero cohomology is in degrees less than or equal to zero). Note that the definition below is shifted by $-\operatorname{dim}_{\mathbb
C}X$  from the definition in \cite{G-M3}, and yields a perverse sheaf which has possibly non-zero stalk cohomology only in degrees between
$-\operatorname{dim}_{\mathbb C}X$ and $-1$, inclusive (unless the space has isolated points).

\vskip .1in

Let $X$ be a complex analytic set, let $X^\circ$ be an open dense subset of the smooth part of $X$, and let $\mathcal L^\bullet$ be a perverse sheaf on $X^\circ$, such that the restriction of $\mathcal L^\bullet$ to each connected component $C$ of $X^\circ$ is a local system on $C$, shifted by $\dim C$ (i.e., shifted into degree $-\dim C$). Let $Y:=X-X^\circ$ and let $j:Y\hookrightarrow X$ denote the (closed) inclusion.

 Then, in $\mathbf
D^b_c(X)$, there is an object,  $\mathbf {IC}_{{}_X}^\bullet(\mathcal L^\bullet)$, called the {\it
intersection cohomology with coefficients in $\mathcal L$}  which is uniquely determined up to
quasi-isomorphism by:

\vskip .1in

\noindent $1$)\hskip .2in $\mathbf {IC}_{{}_X}^\bullet(\mathcal L^\bullet)_{|_{X^\circ}} \ = \ \mathcal L^\bullet$;

\medskip

for all $i$,

\medskip

\noindent $2$)\hskip .2in $\operatorname{dim}{\operatorname{supp}}^{-i}\left(j^*[-1]\mathbf{IC}_{{}_X}^\bullet(\mathcal L^\bullet)\right)  \ = \ \operatorname{dim}\overline{\big\{\mathbf x\in Y\ |\ H^{-(i+1)}\big(\mathbf{IC}_{{}_X}^\bullet(\mathcal L^\bullet)\big)_x\neq 0\big\}} \ \leq i$;

\vskip .1in

\noindent $3$)\hskip .2in $\operatorname{dim}{\operatorname{cosupp}}^i\left(j^![1]\mathbf{IC}_{{}_X}^\bullet(\mathcal L^\bullet)\right)=\operatorname{dim}\overline{\big\{\mathbf x\in Y  \ | \ \mathbb H^{i+1}\big( B^\circ_\epsilon(x), 
B^\circ_\epsilon(x) - x ;\ \mathbf{IC}_{{}_X}^\bullet(\mathcal L^\bullet)\big) \neq 0\big\}} \  \leq i$.

\vskip .2in

$\mathbf {IC}_{{}_X}^\bullet(\mathcal L^\bullet)$ is a perverse sheaf. The uniqueness assertion implies that $\mathbf {IC}_{{}_X}^\bullet(\mathcal L^\bullet)$ is the direct sum of the extension by zero of the intersection cohomology complex of each irreducible component of $X$ using $\mathcal L^\bullet$ restricted to an open dense set of that component.   If $X$ is purely $n$-dimensional, and $\mathcal L$ is a local system (in degree $0$) on an open dense subset of the smooth part of $X$, then one usually writes simply $\mathbf {IC}_{{}_X}^\bullet(\mathcal L)$ in place of $\mathbf {IC}_{{}_X}^\bullet(\mathcal L^\bullet[n])$. See Section 5 for more on
$\mathbf {IC}_{{}_X}^\bullet(\mathcal L^\bullet)$.

\vskip .3in

\noindent{\bf The Category of Perverse Sheaves} (see, also, section 5)

\vskip .2in

The {\it category of perverse sheaves} on $X$, ${\operatorname{Perv}}(X)$, is the full
subcategory of $\mathbf D^b_c(X)$ whose objects are the perverse
sheaves.  Given a Whitney stratification, $\mathcal S$, of $X$, it is also useful to
consider the category ${\operatorname{Perv}}_{{}_{\mathcal S}}(X) := {\operatorname{Perv}}(X)\cap \mathbf D^b_{{}_{\mathcal
S}}(X)$ of perverse sheaves which are constructible with respect to $\mathcal S$.

${\operatorname{Perv}}(X)$ and ${\operatorname{Perv}}_{{}_{\mathcal S}}(X)$ are both Abelian categories in which the short exact
sequences $$
0\rightarrow \Adot \rightarrow \Bdot \rightarrow \mathbf
C^\bullet \rightarrow 0
$$
are
precisely the distinguished triangles

\vskip .2in

\vbox{
$\hskip 2.5in \Adot \longrightarrow \ \Bdot$ \vskip .03in
$\hskip 2.35in \ _{[1]} \ \nwarrow\hskip .15in \swarrow$ \vskip .03in $\hskip
2.75in \mathbf C^\bullet$ .
}

\vskip .1in

If we have complexes $\Adot, \Bdot,$ and $\mathbf C^\bullet$ in 
$\mathbf D^b_c(X)$ (resp. $\mathbf D^b_{{}_{\mathcal
S}}(X)$),
a distinguished triangle $\Adot \rightarrow \Bdot
\rightarrow  \mathbf C^\bullet \rightarrow \Adot[1]$, and $\Adot$ and  $\mathbf C^\bullet$ are perverse, then $\Bdot$ is also in
${\operatorname{Perv}}(X)$ (resp. ${\operatorname{Perv}}_{{}_{\mathcal S}}(X)$).

\vskip .2in

Suppose that the base ring, $R$, is a field. Then, ${\operatorname{Perv}}(X)$ is Noetherian and locally Artinian; in particular, locally, each perverse sheaf has a finite composition series in ${\operatorname{Perv}}(X)$ with uniquely determined simple subquotients.  If the Whitney stratification $\mathcal S$ has a finite number of strata, then 
${\operatorname{Perv}}_{{}_{\mathcal S}}(X)$ is  Artinian and Noetherian; if $X$ is compact, then ${\operatorname{Perv}}(X)$ is also Artinian. When $\Pdot$ is a perverse sheaf with a composition series, we let $[\Pdot]$ denote the formal sum of the simple subquotients of the composition series, counted with multiplicities, i.e., we consider $\Pdot$ in the Grothendieck group of perverse sheaves. If $T$ is a morphism from a perverse sheaf $\Pdot$ to itself, then (locally, at least) $[\ker T] = [\operatorname{coker} T]$; in particular, $\operatorname{supp}(\ker T) = \operatorname{supp}(\operatorname{coker} T)$.

Continue to assume that the base ring is a field. Then, the simple objects in ${\operatorname{Perv}}(X)$ (resp. ${\operatorname{Perv}}_{{}_{\mathcal S}}(X)$) are extensions by zero of
intersection cohomology sheaves on irreducible analytic subvarieties (resp. connected components of
strata) of $X$ with coefficients in irreducible local systems.  To be precise, let $M$ be a
$d$-dimensional connected analytic submanifold (resp. a connected component of a stratum) of $X$ and let $\mathcal L_M$ be
an irreducible local system on $M$; then, the pair $(\overline M , \mathcal L_M)$ is called an
{\it irreducible enriched subvariety} of $X$ (where $\overline M$ denotes the
closure of $M$).  Let $j : \overline M \hookrightarrow X$ denote the
inclusion.  Then, the simple objects of ${\operatorname{Perv}}(X)$ (resp. ${\operatorname{Perv}}_{{}_{\mathcal S}}(X)$)
are those of the form $j_!\mathbf{IC}_{{}_{\overline M}}^\bullet(\mathcal L^\bullet_M[d])$, where $(\overline M , \mathcal L_M)$ is 
an irreducible enriched subvariety (again, we are indexing intersection cohomology so that it is non-zero
only in non-positive dimensions).

\medskip

\noindent{\bf Warning}: The stalk cohomology of kernels, cokernels, and images in the perverse category do {\bf not} correspond in any simple way to kernels, cokernels, and images of the stalk cohomology. In fact, the easiest example of this phenomenon is essentially our example from the warning in Section 1.

Recall that we let $X$ be a space consisting of two complex lines $L_1$ and $L_2$ which intersect in a single point $\mathbf p$. For $i=1,2$,
let $\widetilde{\mathbb C}_{L_i}$ denote the $\mathbb C$-constant sheaf on $L_i$ extended by zero to all of $X$. Now, this time, we shift to obtain perverse sheaves. 

There is a canonical
map,
$\alpha$, from the perverse sheaf
$\Bdot:=\mathbb C^\bullet_X[1]$ to the perverse sheaf which is the direct sum of sheaves
$\Cdot:=\widetilde{\mathbb C}^\bullet_{L_1}[1]\oplus\widetilde{\mathbb C}^\bullet_{L_2}[1]$, such that, on
$L_1-\mathbf p$, $\alpha$ is
$\operatorname{id}\oplus 0$, on $L_2-\mathbf p$ is $0\oplus\operatorname{id}$, and is the diagonal map on the stalk at $\mathbf p$. In fact, this map $\alpha$ is the canonical map from the constant sheaf to the intersection cohomology sheaf (with constant coefficients); see Section 5.
As before (essentially), consider the complex, $\Adot$, which has  $\mathbb C_X$ in degree $-1$, $\widetilde{\mathbb C}_{L_1}\oplus\widetilde{\mathbb
C}_{L_2}$ in degree $0$, zeroes elsewhere, and the coboundary map from degree $-1$ to degree $0$ is $\alpha$.  This complex is isomorphic in ${\operatorname{Perv}}(X)$ to the extension by $0$ of the constant sheaf on the point $\mathbf p$. Let $\gamma$ be the morphism of complexes from $\Adot$ to
$\Bdot$ which is the identity in degree $-1$ and is zero elsewhere. Then,
$$
\begin{CD}0\rightarrow\Adot @>\gamma>>\Bdot@>\alpha>>\Cdot\rightarrow 0\end{CD}
$$
is a short exact sequence in ${\operatorname{Perv}}(X)$. 

However, the induced maps, $\gamma^k_x$ and $\alpha^k_x$, on the degree $k$ stalk cohomology, at a point $x\in X$, are zero maps and injections, respectively.

\bigskip

In fact, if your goal is to detect a zero-morphism via some maps of $R$-modules, the situation is even worse than the above indicates. Suppose that $\Pdot$ is a perverse sheaf on $X$ and, for each $\mathbf p\in X$, let $j_{\mathbf p}$ denote the inclusion of $\mathbf p$ into $X$. Let $T$ be a morphism of the perverse sheaf $\Pdot$ to itself. It is possible for $T$ to induce the zero-map on $H^*(j^*_{\mathbf p}\Pdot)$ and on $H^*(j^!_{\mathbf p}\Pdot)$, for all $\mathbf p\in X$, in all degrees, and for $T$ to still not be the zero-morphism in the category of perverse sheaves.

\vskip 0.3in

\noindent Despite the above warning, we do have the following:

\smallskip

Suppose that the base ring is a field, and that $T$ is a morphism from a perverse sheaf $\Pdot$ to itself. Then, $\operatorname{supp}(\ker T)=\overline{\{x\in X\ |\ \ker T^*_x\neq 0\}}$ and, for all $p\in \operatorname{supp}(\ker T)$, if $s:=\dim_p \operatorname{supp}(\ker T)$, then $H^{-s}(\ker T)_p\cong \ker T^{-s}_p$ and, for generic $x\in \operatorname{supp}(\ker T)$ near $p$, $\ker T^j_x=0$ if $j\neq -s$.

\smallskip

The same statements, with each occurrence of $\ker$ replaced by $\operatorname{coker}$, are also true.

\smallskip

Note that these kernel and cokernel statements are {\bf false} for maps between different perverse sheaves.

\vskip 0.3in

We should remark that it is {\bf possible} to give a criterion, in terms of maps of $R$-modules, for a morphism of perverse sheaves to be the zero-morphism; we shall do this later, after we have the vanishing cycles at our disposal.

\vskip .2in

Finally, we wish to state the {\it Decomposition Theorem} of Beilinson, Bernstein, Deligne, and Gabber.  For this statement,
we must restrict ourselves to $R = \mathbb Q$.  We give the statement as it
appears in \cite{Mac1}, except that in \cite{Mac1} intersection cohomology is defined as a shifted perverse sheaf, and we must adjust by shifting back. Note that, in \cite{Mac1}, the setting is algebraic; the
analytic version appears in \cite{Sa}.

An algebraic map $f: X \rightarrow Y$ is called {\it projective} if it can be
factored as an embedding $X \hookrightarrow Y\times \mathbb P^m$ (for some $m$)
followed by projection $Y\times\mathbb P^m \rightarrow Y$.

\vskip .2in

\noindent {\bf The Decomposition Theorem} [BBD, 6.2.5]:  If $X$ is purely $d$-dimensional and $f: X\rightarrow Y$ is
proper, then there exists a unique set of irreducible enriched subvarieties
$\{(\overline M_\alpha , \mathcal L_\alpha)\}$ in $Y$ and Laurent polynomials
$\{\phi^\alpha = \dots + \phi^\alpha_{-2}t^{-2}+\phi^\alpha_{-1}t^{-1}+\phi^\alpha_0 + \phi_1^\alpha t +
\phi^\alpha_2t^2 + \dots\}$  such that there is a quasi-isomorphism
$$Rf_*\mathbf{IC}^\bullet_{{}_X}(\mathbb Q^\bullet_{X^\circ}[d]) \ \cong
\ \bigoplus_{\alpha,i}\ \mathbf{IC}^\bullet_{{}_{{\overline M}_\alpha}}(\mathcal
L^\bullet_\alpha[d_\alpha])[-i]\otimes \mathbb Q^{\phi^\alpha_i} ,$$
(here, $\mathbf{IC}^\bullet_{{}_{{\overline M}_\alpha}}(\mathcal
L^\bullet_\alpha[d_\alpha])$ actually equals ${j_\alpha}_!\mathbf{IC}^\bullet_{{}_{{\overline M}_\alpha}}(\mathcal
L^\bullet_\alpha[d_\alpha])$, where  $j_\alpha : \overline M_\alpha \hookrightarrow Y$ is the
inclusion and, of course, $d_\alpha$ is the dimension of $M_\alpha$).

Moreover, if $f$ is projective, then the coefficients of $\phi^\alpha$ are palindromic around
$0$ (i.e., $\phi^\alpha(t^{-1}) = \phi^\alpha(t))$ and the even and odd terms are separately unimodal (i.e., if $i
\leqslant 0$, then $\phi^\alpha_{i-2} \leqslant \phi^\alpha_i$).

\vskip .1in

Applying hypercohomology to each side, we obtain:
$$IH^k(X; \mathbb Q) = \bigoplus_{\alpha, i}\ (IH^{k-i}(\overline M_\alpha ; \mathcal
L^\bullet_\alpha[d_\alpha]))^{\phi^\alpha_i} .$$

\vskip .2in

We now wish to describe the category of perverse sheaves on a one-dimensional space; this is a particularly nice case of the
results obtained in \cite{M-V}.  Unfortunately, we will use the notions of vanishing cycles and nearby cycles, which are
not covered until the next section.  Nonetheless, it seems appropriate to place this material here.

\vskip .1in

We actually wish to consider perverse sheaves on the germ of a complex analytic space $X$ at a point $x$.  Hence, we assume
that $X$ is a one-dimensional complex analytic space with irreducible analytic
components $X_1, \dots, X_d$ which all contain $x$, such that $X_i$ is homeomorphic to a complex line and $X_i-\{x\}$ is
smooth for all
$i$.  We wish to describe the category, $\mathcal C$, of perverse sheaves on $X$ with complex coefficients which are
constructible with respect to the stratification
$\{X_1-\{x\},
\dots, X_d-\{x\}, \{x\}\}$.

Since perverse sheaves are topological in nature, we may reduce ourselves to considering exactly the case where $X$
consists  of
$d$ complex lines through the origin in some $\mathbb C^N$.  Let $L$ denote a linear form on $\mathbb C^N$ such that $X\cap
L^{-1}(0) = \{\mathbf 0\}$.

Suppose now that $\Pdot$ is in $\mathcal C$, i.e., $\Pdot$ is perverse on $X$ and constructible with respect
to the stratification which has $\{\mathbf 0\}$ as the only zero-dimensional stratum.  Then $\Pdot_{|_{X-\{\mathbf
0\}}}$ consists of a collection of local systems, $\mathcal L_1, \dots, \mathcal L_d$, in degree $-1$.  These local systems are
completely determined by monodromy isomorphisms $h_i:\mathbb C^{r_i}\rightarrow\mathbb C^{r_i}$ representing looping once around
the origin in
$X_i$.  In terms of nearby cycles, the monodromy automorphism on $H^0(\psi_L\Pdot[-1])_\mathbf 0\cong \bigoplus_i \mathbb C^{r_i}$
is given by $\bigoplus_i h_i$.

The vanishing cycles $\phi_L\Pdot[-1]$ are a perverse sheaf on a point, and so have possibly non-zero cohomology only
in degree $0$; say, $H^0(\phi_L\Pdot[-1])_\mathbf 0\cong\mathbb C^\lambda$.  We have the canonical map
$$r:H^0(\psi_L\Pdot[-1])_\mathbf 0\rightarrow
H^0(\phi_L\Pdot[-1])_\mathbf 0$$ and the variation map $$\operatorname{var}: H^0(\phi_L\Pdot[-1])_\mathbf 0\rightarrow
H^0(\psi_L\Pdot[-1])_\mathbf 0,$$ and
$\operatorname{var}\circ r = \operatorname{id}-\bigoplus_i h_i$.  

Thus, an object in $\mathcal C$ determines a vector space $W:=H^0(\phi_L\Pdot[-1])_\mathbf 0$, a vector space  
$V_i:=\mathbb C^{r_i}$ for each irreducible component $X_i$, an automorphism $h_i$ on $V_i$, and two
linear maps $\alpha: \bigoplus_i V_i\rightarrow W$ and $\beta:W\rightarrow \bigoplus_i V_i$ such that
$\beta\circ\alpha =
\operatorname{id}-\bigoplus_i h_i$.  This situation is nicely represented by a commutative triangle

\vbox{
$$
\begin{CD}\oplus_i V_i@>\hskip .05in\operatorname{id}-\bigoplus_i h_i\hskip .05in>>\oplus_i V_i
\end{CD}$$
\vskip  -.1in
$$
\hbox{}\hskip -0.0in\alpha\searrow\hskip .25in\nearrow\beta
$$
$$
\hbox{}\hskip 1in W \hskip 1in .
$$
}
The category $\mathcal C$ is equivalent to the category of such triangles, where a morphism of triangles is defined in the
obvious way: a morphism is determined by linear maps $\tau_i:V_i\rightarrow V_i^\prime$ and \linebreak $\eta: W\rightarrow W^\prime$
such that

\vbox{$$\begin{CD}
\oplus_i V_i @>\ \alpha\ >> W @>\ \beta\ >> \oplus_i V_i \end{CD}
$$
$$
\hbox{}\hskip -.4in\oplus_i\tau_i\downarrow\hskip .6in\eta\downarrow\hskip .4in\oplus_i\tau_i\downarrow
$$
$$
\begin{CD}\oplus_i V_i^\prime @>\ \alpha^\prime\ >> W^\prime @>\ \beta^\prime\ >> \oplus_i V_i^\prime
\end{CD}$$}
commutes.

\vskip .4in

\noindent\S3. {\bf Nearby and Vanishing Cycles}  

\vskip .1in

Before we start in on the topics of this section, we need to state clearly a convention which is common now, but which we did not use when we began these notes. {\bf When dealing with the derived category, ALL functors $F$ are to be considered the derived functors, whether we write $RF$ or merely $F$.} This notation gets mixed with the old notation throughout the remainder of these notes.

\vskip 0.1in

Historically, there was some confusion surrounding the terminology {\it nearby
(or neighboring) cycles} and {\it vanishing cycles}; now, however, the terminology has stabilized.  In the past, the term ``vanishing cycles" was sometimes used to
describe what are now called the nearby cycles (this is true, for instance, in
\cite{A'C}, \cite{BBD}, and \cite{G-M1}.)

The two different indexing schemes for perverse sheaves also add to this confusion in
statements such as ``the nearby cycles of a perverse sheaf are perverse". 
Finally, a new piece of confusion was added in \cite{K-S}, where the sheaf of
vanishing cycles is shifted by one from the usual definition (we will {\bf not} use
this new, shifted definition).

The point is: one should be very careful when reading works on nearby and
vanishing cycles.

\vskip 0.1in

Now, suppose we have a complex analytic space $X$, a function $f:X\rightarrow \C$, and a complex of sheaves $\Fdot\in \mathbf D^b_c(X)$. What we want from nearby (resp., vanishing) cycles is a complex on $V(f)=f^{-1}(0)$ which has stalk cohomology at $x\in V(f)$ which gives the hypercohomology of the Milnor fiber of $f$ at $x$ (resp., of the pair of a small ball at $x$ modulo the Milnor fiber).

Here, when we write the ``Milnor fiber'', we mean at least up to homotopy. This can be defined as follows for $x\in V(f)$. Assume that $X$ is locally embedded in $\C^{n+1}$. Define the Milnor fiber, $F_{f,x}$, of $f$ at $x$ (at angle $\theta$) to be 
$$
F_{f,x}=F_{f, x}^\theta := B^\circ_\epsilon(x)\cap X\cap f^{-1}(re^{i\theta}),
$$
where $0<\epsilon\ll 1$, $B^\circ_\epsilon(x)$ is an open ball of radius $\epsilon$ centered at $x$, $\theta$ is a real number, and $r>0$ is chosen small compared to $\epsilon$; up to homotopy, this is independent of the choices made. In addition, up to isomorphism, $\mathbb H^*(F_{f,x}; \, \Fdot)$ and $\mathbb H^*(B^\circ_\epsilon(x)\cap X, F_{f,x}; \, \Fdot)$ are independent of the choices made.

The nearby and vanishing cycles are easiest and, perhaps, most intuitively, defined by using the Milnor fiber where $\theta$ above is chosen to be $0$, but -- to deal with the Milnor monodromy -- we must discuss the Milnor fiber for other values of $\theta$. However, first we will discuss the traditional, commonly used, definitions.

\vskip .2in

\noindent{\bf The traditional definitions}

\vskip 0.1in

The following is the classic definition of Deligne \cite{De} for the nearby cycles and is the one that is found most often in the literature.

Consider the diagram

\vbox{
$$\begin{CD} E            @>>>         \widetilde{\mathbb C^*}\\ 
@V\hat\pi VV              @VV\pi V\\
X- f^{-1}(0)           @>\hat f>>        \mathbb C^* 
\end{CD}$$

\hbox{}\hskip 2.24in $i \downarrow$
\vskip -.1in
\hbox{}\hskip 1.57in $f^{-1}(0) \xhookrightarrow[ \ j  \ ]{}X$
}

\vskip .1in

\noindent where:

\vskip .1in

$j : f^{-1}(0) \hookrightarrow X$ is inclusion;

\vskip .1in

$i : X - f^{-1}(0) \hookrightarrow X$ is inclusion;

\vskip .1in

$\hat f =$ restriction of $f$;

\vskip .1in

$\widetilde{\mathbb C^*} =\C=\{a+\theta i\,|\, a, \theta\in\R\}$ is the universal (cyclic) cover of $\mathbb C^*$ via the map $\pi(z)=\exp(z)$, 

\vskip 0.1in

and

\vskip .1in

$E$ denotes the pull-back.

\vskip .2in

The {\it nearby (or neighboring) cycles} of $\Fdot$ along $f$ are defined to be
$$
\psi_f\Fdot := j^*R(i\circ\hat \pi)_*(i\circ \hat\pi)^*\Fdot .
$$
Note that this is a complex of sheaves on $f^{-1}(0)$ and that $\psi_f$ is a functor which takes distinguished
triangles to distinguished triangles.

\smallskip

It is not difficult to see that, for $x\in f^{-1}(0)$, the stalk cohomology of $\psi_f\Fdot$ in degree $k$ is isomorphic to the hypercohomology of the Milnor fiber of $f$ at $x$ with coefficients in $\Fdot$, i.e.,
$$
H^k(\psi_f\Fdot)_x\cong \hyp^k(F_{f, x}; \Fdot),
$$
where we use the Milnor fiber at any angle from our discussion above. To see this (using our notation from above), one uses that by constructibility,
$$
H^k(\psi_f\Fdot)_x\cong \hyp^k(B^\circ_\epsilon(x)\cap X\cap f^{-1}(\D^\circ_\delta); R(i\circ\hat \pi)_*(i\circ \hat\pi)^*\Fdot),
$$
where
$0<\epsilon\ll 1$, $B^\circ_\epsilon(x)$ is an open ball of radius $\epsilon$ centered at $x$, and $\D^\circ_\delta$ is an open disk in $\C$, centered at the origin, of radius $\delta>0$ which is small compared to $\epsilon$. (Such a neighborhood of $x$ is often referred to as a {\it Milnor tube}.) Then,
$$\hyp^k(B^\circ_\epsilon(x)\cap X\cap f^{-1}(\D^\circ_\delta); (i\circ\hat \pi)_*(i\circ \hat\pi)^*\Fdot)\cong \hyp^k(F_{f,x}\times (0,\delta)\times \R; \, r^*\Fdot)\cong \hyp^k(F_{f, x}; \Fdot),
$$
where $r$ is obtained from $i\circ \hat\pi$ via the homeomorphism
$$
(i\circ\hat \pi)^{-1}\big(B^\circ_\epsilon(x)\cap X\cap f^{-1}(\D^\circ_\delta)\big)\cong F_{f,x}\times (0,\delta)\times \R.
$$

\medskip

As $\psi_f\big(\Fdot[k]\big)= \big(\psi_f\Fdot\big)[k]$, we may write
$\psi_f\Fdot[k]$ unambiguously. In fact, it is frequently useful to consider the functor where one first shifts the complex by
$k$ and then takes the nearby cycles; thus, we introduce the notation $\psi_f[k]$ to be the functor such that
$\psi_f[k]\Fdot=\psi_f\Fdot[k]$ (and which has the corresponding action on morphisms). 

If $\Pdot$ is a perverse sheaf on $X$, then $\psi_f[-1]\Pdot$ is perverse on $f^{-1}(0)$. (Actually, to conclude that 
$\psi_f[-1]\Pdot$ is  perverse, we only need to assume that
${\Pdot}_{|_{X - f^{-1}(0)}}$ is  perverse.) Because $\psi_f[-1]$ takes perverse
sheaves to perverse sheaves, it is useful to include the shift by $-1$ in many statements about $\psi_f$. Consequently, we also
want to shift $j^*\Fdot$ by $-1$ in many statements, and so we write $j^*[-1]$ for the functor which first shifts by $-1$ and
then pulls-back by $j$.

\medskip

This definition of the nearby cycles is elegant in that it does not refer to a specific angle in providing for a functor which yields the stalk cohomology of the Milnor fiber (at an arbitrary angle).

\bigskip

As there is a canonical map $\Fdot \rightarrow Rg_*g^*\Fdot$ for any
map $g : Z\rightarrow X$, there is a map
$$\Fdot \rightarrow R(i\circ\hat\pi)_*(i\circ\hat\pi)^*\Fdot$$
and, hence, a natural map, $\operatorname{comp}={\operatorname{comp}}_f$, called the {\it comparison map} (or {\it specialization map}):
$$\begin{CD}j^*[-1]\Fdot \ @>\ \operatorname{comp}\ >>
\ j^*[-1]R(i\circ\hat\pi)_*(i\circ\hat\pi)^*\Fdot = \psi_f[-1]\Fdot .\end{CD}$$

\smallskip

There is an automorphism, the {\it monodromy automorphism} 
$T_f:\psi_f[-1]\Fdot\rightarrow \psi_f[-1]\Fdot$, which comes from the deck transformation obtained in Deligne's
definition of $\psi_f\Fdot$ (and, hence, $\psi_f[-1]\Fdot$) by traveling once around the origin in $\mathbb C$. In fact, $T_f$ is
a natural isomorphism from the functor $\psi_f[-1]$ to itself; thus, strictly speaking, when we write
$T_f:\psi_f[-1]\Fdot\rightarrow \psi_f[-1]\Fdot$, we should include $\Fdot$ in the notation for $T_f$ -- however, we shall
normally omit the explicit reference to  $\Fdot$ if the complex is clear. The comparison map 
$$\begin{CD}j^*[-1]\Fdot \ @>\ \operatorname{comp}\ >>\ \psi_f[-1]\Fdot\end{CD}$$ is $T_f$-equivariant, i.e., $\operatorname{comp} = T_f\circ \operatorname{comp}$.

\smallskip

Since we have a map $\operatorname{comp}:j^*[-1]\Fdot \rightarrow \psi_f[-1]\Fdot$, the third vertex of a distinguished triangle is
defined up to isomorphism in the derived category.  We can define the sheaf of shifted {\it vanishing cycles}, $\phi_f[-1]\Fdot$, of $\Fdot$ along $f$ to be this
third vertex, and so, there is a distinguished triangle

\vbox{
\hbox{}$\hskip 1.8in j^*[-1]\Fdot \xrightarrow{\ \textnormal{comp}\ } \ \psi_f[-1]\Fdot$ \vskip .06in
$\hskip 2in\ _{[1]} \ \nwarrow\hskip .4in \swarrow \textnormal{\footnotesize can}$
\vskip .03in
$\hskip 2.3in \phi_f[-1]\Fdot$,
}
\medskip
\noindent which we refer to as the {\it Milnor triangle}.

\medskip

The above is the definition of $\phi_f[-1]\Fdot$ found in many sources.  Unfortunately, this definition of $\phi_f[-1]$ makes it technically difficult to define the functor $\phi_f[-1]$ together with a monodromy automorphism, $\widetilde T_f$, on  $\phi_f[-1]$ which is compatible with the monodromy on $\psi_f[-1]$; see the discussion of Sch\"urmann on pages 25-26 of \cite{Sch}, where he discusses how to fix this issue. When this issue is addressed, the map $\psi_f[-1]\xrightarrow{\ \textnormal{can}\ }\phi_f[-1]$ in the triangle is a natural map of functors; it is called the {\it canonical map}, and there exists a natural automorphism of the Milnor triangle (functors):

\vbox{$$\begin{CD}
 j^*[-1]@>\ \textnormal{comp} \ >> \psi_f[-1] @>\ \textnormal{can} \ >> \phi_f[-1]\xrightarrow{[1]} \end{CD}
$$
$$
\hbox{}\hskip -.35in{\scriptstyle\textnormal{id}}\ \downarrow\hskip .6in\ \ {\scriptstyle T_f} \ \downarrow\hskip .65in\ {\scriptstyle \widetilde T_f}\downarrow
$$
$$
\begin{CD}j^*[-1]@>\ \textnormal{comp} \ >> \psi_f[-1] @>\ \textnormal{can} \ >> \phi_f[-1]\xrightarrow{[1]} \end{CD}
$$}

\bigskip

There is another difficulty with the definition of $\phi_f[-1]\Fdot$ as the mapping cone of the comparison map; it makes it difficult to define the {\it variation morphism}, var, from $\phi_f[-1]$ to $\psi_f[-1]$, and to derive the distinguished triangle of functors which we refer to as the {\it variation triangle} (or {\it dual Milnor triangle}):

\vbox{
\hbox{}$\hskip 1.8in \phi_f[-1]\Fdot \xrightarrow{\ \textnormal{var}\ } \ \psi_f[-1]\Fdot$ \vskip .06in
$\hskip 2in\ _{[1]} \ \nwarrow\hskip .4in \swarrow  \textnormal{\footnotesize pmoc}$
\vskip .03in
$\hskip 2.4in j^![1]\Fdot$,
}

\medskip

\noindent where, for lack of a better morphism name, we have used pmoc (i.e., comp written backwards).

On pages 369-370 of \cite{Sch}, Sch\"urmann also defines the variation map, gives the proof that the variation triangle above exists, and shows that $\operatorname{var}\circ\operatorname{can}=\operatorname{id}-T_f$. Furthermore, there exists a natural automorphism of the variation triangle (functors):

\vbox{$$\begin{CD}
\phi_f[-1]@>\ \textnormal{var} \ >> \psi_f[-1] @>\ \textnormal{pmoc} \ >>  j^![1]\xrightarrow{[1]} \end{CD}
$$
$$
\hbox{}\hskip -.3in{\scriptstyle \widetilde T_f}\ \downarrow\hskip .55in\ \ {\scriptstyle T_f} \ \downarrow\hskip .65in\ {\scriptstyle\textnormal{id}}\downarrow
$$
$$
\begin{CD}
\phi_f[-1]@>\ \textnormal{var} \ >> \psi_f[-1] @>\ \textnormal{pmoc} \ >>  j^![1]\xrightarrow{[1]} \end{CD}
$$}

\bigskip

\noindent {\bf The Kashiwara and Schapira approach}

\medskip

 In Definition 8.6.2, Kashiwara and Schapira (K-S) \cite{K-S} give an elegant definition of $\phi_f[-1]$ (but they call it $\phi_f$) which allows them to easily define the monodromy automorphism and the variation map and prove fundamental results, but that definition is a bit hard to understand. Still, we suggest the reader make the effort to follow the K-S approach; it is very beautiful. (Their derivation uses Poincar\'e-Verdier duality and so does not work in the more general setting of Sch\"urmann; see \cite{Sch}, Remark 5.4.7.) Throughout our discussion of the K-S approach, we will replace their $\phi_f$ with $\phi_f[-1]$ to agree with all other sources. We also want to indicate some small issues with their diagram 8.6.5.

The idea in the K-S approach is fairly simple; they want the vanishing cycles to be isomorphic to the mapping cone of the comparison map in a canonical, natural way - a way that lets them easily define the variation map, demonstrate the compatibility of the monodromy automorphisms, and show that 
$$\operatorname{var}\circ\operatorname{can}=\operatorname{id}-T_f \hskip 0.2in \textnormal{and} \hskip 0.2in \operatorname{can}\circ\operatorname{var}=\operatorname{id}-\widetilde T_f.
$$
 
Kashiwara and Schapira begin with the standard Deligne definition of the nearby cycles (with mild changes). They then use Poincar\'e-Verdier duality to
show that the nearby cycles $\psi_f\Fdot$ are isomorphic to the restriction to $f^{-1}(0)$ of  $R\mathbf {Hom}^\bullet(f^*p_!{\mathbf A}_{\widetilde{\C^*}},\Fdot)$; see their formula 8.6.2.

Then, they define a complex $\Kdot$ with only two non-zero terms {\bf precisely} so that the restriction to $f^{-1}(0)$ of  $R\mathbf {Hom}^\bullet(f^*\Kdot, \Fdot)$  yields a complex isomorphic to the mapping cone of the comparison map.

For the reader attempting to follow the K-S construction/derivation, we want to clarify some issues in their diagram 8.6.5. First, we let $\hat p$ denote the covering projection $\hat p: \widetilde{\C^*}\rightarrow \C^*$ and let $u:\C^*\hookrightarrow\C$ denote the inclusion, so that the K-S map $p$ factors as $p=u\hat p$. We let $\widehat{\operatorname{tr}}$ and $w$ denote the natural maps 
$$\widehat{\operatorname{tr}}: \hat p_!\hat p^!{\mathbf A}_{\C^*}\rightarrow{\mathbf A}_{\C^*} \hskip 0.2in \textnormal{and}\hskip 0.2in w:u_!u^!{\mathbf A}_{\C}\rightarrow {\mathbf A}_{\C}.
$$
With this notation, we have an isomorphism $p_!p^!{\mathbf A}_\C\cong u_!\hat p_!\hat p^!{\mathbf A}_{\C^*}$ and the trace morphism, $\operatorname{tr}$, of K-S is naturally isomorphic to $wu_!(\widehat{\operatorname{tr}})$. Now, in diagram 8.6.5 of \cite{K-S}, there are two mild errors in the top row which we can ``fix''. The last non-zero entry in the top row should be $u_!{\mathbf A}_{\C^*}$ and the map written as $\operatorname{tr}$ in the top row should be $u_!(\widehat{\operatorname{tr}})$. Finally, the last downward arrow (which is unlabeled) should be our map $w$.

Kashiwara and Schapira give a more down-to-earth characterization of the vanishing cycle functor in Exercise VIII.13:
$$
\phi_f[-1]\Fdot \cong \big(R\Gamma_{\operatorname{Re}f\geq 0}(\Fdot)\big)_{{\big |}_{f^{-1}(0)}},
$$ 
where we have included a shift to make this characterization agree with all other sources. Note that, like our earlier discussion of the Milnor fiber, this characterization uses angles in $\C^*$ in the guise of ${\operatorname{Re}f\geq 0}$. This exercise is not difficult since it is relatively easy to show that $\big(R\Gamma_{\operatorname{Re}f\geq 0}(\Fdot)\big)_{{\big |}_{f^{-1}(0)}}$ is isomorphic to the mapping cone of the comparison map and so is isomorphic to $\phi_f[-1]\Fdot$ in the derived category.

\bigskip
 
 \noindent{\bf The ``intuitive'' definitions}
 
 \smallskip
 
One can define the nearby and vanishing cycles functors in a way that corresponds more obviously to how the Milnor fiber is described topologically. 
 
 Let $L:=f^{-1}([0, \infty))$, let $\ell:f^{-1}(0)\hookrightarrow L$  and $m:L\hookrightarrow X$ be the closed inclusions, and let $n:L-f^{-1}(0)\hookrightarrow L$ be the open inclusion. Note that $j=m\ell$.
 
 Then we have natural isomorphisms 
 $$\psi_f[-1]\simeq\ell^*n_*n^*m^*[-1]\hskip 0.2in \textnormal{ and }\hskip 0.2in  \phi_f[-1]\simeq \ell^*\ell_!\ell^!m^*\simeq \ell^!m^*.$$
 With these definitions, it is easy to see that, for all $x\in V(f)$ and using $\theta=0$ in the definition of the Milnor fiber $F^\theta_{f, x}$, we have isomorphisms for the stalk cohomologies given by
  $$
 H^k(\psi_f[-1]\Fdot)_x\cong \mathbb H^{k-1}(F^0_{f,x};\,\Fdot)\hskip 0.2in\textnormal{and}\hskip 0.2in H^k(\phi_f[-1]\Fdot)_x\cong \mathbb H^{k}(B_\epsilon^\circ(x)\cap X, F^0_{f,x};\,\Fdot).
 $$
 
 \bigskip

Beginning with the canonical distinguished triangle
$$
\rightarrow \ell_!\ell^!\rightarrow \operatorname{id}\rightarrow n_*n^*\xrightarrow{[1]},
$$
we find there is a distinguished triangle
$$
\rightarrow \ell^*\ell_!\ell^!m^*\rightarrow \ell^*m^*\rightarrow \ell^*n_*n^*m^*\xrightarrow{[1]}.
$$
Turning the above triangle, and using our definitions of the nearby and vanishing cycles and that $j=m\ell$, we arrive at the Milnor triangle:

$$\begin{CD}j^*[-1]@>\ {\operatorname{comp}} \ >> \psi_f[-1]@>\ {\text{can}}\ >>\phi_f[-1]\xrightarrow{[1]}\end{CD}.$$

\bigskip

It is an easy exercise to show that $\psi_f\Fdot$ depends only on ${\Fdot}_{|_{X - f^{-1}(0)}}$ (and $f$), and so it is utterly trivial that the natural maps $\alpha: \operatorname{id}\rightarrow i_*i^*$ and $\beta: i_!i^!\rightarrow\operatorname{id}$ induce natural isomorphisms
$$
\psi_f[-1]\xrightarrow[\ \simeq  \ ]{\ \gamma:=\psi_f[-1]\alpha\ }\psi_f[-1]i_*i^* \hskip 0.2in \textnormal{ and }\hskip 0.2in \psi_f[-1]i_!i^! \xrightarrow[\ \simeq \ ]{\ \delta :=\psi_f[-1]\beta\ }\psi_f[-1].
$$

\bigskip

Applying the Milnor triangle to $i_!i^!$ and noting that $j^*i_!i^!=0$, we find that there is a natural isomorphism induced by ${\text{can}}$:
$$
\psi_f[-1]i_!i^!\xrightarrow[\ \simeq \ ]{\tau}\phi_f[-1]i_!i^!,
$$
and so a natural isomorphism $\nu$:
$$
\psi_f[-1]\xrightarrow[\simeq]{\nu=\tau\delta^{-1} }\phi_f[-1]i_!i^!.
$$

By applying the shifted vanishing cycle functor to the canonical morphism $\beta: i_!i^!\rightarrow \operatorname{id}$, we obtain a natural map 
$$\phi_f[-1]i_!i^!\xrightarrow[]{\rho:=\phi_f[-1]\beta}\phi_f[-1]$$
such that ${\operatorname{can}}=\rho\circ \nu$.

\bigskip

\vskip 0.2in

\vbox{
\noindent{\bf Nearby and vanishing cycles at arbitrary angles}

\vskip .2in

To treat the Milnor monodromy automorphisms on $\psi_f[-1]$ and $\phi_f[-1]$, we must repeat the above definitions of $\psi_f[-1]$ and $\phi_f[-1]$ but now at arbitrary angles $\theta$. We also need this discussion in order to obtain a natural isomorphism $\omega$, the {\bf Wang isomorphism}:}
$$
\phi_f[-1]i_*i^*\xrightarrow[ \ \simeq \ ]{\omega}\psi_f[-1],
$$
which is crucial in our definition of the {\it variation map}.

\vskip 0.2in

We continue to let $j$ denote the inclusion of $f^{-1}(0)$ into $X$.

\medskip

Let 
$$B_\theta:=f^{-1}\left(e^{i\theta}\big\{v\in\mathbb C \ | \ \operatorname{Re} v\geq 0\big\}\right),$$
and let $p_\theta$ denote the inclusion of $B_\theta$ into $X$. Let 
$$A_\theta:=f^{-1}\left(e^{i\theta}\big\{v\in\mathbb C \ | \ \operatorname{Re} v> 0\big\}\right),$$
and let $q_\theta$ denote the inclusion of $A_\theta$ into $X$. Note that $A_\theta = X-B_{\theta+\pi}$.

Let 
$$L_\theta:=f^{-1}\left(e^{i\theta}[0, \infty)\right),
$$
and let $l_\theta: f^{-1}(0)\hookrightarrow L_\theta$, $m_\theta: L_\theta\hookrightarrow X$, and $n_\theta: L_\theta-f^{-1}(0)\hookrightarrow L_\theta$ denote the inclusions. Note that $j=m_\theta l_\theta$.

Then, the shifted {\it nearby cycles and vanishing cycles of $\Fdot$ along $f$ at angle $\theta$} are defined by
$$
\psi_f^\theta[-1]\Fdot \ := \ j^*{q_{\theta}}_*q_{\theta}^*\Fdot[-1] \ \cong \ l_\theta^*{n_\theta}_*n_\theta^*m_\theta^*\Fdot[-1]
$$
and
$$
\phi_f^\theta[-1]\Fdot \ := \ j^*{p_{\theta+\pi}}_!p_{\theta+\pi}^!\Fdot \ \cong \  \big(R\Gamma_{{}_{B_{\theta+\pi}}}(\Fdot)\big)_{|_{f^{-1}(0)}} \ \cong \ l_\theta^*{l_\theta}_!l_\theta^!m_\theta^*\Fdot  \ \cong \ l_\theta^!m_\theta^*\Fdot.
$$
respectively.

\smallskip

Then, one can define $\psi_f[-1]\Fdot := \psi^0_f[-1]\Fdot$ and $\phi_f[-1]\Fdot:=\phi^0_f[-1]\Fdot$, and consider the natural isomorphisms $T^\theta_f: \psi_f[-1]\Fdot \rightarrow \psi_f^\theta[-1]\Fdot$ and $\widetilde T^\theta_f: \phi_f[-1]\Fdot \rightarrow \phi_f^\theta[-1]\Fdot$ which correspond to ``rotating'' $B_\theta$ or $L_\theta$ an angle of $\theta$ counterclockwise around the origin over $\mathbb C$. Then, $T_f:=T^{2\pi}_f$ and $\widetilde T_f:=\widetilde T^{2\pi}_f$ are the usual monodromy actions on the shifted nearby and vanishing cycles.

\medskip

In this set-up, the comparison map $j^*[-1]\rightarrow \psi_f[-1]$ is induced by the map $\operatorname{id}\rightarrow {q_0}_*{q_0}^*$ and/or $\operatorname{id}\rightarrow {n_0}_*{n_0}^*$, while the canonical map $\psi_f[-1]\rightarrow \phi_f[-1]$ is induced by ${q_0}_*q_0^*[-1]\rightarrow {p_\pi}_!p_\pi^!$ and/or ${n_0}_*n_0^*[-1]\rightarrow {l_0}_!l_0^!$.

\bigskip

Now it is time to produce the variation map. By our construction, one clearly has the automorphism of distinguished triangles:

\smallskip

\vbox{$$\begin{CD}
j^*[-1]@>\ \textnormal{comp} \ >> \psi_f[-1] @>\ \textnormal{can} \ >> \phi_f[-1]\xrightarrow{[1]} \end{CD}
$$
$$
\hbox{}\hskip -.35in{\scriptstyle\textnormal{id}}\ \downarrow\hskip .6in\ \ {\scriptstyle T_f} \ \downarrow\hskip .65in\ {\scriptstyle \widetilde T_f}\downarrow
$$
$$
\begin{CD}j^*[-1]@>\ \textnormal{comp} \ >> \psi_f[-1] @>\ \textnormal{can} \ >> \phi_f[-1]\xrightarrow{[1]} \end{CD}
$$}

\bigskip

Now, considering the commutative diagram between two distinguished triangles

\vbox{$$\begin{CD}
j^*[-1]\Fdot@>\ \textnormal{comp} \ >> \psi_f[-1]\Fdot @>\ \textnormal{can} \ >> \phi_f[-1]\Fdot\xrightarrow{[1]} \end{CD}
$$
$$
\hbox{}\hskip -0.95in{\scriptstyle\textnormal{id}}\ \downarrow\hskip .48in\ \ {\scriptstyle \operatorname{id}-T_f} \ \downarrow\hskip .65in
$$
$$
\begin{CD}\hbox{}\hskip 0.32in 0\phantom{v}@>\  \ \phantom{comp}\  \ \ >>\psi_f[-1]\Fdot @>\ \textnormal{id} \ >> \psi_f[-1]\Fdot\xrightarrow{[1]} \end{CD}\ ,
$$}

\medskip

\noindent we see that there exists a map, which we {\bf could} call the variation, $\operatorname{var}:\phi_f[-1]\Fdot \rightarrow \psi_f[-1]\Fdot$ which makes the diagram commute, i.e., such that $\operatorname{var}\circ\operatorname{can}=\operatorname{id} -T_f$. However, this has the usual problem when completing a commutative diagram between distinguished triangles; it would not uniquely define $\operatorname{var}$ and would not yield a natural map of functors which is compatible with the two monodromies such that $\operatorname{var}\circ\operatorname{can}=\operatorname{id} -T_f$ and $\operatorname{can}\circ\operatorname{var}=\operatorname{id} -\widetilde T_f$. We know from Kashiawara-Schapira and Sch\"urmann that such a natural transformation exists, but we want to define it using our ``simple'' definitions of the nearby and vanishing cycles (at angle $0$). And so, we must proceed more carefully.

\bigskip

\noindent{\bf Construction of the Wang isomorphism}

\medskip

Our primary tool for defining the variation will be the Wang isomorphism, which we construct now using functor-level excision.

\medskip

Recall that $A_\pi=f^{-1}\left(\big\{v\in\mathbb C \ | \ \operatorname{Re} v<0\big\}\right)$
and that $q_\pi$ is the (open) inclusion of $A_\pi$ into $X$. Also recall that $L_\pi:=f^{-1}\left( (-\infty,0])\right)$ and  $m_\pi: L_\pi\hookrightarrow X$ is the (closed) inclusion. 

We introduce  the (closed) inclusion $\rho:L_\pi-f^{-1}(0)\hookrightarrow A_\pi$, the (open) inclusion $\sigma: A_\pi\hookrightarrow X-f^{-1}(0)$, and the (open) inclusion $\tau: A_\pi-L_\pi\hookrightarrow A_\pi$. Note that $\sigma\rho$ is the inclusion of a closed set, $q_\pi=i\sigma$, and that $A_\pi-L_\pi$ has two connected components. The functor-level excision that we referred to is the natural isomorphism $(\sigma\rho)_!(\sigma\rho)^!\cong \sigma_*\rho_!\rho^!\sigma^*$.

It is relatively easy to show that we have natural isomorphisms
$$\phi_f[-1]i_*i^*\cong j^*{m_\pi}_!{m_\pi}^!i_*i^*\cong j^*i_*(\sigma\rho)_!(\sigma\rho)^!i^*\cong j^*i_*\sigma_*\rho_!\rho^!\sigma^*i^*\cong j^*{q_\pi}_*\rho_!\rho^!{q_\pi}^*.
$$

Now, the distinguished triangle 
$$
\rho_!\rho^!\rightarrow\operatorname{id}\rightarrow \tau_*\tau^*\xrightarrow{[1]}
$$
yields a distinguished triangle
$$
\phi_f[-1]i_*i^*\cong j^*{q_\pi}_*\rho_!\rho^!{q_\pi}^*\rightarrow j^*{q_\pi}_*{q_\pi}^*\rightarrow j^*{q_\pi}_*\tau_*\tau^*{q_\pi}^*\xrightarrow{[1]}.
$$

\smallskip

\noindent Shifting, rotating, and using that $j^*{q_\pi}_*{q_\pi}^*[-1]\rightarrow j^*{q_\pi}_*\tau_*\tau^*{q_\pi}^*[-1]$ is naturally isomorphic to 
$$\psi_f[-1]\xrightarrow{(\operatorname{id}, T_f)} \psi_f[-1]\oplus \psi_f[-1],
$$

\medskip

\noindent we arrive at a natural map $\eta: \psi_f[-1]\oplus \psi_f[-1]\rightarrow \phi_f[-1]i_*i^*$ and a distinguished triangle

$$
\psi_f[-1]\xrightarrow{(\operatorname{id}, T_f)} \psi_f[-1]\oplus \psi_f[-1]\xrightarrow{\eta} \phi_f[-1]i_*i^*\xrightarrow{[1]}.
$$
\vskip 0.2in

Below, we use notation which looks like we are assuming that the complexes of sheaves have ``elements'' on which there is a group structure; we are {\bf not} assuming this, and the translations into the correct statements using the group structure on the morphisms (i.e., using the additive structure on the triangulated category) is trivial. We write things in this manner because we feel that they are easier to read and understand.

With this understanding, we proceed. Since  $\eta\circ (\operatorname{id}, T_f)=0$, we have that 
$$
\eta(x,y)=\eta(x,y)-\eta(x, T_f x)=\eta(0, y-T_f x)
$$
and so there is a morphism of distinguished triangles

\vbox{$$\psi_f[-1]\xrightarrow{ \ (\operatorname{id}, T_f) \ }\psi_f[-1]\oplus\psi_f[-1]\xrightarrow{(x,y)\ \mapsto \ y-T_fx}  \psi_f[-1]\xrightarrow{[1]}$$
$$
\hbox{}\hskip -.45in{\scriptstyle\textnormal{id}}\ \downarrow\hskip 1in\ \ {\scriptstyle\textnormal{id}} \ \downarrow \  \hskip 1in\ {\scriptstyle z\mapsto \eta(0, z)}\downarrow
$$
$$
\hbox{}\hskip 0.2in\psi_f[-1]\xrightarrow{\ (\operatorname{id}, T_f)\ }\psi_f[-1]\oplus\psi_f[-1]\xrightarrow{\ \phantom{abcde}\eta\phantom{efghi} \ }  \phi_f[-1]i_*i^*\xrightarrow{[1]}.
$$}

\smallskip

\noindent It follows that the natural map $\zeta: \psi_f[-1]\rightarrow \phi_f[-1]i_*i^*$ given by $\zeta(z)=\eta(0, z)$ is an isomorphism. 

Recall that we use $\alpha: \operatorname{id}\rightarrow i_*i^*$ to denote the canonical natural map. We define the {\bf Wang natural isomorphism $\omega$} from $\phi_f[-1]i_*i^*$ to $\psi_f[-1]$  to be $\zeta^{-1}$ and
we define the {\bf variation natural transformation}, $\operatorname{var}$, by
$$
\operatorname{var}= \omega\circ\phi_f[-1](\alpha): \phi_f[-1]\rightarrow \psi_f[-1].
$$

\smallskip

\noindent It follows from the construction that
$$
\operatorname{var}\circ \widetilde T_f= T_f\circ\operatorname{var}, \ \ \operatorname{id}-T_f= \operatorname{var}\circ \operatorname{can}, \textnormal{ and } \operatorname{id}-\widetilde T_f= \operatorname{can}\circ \operatorname{var}.
$$

\bigskip

Recall now our earlier (simple) natural isomorphism:
$$
\psi_f[-1]\xrightarrow[\simeq]{\ \nu\ }\phi_f[-1]i_!i^!.
$$

The {\bf Milnor triangle} follows from applying $\phi_f[-1]$ to the canonical distinguished triangle 
$$
i_!i^!\rightarrow\operatorname{id}\rightarrow j_*j^*\xrightarrow{[1]}
$$
to obtain
$$
\phi_f[-1]i_!i^!\rightarrow\phi_f[-1]\rightarrow \phi_f[-1]j_*j^*\xrightarrow{[1]},
$$
then turning the triangle and using the isomorphisms $\nu^{-1}$ and $\phi_f[-1]j_*j^*\cong j^*$ to produce (again)
$$
j^*[-1]\xrightarrow{ \ \operatorname{comp} \ } \psi_f[-1]\xrightarrow{\text{ \ can \ }}\phi_f[-1]\xrightarrow{[1]}.
$$

\bigskip

In an analogous fashion, the {\bf variation triangle} follows from applying $\phi_f[-1]$ to the canonical distinguished triangle 
$$
j_!j^!\rightarrow\operatorname{id}\rightarrow i_*i^*\xrightarrow{[1]}
$$
and using the Wang isomorphism to obtain
$$
\phi_f[-1]\xrightarrow{\ \textnormal{var} \ }\psi_f[-1] \xrightarrow{\ \textnormal{pmoc} \ }j^![1]\xrightarrow{[1]}.
$$

\medskip

Furthermore, applying $\phi_f[-1]$ to the canonical distinguished triangle
$$
i_!i^!\rightarrow i_*i^*\rightarrow j_*j^*i_*i^*\xrightarrow{[1]}
$$
and using the isomorphisms $\nu^{-1}$ and $\omega$, we obtain the {\bf Wang triangle}
$$
j^*[-1]i_*i^*\rightarrow\psi_f[-1]\xrightarrow{\operatorname{id}-T_f}\psi_f[-1]\xrightarrow{[1]}.
$$

\smallskip

\noindent The associated long exact sequences on stalk cohomology are the Wang sequences for the Milnor fibrations.

\vskip 0.2in

Using the natural isomorphism $j^!i_![1]\cong j^*i_*$, we also have the distinguished triangle 

$$
\psi_f[-1]\xrightarrow{\operatorname{id}-T_f}\psi_f[-1]\rightarrow j^![1]i_!i^!\xrightarrow{[1]}.
$$

\smallskip

\noindent The associated long exact sequences on costalk cohomology are the Wang sequences for the Milnor fibrations with boundary.

\vskip .2in

From the constructions, it is clear that taking nearby cycles and vanishing cycles are local operations. That is, if $l:W\hookrightarrow X$ is the inclusion of an open set, $\hat l:W\cap f^{-1}(0)\hookrightarrow f^{-1}(0)$ is the corresponding inclusion, and $\hat f :=f_{|_{W}}$, then there are natural isomorphisms
$$
{\hat l}^*\psi_f \ \cong \ \psi_{\hat f}\,l^*\hskip .2in \text{and}\hskip .2in {\hat l}^*\phi_f \ \cong \ \phi_{\hat f}\,l^*.
$$

\medskip

Let $\mathcal S = \{S_\alpha\}$ be a Whitney stratification of $X$ and suppose $\Fdot \in \mathbf D^b_{{}_{\mathcal S}}(X)$.  Given an analytic map $f: X
\rightarrow \mathbb C$, define a {\it (stratified) critical point} of $f$ (with
respect to $\mathcal S$) to be a point $x \in S_\alpha \subseteq X$ such that
$f_{|_{S_\alpha}}$ has a critical point at $x$; we denote the set of such critical
points by  $\Sigma_{{}_{\mathcal S}}f$.

For any Whitney stratification, $\mathcal S$, with respect to which $\Fdot$ is
constructible, the support of $\mathbf H^*(\phi_f\Fdot)$ is contained in the stratified
critical locus of $f$, $\Sigma_{{}_{\mathcal S}}f$. In addition, if  $\mathcal S$ is a Whitney stratification with respect to
which $\Fdot$ is constructible and such that $f^{-1}(0)$ is a union of strata, then --  by \cite{BMM} and \cite{P} -- it
follows that $\mathcal S$ also satisfies Thom's $a_f$ condition; by Thom's second isotopy lemma, this implies that the entire
situation locally trivializes over strata, and hence both $\psi_f\Fdot$ and $\phi_f\Fdot$ are constructible with respect to
$\{S\in\mathcal S \ |\ S\subseteq f^{-1}(0)\}$.
\vskip .1in

If $\Pdot$ is a perverse sheaf on $X$, then $\phi_f[-1]\Pdot$ is a perverse sheaf on $f^{-1}(0)$.

\vskip .2in

The reader may wonder how the variation morphism is related to the classical variation isomorphism for isolated affine hypersurface singularities (see, for instance \cite{Lo}). Let $m_x$ denote the inclusion of $x$ into $f^{-1}(0)$, and let $l_x$ denote the inclusion of $x$ into $X$, so that $l_x=j\circ m_x$. Then, the variation morphism yields a natural map from $m_x^!\phi_f[-1]\Fdot$ to $m_x^!\psi_f[-1]\Fdot$, and so gives corresponding maps on cohomology; these maps are, in a sense, the generalization of the classical variation maps.

In particular, suppose that $x$ is an isolated point in the support of $\phi_f[-1]\Fdot$. Then, there is a natural isomorphism 
$$\begin{CD}\gamma_x:m_x^*\phi_f[-1]\Fdot@>\cong>> m_x^!\phi_f[-1]\Fdot.\end{CD}$$ By composing with $$m_x^*(\operatorname{can}): m^*_x\psi_f[-1]\Fdot\rightarrow m^*_x\phi_f[-1]\Fdot$$  and 
$$m_x^!(\operatorname{var}): m_x^!\phi_f[-1]\Fdot\rightarrow m_x^!\psi_f[-1]\Fdot,$$ we obtain a morphism
$$
v_x:=m_x^!(\operatorname{var})\circ \gamma_x\circ m_x^*(\operatorname{can}): m^*_x\psi_f[-1]\Fdot\rightarrow m_x^!\psi_f[-1]\Fdot.
$$
Using a local embedding into affine space, a closed ball $B_\epsilon(x)$, centered at $x$, of radius $\epsilon$, and assuming $0< |\xi|\ll\epsilon\ll<1$, we obtain induced maps on cohomology
$$
\begin{CD} \mathbb H^{i-1}(B_\epsilon(x)\cap X \cap f^{-1}(\xi) ; \Fdot) @>v^i_x>>  \mathbb H^{i-1}(B_\epsilon(x)\cap X\cap  f^{-1}(\xi), \partial B_\epsilon(x)\cap X \cap f^{-1}(\xi) ; \Fdot).\end{CD}
$$
Writing $F_{f, x}$ for the (compact) Milnor fiber  $B_\epsilon(x)\cap X \cap f^{-1}(\xi)$, this gives us
$$
\begin{CD}\mathbb H^{i-1}(F_{f, x}; \Fdot) @>v^i_x>>  \mathbb H^{i-1}(F_{f, x}, \partial F_{f, x}; \Fdot).\end{CD}
$$

If $l_x^![1]\Fdot\cong m_x^!j^![1]\Fdot$ is zero in degree $i$ and degree $i-1$, then $m_x^!(\operatorname{var})$ is an isomorphism in degree $i$. If $l_x^*[-1]\Fdot\cong m_x^*j^*[-1]\Fdot$ is zero in degree $i$ and $i+1$, then $m_x^*(\operatorname{can})$ is an isomorphism in degree $i$. Thus, when the conditions of the previous two sentences are satisfied, $v_x^i$ is an isomorphism.

In the classical case, $X$ is an open subset of $\mathbb C^{n+1}$, $\Fdot$ is the constant sheaf $\mathbb Z_X$, and one looks at $v^{n+1}_x$, where $n\geq 1$. In this case,  all of the maps composed to yield $v_x:=m_x^!(\operatorname{var})\circ \gamma_x\circ m_x^*(\operatorname{can})$ are isomorphisms in degree $n+1$. $v_x^{n+1}$ is the classical variation isomorphism on cohomology.

\vskip .2in

One might wonder if there is a distinguished triangle analogous to the Wang triangle, but involving the vanishing cycles and $ \operatorname{id}-\widetilde{T}_f$. 

To answer this, start with the two distinguished triangles
$$
\begin{CD}\phi_f[-1]\Fdot@>\ \operatorname{var}\ >>\psi_f[-1]\Fdot@>\ \operatorname{pmoc}\ >> j^![1]\Fdot\longrightarrow\phi_f\Fdot\end{CD}
$$
and
$$
\begin{CD}\psi_f[-1]\Fdot@>\ \  \operatorname{can}\ \ >>\phi_f[-1]\Fdot\longrightarrow j^*\Fdot@>\ \operatorname{comp}\ >>\psi_f\Fdot,\end{CD}
$$

\medskip

\noindent and apply the octahedral lemma to conclude that there exists a complex $w_f\Fdot$ and two distinguished triangles 
$$
\phi_f[-1]\Fdot\xrightarrow{\ \operatorname{id}-\widetilde{T}_f\ }\phi_f[-1]\Fdot\longrightarrow w_f\Fdot\xrightarrow{[1]}
$$
and
$$
j^*[-1]\Fdot\xrightarrow{\tau_f:=\operatorname{pmoc}\circ\operatorname{comp}}j^![1]\rightarrow w_f\Fdot\xrightarrow{[1]}.
$$

\bigskip

Thus the application of the octahedral lemma tells us that the mapping cone of $\operatorname{id}-\widetilde{T}_f$ is isomorphic to
the mapping cone of $\tau_f$. Note that, while $j^*[-1]\Fdot$ and $j^![1]\Fdot$ depend only on $f^{-1}(0)$ (and $\Fdot$),
$\tau_f$ may change if $f$ (or some factor of $f$) is raised to a power.  We refer to the morphism $\tau_f$ from $j^*[-1]\Fdot$ to $j^![1]\Fdot$ as the {\it Wang morphism of $f$}. 

\vskip .1in

Also note that the Wang morphism is {\bf not} (usually) the natural morphism from $j^*[-1]\Fdot$ to $j^![1]\Fdot$ obtained by restricting the composition of natural maps
$$
j_*j^*[-1]\Fdot\longrightarrow Ri_!i^!\Fdot\longrightarrow Ri_*i^*\Fdot\longrightarrow j_!j^![1]\Fdot.
$$
The map above induces the zero map on stalk cohomology at all points in $f^{-1}(0)$; the Wang morphism induces isomorphisms on the stalk cohomology at all points of $f^{-1}(0)-\operatorname{supp}\phi_f[-1]\Fdot$.

\vskip 0.2in

Let $\Pdot$ be a perverse sheaf on a $d$-dimensional space $X$. Suppose that $f:X\rightarrow\mathbb C$ is such that $Y:=f^{-1}(0)$ contains no $d$-dimensional irreducible component of $X$. As above, let $j:Y\hookrightarrow X$ and $i: X-Y\hookrightarrow X$ denote the inclusions. Then $\phi_f[-1]\Pdot$ and $\psi_f[-1]\Pdot$ are perverse on a space of dimension $d-1$. Consequently, for all $k\leq -d$, $\mathbf H^k(\phi_f[-1]\Pdot)=0$ and $\mathbf H^k(\psi_f[-1]\Pdot)=0$ and so, by the variation triangle, $\mathbf H^k(j^![1]\Pdot)=0$ for all $k\leq -d-1$, i.e., $\mathbf H^k(j^!\Pdot)=0$ for all $k\leq -d$. It follows that $\mathbb H^k(X; j_!j^!\Pdot)=0$ for all $k\leq -d$. From the distinguished triangle relating $j_!j^!$ and $i_*i^*$, we find that the natural map 
$$\mathbb H^{-d}(X; \Pdot)\to\mathbb H^{-d}(X; i_*i^*\Pdot)\cong \mathbb H^{-d}(X-Y; \Pdot_{|_{X-Y}})\cong H^{0}\big(X-Y; \mathbf H^{-d}\big(\Pdot_{|_{X-Y}}\big)\big)$$
 is an injection.

If $f$ is chosen so that $\mathbf H^*\big(\Pdot_{|_{X-Y}}\big)$ is locally constant in degree $-d$ and zero in other degrees, then we conclude that $\mathbb H^{-d}(X; \Pdot)$ injects into the global sections of $\mathbf H^{-d}\big(\Pdot_{|_{X-Y}}\big)$.

\vskip 0.2in

Suppose we have 
$$\begin{CD}X @>\pi>> Y @>f>> \mathbb C,\end{CD}$$ 

\medskip

\noindent where $\pi$ is proper and $\hat\pi :
\pi^{-1}f^{-1}(0) \rightarrow f^{-1}(0)$ is the restriction of $\pi$.  Then, for all
$\Adot \in \mathbf D^b_c(X)$, there are natural isomorphisms
$$R\hat\pi_*(\psi_{f\circ\pi}\Adot) \cong \psi_f(R\pi_*\Adot)
\text{ \ and \ } R\hat\pi_*(\phi_{f\circ\pi}\Adot) \cong \phi_f(R\pi_*\mathbf A^
\bullet)  .$$

\vskip .1in

Suppose that $\mathbf Q^\bullet$ is  perverse on $X - f^{-1}(0)$, $i : X - f^{-1}(0) \rightarrow X$ is the inclusion, and that $Ri_*\mathbf Q^\bullet$ is constructible (this is automatic in the algebraic setting). Then it is easy to
see that $Ri_*\mathbf Q^\bullet$ satisfies the cosupport condition; moreover, by
combining the fact that $\psi_f(Ri_*\mathbf Q^\bullet)[-1]$ is  perverse on $f^{-1}(0)$ with
 the Wang sequences on stalk cohomology, one can prove that $Ri_*\mathbf
Q^\bullet$ also satisfies the support condition - hence, $Ri_*\mathbf Q^\bullet$ is perverse.
In an analogous fashion, one obtains that $Ri_!\mathbf
Q^\bullet$ is perverse (if the base ring is a field, this can be obtained by dualizing) provided that $Ri_!\mathbf
Q^\bullet$ is constructible. Note that if $\mathbf Q^\bullet$ is the restriction to $X - f^{-1}(0)$ of a constructible complex on all of $X$, then the constructibility of $Ri_*\mathbf Q^\bullet$ and $Ri_!\mathbf Q^\bullet$ are automatic.

\bigskip

The functors $\psi_f[-1]$ and $\mathcal D$ commute, up to natural isomorphism, as do
$\phi_f[-1]$ and $\mathcal D$.  That is, there are natural isomorphisms $\alpha: \psi_f[-1]\circ\mathcal D\rightarrow \mathcal D\circ\psi_f[-1]$ and $\beta: \phi_f[-1]\circ\mathcal D\rightarrow \mathcal D\circ\phi_f[-1]$.

Furthermore, via these natural isomorphisms, the natural transformations $\operatorname{can}$ and $\operatorname{var}$ between $\psi_f$ and $\phi_f$ are dual; more precisely, for all
$\Adot \in \mathbf D^b_c(X)$, the relationship between ${\operatorname{can}}_{\mathcal D\Adot}:\psi_f[-1]\mathcal D\Adot\rightarrow \phi_f[-1]\mathcal D\Adot$ and ${\operatorname{var}}_{\Adot}:\phi_f[-1]\Adot\rightarrow\psi_f[-1]\Adot$ is given by 
$${\operatorname{can}}_{\mathcal D\Adot}= {\beta_{\Adot}^{-1}}\circ\mathcal D({\operatorname{var}}_{\Adot})\circ\alpha_{\Adot}.$$

\medskip

\noindent{\bf Characterizing zero morphisms between perverse sheaves}
\smallskip

Suppose that $T:\Pdot\rightarrow\Qdot$ is a morphism of perverse sheaves on $X$. Then, $T$ is the zero morphism if and only if, for all $\mathbf p\in X$, for all germs of complex analytic $f: (X, \mathbf p)\rightarrow (\mathbb C, \mathbf 0)$ such that $\dim_{\mathbf p}\operatorname{supp}\phi_f[-1]\Pdot\leq 0$ and $\dim_{\mathbf p}\operatorname{supp}\phi_f[-1]\Qdot\leq 0$,  the induced map $\phi_f[-1]T$ from $H^0(\phi_f[-1]\Pdot)_{\mathbf p}$ to $H^0(\phi_f[-1]\Qdot)_{\mathbf p}$ is zero. One direction is trivial; for lack of a convenient reference, we will sketch a proof of the other direction.

Consider the short exact sequences in ${\operatorname{Perv}}(X)$:
$$
\begin{CD}0\rightarrow  \ \operatorname{ker} T\rightarrow \Pdot @>\alpha>> \ \operatorname{im} T \ \rightarrow 0\hskip .2in \text{and}\hskip .2in 0\rightarrow  \ \operatorname{im} T @>\beta>>  \ \Qdot \rightarrow\operatorname{coker} T\rightarrow 0,\end{CD}
$$
where $\beta\circ\alpha = T$. Whitney stratify $X$ by a stratification with respect to which $\Pdot$, $\Qdot$, $\operatorname{im} T$, $\operatorname{ker} T$, and $\operatorname{coker} T$ are constructible. Assume that $\operatorname{im} T\neq 0$, but that, for all germs of complex analytic $f: (X, \mathbf p)\rightarrow (\mathbb C, \mathbf 0)$ such that $\dim_{\mathbf p}\operatorname{supp}\phi_f[-1]\Pdot\leq 0$ and $\dim_{\mathbf p}\operatorname{supp}\phi_f[-1]\Qdot\leq 0$,  the induced map $\phi_f[-1]T$ from $H^0(\phi_f[-1]\Pdot)_{\mathbf p}$ to $H^0(\phi_f[-1]\Qdot)_{\mathbf p}$ is zero. We wish to derive a contradiction.

Let $S$ be a top-dimensional stratum contained in $\operatorname{supp}(\operatorname{im}T)$, and let $\mathbf p\in S$. Thus, the stalk cohomology of $\operatorname{im}T$ at $\mathbf p$ is a non-zero $R$-module $M$ in degree $-\dim S$. We will show that, in fact, $M$ must be zero.

Let $f$ be the germ of a Morse function for the stratification at $\mathbf p$; this means that $f$ extends to a complex analytic germ $\tilde f$ on an ambient affine space such that ${\tilde f}_{|_S}$ has a Morse singularity, and that $d_{\mathbf p}\tilde f$ is not a degenerate covector, i.e., is not a limiting covector from a bigger stratum. It follows that, at $\mathbf p$, the dimensions of the supports of $\phi_f[-1]\Pdot$, $\phi_f[-1]\Qdot$, $\phi_f[-1]\operatorname{im}T$, $\phi_f[-1]\operatorname{ker} T$, and $\phi_f[-1]\operatorname{coker} T$ are all, at most, $0$, and $H^0(\phi_f[-1]\operatorname{im}T)_{\mathbf p}\cong M$.

Hence, we have the short exact sequences of $R$-modules:

$$
\begin{CD}0\rightarrow  \ H^0(\phi_f[-1]\operatorname{ker} T)_{\mathbf p}\rightarrow H^0(\phi_f[-1]\Pdot)_{\mathbf p} @>\phi_f[-1]\alpha>> \ H^0(\phi_f[-1]\operatorname{im} T)_{\mathbf p} \ \rightarrow 0\end{CD}$$
and
$$
\begin{CD}0\rightarrow   H^0(\phi_f[-1]\operatorname{im} T)_{\mathbf p} @>\phi_f[-1]\beta>>  H^0(\phi_f[-1]\Qdot)_{\mathbf p} \rightarrow H^0(\phi_f[-1]\operatorname{coker} T)_{\mathbf p}\rightarrow 0,\end{CD}
$$

\medskip

\noindent where $\phi_f[-1]\beta\circ\phi_f[-1]\alpha = \phi_f[-1]T$. 

\medskip

Thus, the image of $\phi_f[-1]\beta\circ\phi_f[-1]\alpha = \phi_f[-1]T$ is isomorphic to $H^0(\phi_f[-1]\operatorname{im} T)_{\mathbf p}\cong M$. But, by assumption, this image is zero, which contradicts that $M$ was non-zero.

\vskip .2in

\noindent{\bf The Unipotent and Complementary Nearby and Vanishing Cycles}

\medskip

Suppose that our base ring is a field, that $\Pdot$ is a perverse sheaf on $X$, and that we have an analytic function $f:X\rightarrow\mathbb C$. Then, (globally, in the algebraic setting; locally, in the analytic setting) there exists $n$ such that, in ${\operatorname{Perv}}(f^{-1}(0))$, for all $i\geq n$,  
$$\operatorname{ker}\left(\operatorname{id}-T_f\right)^n \ = \ \operatorname{ker}\left(\operatorname{id}-T_f\right)^i \ \subseteq \ \psi_f[-1]\Pdot$$
and 
$$\operatorname{ker}\left(\operatorname{id}-\widetilde T_f\right)^n=\operatorname{ker}\left(\operatorname{id}-\widetilde T_f\right)^i\ \subseteq \ \phi_f[-1]\Pdot.$$

\medskip

The perverse sheaves $\operatorname{ker}\left(\operatorname{id}-T_f\right)^n$ and $\operatorname{ker}\big(\operatorname{id}-\widetilde T_f\big)^n$ are maximal perverse subsheaves of $\psi_f[-1]\Pdot$ and $\phi_f[-1]\Pdot$, respectively, on which $T_f$ and $\widetilde T_f$, respectively, act unipotently; they are denoted by $\psi^u_f[-1]\Pdot$ and $\phi^u_f[-1]\Pdot$, and are called the {\it unipotent} nearby and vanishing cycles, respectively.

We refer to the images $\psi_f^\perp[-1]\Pdot:=\operatorname{im}\left(\operatorname{id}-T_f\right)^n$ and $\phi_f^\perp[-1]\Pdot:=\operatorname{im}\left(\operatorname{id}-\widetilde T_f\right)^n$ as the {\it complementary} nearby and vanishing cycles, respectively.

\medskip

It is easy to show that the maps $\operatorname{can}:\psi_f[-1]\Pdot\rightarrow \phi_f[-1]\Pdot$ and $\operatorname{var}:\phi_f[-1]\Pdot\rightarrow \psi_f[-1]\Pdot$ restrict to isomorphisms
$$
\begin{CD}\psi_f^\perp[-1]\Pdot@>{\operatorname{can}}^\perp>\cong>\phi_f^\perp[-1]\Pdot\hskip .4in\text{and}\hskip .4in \phi_f^\perp[-1]\Pdot@>{\operatorname{var}}^\perp>\cong>\psi_f^\perp[-1]\Pdot\end{CD}
$$
and, hence, that $\operatorname{id}-T_f$ and $\operatorname{id}-\widetilde T_f$ restrict to isomorphisms on $\psi_f^\perp[-1]\Pdot$ and $\phi_f^\perp[-1]\Pdot$, respectively.

\smallskip

Furthermore, letting $j:f^{-1}(0)\rightarrow X$ denote the inclusion, restriction yields distinguished triangles
$$
\begin{CD}j^*[-1]\Pdot@>\ {\operatorname{comp}}^u \ >> \psi^u_f[-1]\Pdot@>\ {\operatorname{can}}^u\ >>\phi^u_f[-1]\Pdot\longrightarrow j^*\Pdot\end{CD}
$$ 
and
$$\begin{CD}
\phi^u_f[-1]\Pdot@>\ {\operatorname{var}}^u\ >>\psi^u_f[-1]\Pdot@>\ {\operatorname{pmoc}}^u\ >> j^![1]\Pdot\longrightarrow\phi^u_f\Pdot.\end{CD}
$$

\smallskip

There are direct sum splittings $\psi_f[-1]\Pdot\cong \psi^u_f[-1]\Pdot\oplus\psi^\perp_f[-1]\Pdot$, $\phi_f[-1]\Pdot\cong \phi^u_f[-1]\Pdot\oplus\phi^\perp_f[-1]\Pdot$ and, via these isomorphisms, $\operatorname{can}=\big({\operatorname{can}}^u, {\operatorname{can}}^\perp\big)$ and $\operatorname{var}=\big({\operatorname{var}}^u, {\operatorname{var}}^\perp\big)$.

\smallskip

If the base field is algebraically closed, then $\psi^\perp_f[-1]\Pdot$ and $\phi^\perp_f[-1]\Pdot$ can be further decomposed as direct sums of generalized eigenspaces of their respective monodromies and, of course, these direct sum decompositions are compatible with the monodromy maps.

\smallskip

Suppose one is in the special case where $\operatorname{id}-\widetilde T_f$ is an isomorphism. Then $\phi^u_f[-1]\Pdot=0$ and, using the discussion above, we have
$$
\psi_f[-1]\Pdot\cong \psi^u_f[-1]\Pdot\oplus\psi^\perp_f[-1]\Pdot\cong j^*[-1]\Pdot\oplus \phi^\perp_f[-1]\Pdot \cong j^*[-1]\Pdot\oplus \phi_f[-1]\Pdot.
$$

\bigskip

\noindent{\bf Gabber's result}

\smallskip

Suppose that our base ring is a perfect field (e.g., a field of characteristic $0$, a finite field, or an algebraically closed field), that $\Pdot$ is a perverse sheaf on $X$, and that we have an analytic function $f:X\rightarrow\mathbb C$. Then, globally, in the algebraic setting, and locally, in the analytic setting, the derived category/perverse sheaf version of the Jordan-Chevalley Decomposition tells us that the monodromy map $T_f$ can be factored, in a unique way, as $T_f=F(\operatorname{id}+N)$, where $F$ is semisimple, $N$ is nilpotent, and $F$ and $N$ commute.

Given the nilpotent operator $N:\psi_f[-1]\Pdot\rightarrow\psi_f[-1]\Pdot$, it is an exercise in linear algebra to show that there exists a unique (bounded) increasing filtration $W^\bullet$ of $\psi_f[-1]\Pdot$ such that, for all $i$, $N(W^i)\subseteq W^{i-2}$, and such that $N^i$ induces an isomorphism from the $i$-th graded object ${\operatorname{Gr}}^i\psi_f[-1]\Pdot$ to the $(-i)$-th graded object ${\operatorname{Gr}}^{-i}\psi_f[-1]\Pdot$.

Now, suppose that we consider the case where $\Pdot$ is replaced by an intersection cohomology complex $\Idot$, with, possibly, non-constant coefficients. Gabber's theorem deals with the algebraic setting and $\ell$-adic perverse sheaves. Locally, for analytic spaces and maps, and using $\mathbb Q$ for the base ring, the result is due to M. Saito: 

\bigskip

\noindent{\bf Theorem}: The graded pieces ${\operatorname{Gr}}^i\psi_f[-1]\Idot$ are semisimple in ${\operatorname{Perv}}(f^{-1}(0))$, i.e., they are direct sums of intersection cohomology complexes.

\vskip .3in

\noindent {\bf The Sebastiani-Thom Isomorphism}

\vskip .2in

Let $f:X\rightarrow\mathbb C$ and $g:Y\rightarrow\mathbb C$ be complex analytic functions. Let $\pi_1$ and $\pi_2$ denote the
projections of $X\times Y$ onto $X$ and $Y$, respectively.  Let $\Adot$ and $\Bdot$ be bounded, constructible complexes of
sheaves of $R$-modules on
$X$ and $Y$, respectively. In this situation, $\Adot\lboxtimes\Bdot := \piten$. Let us adopt the similar notation
$f\boxplus g:= \map$.

\bigskip

Let  $p_1$ and $p_2$ denote the projections of $V(f)\times V(g)$ onto $V(f)$ and $V(g)$, respectively, and let
$k$ denote the inclusion of $V(f)\times V(g)$ into $V(f\boxplus g)$.

\bigskip

\noindent{\bf Theorem (Sebastiani-Thom Isomorphism)}: {\it There is a natural
isomorphism
$$
k^*\phi_{{}_{f\boxplus g}}[-1]\big(\Adot\lboxtimes\Bdot\big)\ \cong\ \phi_{{}_f}[-1]\Adot\lboxtimes \phi_{{}_g}[-1]\Bdot,
$$
and this isomorphism commutes with the corresponding monodromies. 

\smallskip

Moreover, if we let $\mathbf p:=(\mathbf x, \mathbf y)\in
X\times Y$ be such that
$f(\mathbf x) = 0$ and $g(\mathbf y)=0$, then, in an open neighborhood of $\mathbf p$, the
complex 
$\phi_{{}_{f\boxplus g}}[-1]\big(\Adot\lboxtimes\Bdot\big)$ has support contained in $V(f)\times V(g)$, and, in any open set
in which we have this containment, there are natural isomorphisms
$$
\phi_{{}_{f\boxplus g}}[-1]\big(\Adot\lboxtimes\Bdot\big)\ \cong\ k_!(\phi_{{}_f}[-1]\Adot\lboxtimes \phi_{{}_g}[-1]\Bdot)\
\cong\ k_*(\phi_{{}_f}[-1]\Adot\lboxtimes \phi_{{}_g}[-1]\Bdot).
$$}

\bigskip

\noindent{\bf A'Campo's Theorem}: Suppose that the base ring is a field; if  $m_x$ denotes the maximal ideal of $X$ at $x$ and $f \in m_x^2$, then the
Lefschetz number of the map 
$$\begin{CD}\mathbf H^*(\psi_f[-1]\Adot)_x @>T_f>> \mathbf H^*(\psi_f[-1]\Adot)_x\end{CD}$$
equals $0$, i.e., $$\begin{CD}\sum_i(-1)^i\ \operatorname{Trace}\{\mathbf H^i(\psi_f[-1]\Adot)_x @>T_f>> \mathbf H^i(\psi_f[-1]\Adot)_x\} \ = \
0.\end{CD}$$

\bigskip

\noindent{\bf The Monodromy Theorem}: Let $\Adot\in \mathbf D^b_c(X)$ and $f: X \rightarrow \mathbb C$.  The monodromy
automorphism $$\begin{CD}\psi_f[-1]\Adot @>T_f>> \psi_f[-1]\Adot\end{CD}$$ 

\medskip

\noindent induces a map on
stalk cohomology which is quasi-unipotent, i.e., letting $T_f$ also denote the map on stalk
cohomology, this means that there exist integers $k$ and $j$ such that $(\operatorname{id} - T_f^k)^j =
0$. Over $\mathbb C$, this is equivalent to all of the eigenvalues of $T_f$ being roots of unity.

\bigskip

If $\Pdot$ is a  perverse sheaf, then we may use the
Abelian structure of the category ${\operatorname{Perv}}(X)$ to investigate the map
$$\begin{CD}\psi_f[-1]\Pdot @>T_f>> \psi_f[-1]\Pdot.\end{CD}$$  

\medskip

\noindent The Monodromy Theorem implies that this morphism can be factored (globally, in the algebraic setting; locally, in the analytic setting)
into $T_f = F(\operatorname{id} + N)$, where $F$ has finite order, $N$ is nilpotent, and $F$ and $N$ commute. 

\bigskip

\noindent\S4. {\bf Some Quick Applications}  

\vskip .3in

The applications of perverse sheaves are widespread and are frequently quite deep -
particularly for those applications which rely on the decomposition theorem.  For
beautiful discussions of these applications, we highly recommend \cite{Mac1} and
\cite{Mac2}.  We shall not describe any of these applications here; rather we shall
give some fairly easy results on general Milnor fibers.  These results are ``easy''
now that we have all the machinery of the first three sections at our disposal. 
While the applications below could undoubtedly be proved without the general theory
of perverse sheaves, with this theory in hand, the results and their proofs can be
presented in a unified manner and, what is more, the proofs become mere exercises.

\vskip .3in

Consider the classical case of
the Milnor fiber of a non-zero map $f: (\mathbb C^{n+1} , \mathbf 0) \rightarrow (\mathbb C ,
0)$.  Let $X = \mathbb C^{n+1}$ and let $s = \text{ dim }\Sigma f$.  Then, as $X$ is a
manifold, $\mathbb Z_X^\bullet[n+1]$ is a perverse sheaf and so $\phi_f\mathbb
Z_X^\bullet[n]$ is perverse on $f^{-1}(0)$ with support only on $\Sigma f$.  It
follows that the stalk cohomology of $\phi_f\mathbb Z_X^\bullet[n]$ is non-zero only for
degrees $i$ with $-s \leqslant i \leqslant 0$; that is, we recover the well-known
result that the reduced cohomology of the Milnor fiber can be non-zero only in
degrees $i$ such that $n-s \leqslant i \leqslant n$.

A much more general case is just as easy to derive from the machinery that we have. 
Suppose that $X$ is a purely $(n+1)$-dimensional local complete intersection with arbitrary
singularities.  Let $\mathcal S$ be a Whitney stratification of $X$.  Let  $\mathbf p\in
X$ be such that $\text{dim}_\mathbf pf^{-1}(0) = n$, and let $F_{{}_{f, \mathbf p}}$ denote
the Milnor fiber of $f$ at $\mathbf p$.  Then, as $X$ is a
local complete intersection, $\mathbb Z_X^\bullet[n+1]$ is a perverse sheaf and so
$\phi_f\mathbb Z_X^\bullet[n]$ is perverse on $f^{-1}(0)$ with support only on 
$\Sigma_{{}_{\mathcal S}}f$.  It
follows that the stalk cohomology of $\phi_f\mathbb Z_X^\bullet[n]$ is non-zero only for
dimensions $i$ with $-\text{\ dim}_\mathbf p\Sigma_{{}_{\mathcal S}}f \leqslant i \leqslant
0$.  Hence, the reduced cohomology of $F_{{}_{f, \mathbf p}}$ can be non-zero only in
degrees $i$ such that $(n-\text{\ dim}_\mathbf p\Sigma_{{}_{\mathcal S}}f) \leqslant i \leqslant
n$.

While this general statement could no doubt be proved by induction on hyperplane sections,
the above proof via general techniques avoids the re-working of many technical lemmas on
privileged neighborhoods and generic slices.

\vskip .2in

Another application relates to the homotopy-type of the complex link of a space at a
point; for instance, for an $s$-dimensional local complete intersection, the complex link
has the homotopy-type of a bouquet of spheres of real dimension $s-1$.  In terms of
vanishing cycles and perverse sheaves, we only obtain this result up to cohomology:  let
$(X,x)$ be a germ of an analytic space embedded in some $\mathbb C^n$, and assume $s
:= \text{ dim }X = \text{ dim}_xX$.  Suppose that we have a perverse sheaf,
$\Pdot$, on $X$ (e.g., the shifted constant sheaf, if $X$ is a local complete
intersection).  Let $l$ be a generic linear form, and consider $\phi_{l - l(x)}[-1]\Pdot$; this is a perverse sheaf on an $(s-1)$-dimensional space and, as $l$ is generic,
it is supported at the single point $x$ (because the hyperplane slice $l = l(x)$ can be
chosen to transversely intersect all the strata of any stratification with respect to
which $\Pdot$ is constructible - except, possibly, the point-stratum $x$
itself).  Hence, $H^*(\phi_{l - l(x)}[-1]\Pdot)_x$ is (possibly) non-zero only
in degree $0$.  In the case of the shifted constant sheaf $\mathbb Z^\bullet_X[s]$ on a local complete intersection,
this gives the desired result.

\vskip .2in

For our final application, we wish to investigate functions with one-dimensional
critical loci; we must first set up some notation.  

Let $\mathcal U$ be an open neighborhood
of the origin in $\mathbb C^{{}^{n+1}}$ and suppose that $f: (\mathcal U, \mathbf 0) \rightarrow
(\mathbb C, 0)$ has a one-dimensional critical locus at the origin, i.e., $\text{dim}_\mathbf
0\Sigma f = 1$.  Suppose also that $\mathcal U$ is small enough so that, inside $\mathcal U$, $\Sigma f\subseteq V(f)$. The reduced cohomology of the Milnor fiber, $F_{{}_{f, \mathbf 0}}$, of
$f$ at the origin is possibly non-zero only in dimensions $n-1$ and $n$.  We wish to show
that the $(n-1)$-st cohomology group embeds inside another group which is fairly easy to
describe; thus, we obtain a bound on the $(n-1)$-st Betti number of the Milnor fiber of $f$.

For each component $\nu$ of $\Sigma f$, one may consider a generic hyperplane slice, $H$, at points
$\mathbf p \in \nu-\{\mathbf 0\}$ close to the origin; then, the restricted function, $f_{|_H}$, will have an
isolated critical point at $\mathbf p$.  By shrinking the neighborhood $\mathcal U$ if
necessary, we may assume that each $\nu-\{\mathbf 0\}$ is homotopy-equivalent to a circle and that the Milnor number of the isolated singularity of
$f_{|_H}$ at $\mathbf p$ is independent of the point $\mathbf p \in
\nu-\{\mathbf 0\}$; denote this value by $\mu_\nu^\circ$. As $\nu-\{\mathbf 0\}$ is
homotopy-equivalent to a circle, there is a single monodromy map from the Milnor fiber of 
$f_{|_H}$ at $\mathbf p\in\nu-\{\mathbf 0\}$ to itself, which induces a map on the middle-dimensional cohomology, i.e., a map $h_\nu: \mathbb Z^{{}^{\mu^\circ_\nu}} 
\rightarrow 
\mathbb Z^{{}^{\mu^\circ_\nu}}$.  We wish to show that  $\widetilde H^{n-1}(F_{{}_{f, \mathbf
0}})$ (with integer coefficients) injects into $\oplus_\nu \text{ker}(\operatorname{id}-h_\nu)$.

 As $\mathbb Z_{{}_\mathcal U}^\bullet[n+1]$ is perverse,
$\phi_f[-1]\mathbb Z_{{}_\mathcal
U}^\bullet[n+1]$ is perverse with support equal to $\Sigma f$. Therefore the restriction
$$\mathbf K^\bullet:=\big(\phi_f[-1]\mathbb Z_{{}_\mathcal
U}^\bullet[n+1]\big)_{{\big |}_{\Sigma f}}$$
is perverse.
Let $j:\{\mathbf 0\}\hookrightarrow\Sigma f$ denote the closed inclusion, and let $i:\Sigma f -\{\mathbf 0\}\hookrightarrow \Sigma f$ denote the    open inclusion. 

We always have the distinguished triangle

\vskip .2in

\vbox{
$\hskip 2.35in j_!j^!\mathbf K^\bullet \longrightarrow \ \mathbf K^\bullet$ \vskip .03in
$\hskip 2.35in \ _{[1]} \ \nwarrow\hskip .25in \swarrow$ \vskip .03in $\hskip
2.6in Ri_*i^*\mathbf K^\bullet$
}

\vskip .1in

We wish to examine the associated stalk cohomology long exact sequence at the origin.

\vskip .1in

First, we have that $H^{-1}((j_!j^!\mathbf K^\bullet)_\mathbf 0) = H^{-1}(j^!\mathbf
K^\bullet)$ and so, by the cosupport condition for perverse sheaves,
$H^{-1}((j_!j^!\mathbf K^\bullet)_\mathbf 0) = 0$.

\vskip .1in

Now, we need to look more closely at the sheaf $Ri_*i^*\mathbf K^\bullet$.  $i^*\mathbf
K^\bullet$ is the restriction of $\mathbf K^\bullet$ to $\Sigma f-\{\mathbf 0\}$; at points near the origin, this
sheaf has stalk cohomology only in degree $-1$. 
Moreover, the cohomology sheaf $\mathbf H^{-1}(i^*\mathbf K^\bullet)$ is locally constant, and so consists of local systems on each $\nu-\{\mathbf 0\}$.  It follows that 
$$
H^{-1}\big(Ri_*i^*\mathbf K^\bullet\big)_{\mathbf 0}\cong \hyp^{-1}(\Sigma f-\{\mathbf 0\}; i^*\mathbf K^\bullet)\cong \oplus_\nu \text{ker}(\text{id} - h_\nu).
$$

Thus, when we consider the long exact sequence on stalk cohomology associated to our
distinguished triangle, we find -- starting in dimension $n-1$ -- that it begins
$$
0 \rightarrow \widetilde H^{n-1}(F_{{}_{f, \mathbf 0}}) \rightarrow \oplus_\nu \text{ker}(\operatorname{id} - h_\nu)
\rightarrow \dots .
$$
The desired conclusion follows.

\bigskip

\noindent\S5. {\bf Truncation and Perverse Cohomology}  

\medskip

This section is taken entirely from \cite{BBD}, \cite{G-M3}, and \cite{K-S}.

\medskip

There are (at least) two forms of truncation associated to an object $\Fdot\in \mathbf
D^b_c(X)$ -- one form of truncation is related to the ordinary cohomology of the complex, while the
other form leads to something called the {\it perverse cohomology} or {\it perverse projection}.  These two types
of truncation bear little resemblance to each other, except in the general framework of a
$t$-structure on $\mathbf D^b_c(X)$.

\vskip .1in

Loosely speaking, a $t$-structure on $\mathbf D^b_c(X)$ consists of two full subcategories, denoted 
$\mathbf D^{{}^{\leqslant 0}}(X)$ and $\mathbf D^{{}^{\geqslant 0}}(X)$, such that for any $\Fdot\in \mathbf D^b_c(X)$, there exist $\mathbf E^\bullet\in \mathbf D^{{}^{\leqslant 0}}(X)$, 
$\mathbf G^\bullet\in \mathbf D^{{}^{\geqslant 0}}(X)$, and a distinguished triangle
\vskip .2in

\vbox{
$\hskip 2.5in \mathbf E^\bullet \longrightarrow \ \Fdot$ \vskip .03in
$\hskip 2.35in \ _{[1]} \ \nwarrow\hskip .15in \swarrow$ \hskip .5in ; \vskip .03in $\hskip
2.65in \mathbf G^\bullet[-1]$
}

\vskip .1in

\noindent moreover, such $\mathbf E^\bullet$ and $\mathbf G^\bullet$ are required to be unique up to
isomorphism in $\mathbf D^b_c(X)$.

\medskip

Given a $t$-structure as above, and using the same notation, we write $\mathbf E^\bullet =
\tau_{\leqslant 0}\Fdot$ (the {\it truncation of $\Fdot$ below $0$}) and  $\mathbf
G^\bullet =
\tau^{\geqslant 0}\left(\Fdot[1]\right)$ (the {\it truncation of $\Fdot[1]$ above
$0$}); these are the basic truncation functors associated to the
$t$-structure.

In addition, we write $\mathbf D^{{}^{\leqslant n}}(X)$ for 
$$\mathbf D^{{}^{\leqslant 0}}(X)[-n] := \left\{\Fdot[-n]\ |\ \Fdot\in 
\mathbf D^{{}^{\leqslant 0}}(X) \right\} ,
$$
and we analogously write  $\mathbf D^{{}^{\geqslant n}}(X)$ for 
$\mathbf D^{{}^{\geqslant 0}}(X)[-n]$.  Thus $\Adot\in \mathbf D^{{}^{\leqslant n}}(X)$ (resp., $ \mathbf D^{{}^{\geqslant n}}(X)$)  if and only if $\Adot[n]\in \mathbf D^{{}^{\leqslant 0}}(X)$ (resp., $ \mathbf D^{{}^{\geqslant 0}}(X)$).

Also, we define 
$\tau_{\leqslant n}\Fdot$ by 
$$\tau_{\leqslant n}\Fdot \ = \ \left(\tau_{\leqslant 0}(\Fdot[n])\right)[-n] \ = \ ([-n]\circ \tau_{\leqslant 0}\circ [n])\Fdot ,$$ 
and
we analogously define 
$\tau^{\geqslant n}\Fdot$ as  $([-n]\circ \tau^{\geqslant 0}\circ [n])\Fdot$.

\vskip 0.1in

It follows that
$$
\tau_{\leqslant n}(\Fdot[k])= (\tau_{\leqslant n+k}\Fdot)[k]\hskip 0.2in\text{ and } \hskip 0.2in \tau^{\geqslant n}(\Fdot[k])= (\tau^{\geqslant n+k}\Fdot)[k].
$$

\vskip .1in

Note that $\tau_{\leqslant n}\Fdot\in \mathbf D^{{}^{\leqslant n}}(X)$,  
$\tau^{\geqslant n}\Fdot\in \mathbf D^{{}^{\geqslant n}}(X)$ and, for all $n$, we have a
distinguished triangle
\vskip .2in

\vbox{
$\hskip 2.5in \tau_{\leqslant n}\Fdot \longrightarrow \ \Fdot$ \vskip .03in
$\hskip 2.5in \ _{[1]} \ \nwarrow\hskip .20in \swarrow$ \hskip .5in . \vskip .03in $\hskip
2.7in \tau^{\geqslant n+1}\Fdot$
}

\bigskip

Writing $\simeq$ to denote natural isomorphisms between functors: for all $a$ and $b$, 
$$\tau_{\leqslant b}\circ\tau^{\geqslant a}\simeq\tau^{\geqslant a}\circ\tau_{\leqslant b} ,$$
$$\tau_{\leqslant b}\circ\tau_{\leqslant a}\simeq\tau_{\leqslant a}\circ\tau_{\leqslant b} ,$$ 
and
$$\tau^{\geqslant b}\circ\tau^{\geqslant a}\simeq\tau^{\geqslant a}\circ\tau^{\geqslant b} .$$
 
If $a\geqslant b$, then 
$$\tau_{\leqslant b}\circ\tau_{\leqslant a}\simeq\tau_{\leqslant b} ,$$
and  
$$\tau^{\geqslant a}\circ\tau^{\geqslant b}\simeq\tau^{\geqslant a} .$$  

Also, if $a>b$, then 
$$\tau_{\leqslant b}\circ\tau^{\geqslant a} = \tau^{\geqslant a}\circ\tau_{\leqslant b} = 0 .$$

\vskip .2in

The {\it heart} of the $t$-structure is defined to be the full subcategory $\mathcal C := \mathbf
D^{{}^{\leqslant 0}}(X) \cap  \mathbf D^{{}^{\geqslant 0}}(X)$; this is always an Abelian category. 
We wish to describe the kernels and cokernels in this category.

\vskip .1in

Let $\mathbf E^\bullet, \Fdot \in \mathcal C$ and let $f$ be a morphism from $\mathbf E^\bullet$
to $\Fdot$.  We can form a distinguished triangle in $\mathbf D^b_c(X)$

\vskip .1in

\vbox{
$\begin{CD}\hskip 2.5in \mathbf E^\bullet @>\ \ f\ \ >> \ \Fdot\end{CD}$ \vskip .03in
$\hskip 2.35in \ _{[1]} \ \nwarrow\hskip .35in \swarrow$ \hskip .5in , \vskip .03in $\hskip
2.80in \mathbf G^\bullet$
}
\noindent where $\mathbf G^\bullet$ need not be in $\mathcal C$.  Then, up to natural isomorphism,
$$\operatorname{coker} f = \tau^{\geqslant 0}\mathbf G^\bullet \text{  and  }
\operatorname{ker} f = \tau_{\leqslant 0}(\mathbf G^\bullet[-1]) .$$

\vskip .2in

We define cohomology associated to a $t$-structure as follows.  Define ${}^tH^0(\Fdot)$ to
be $\tau^{\geqslant 0}\tau_{\leqslant 0}\Fdot$; this is naturally isomorphic to 
$\tau_{\leqslant 0}\tau^{\geqslant 0}\Fdot$.  Now, define 
${}^tH^n(\Fdot)$ to
be $${}^tH^0(\Fdot[n]) = \left(\tau^{\geqslant n}\tau_{\leqslant n}\Fdot\right)[n].$$  Note that this cohomology does not
give back modules or even sheaves of modules, but rather gives back complexes which are objects in the heart of the
$t$-structure.

\bigskip

If $\Fdot\in \mathbf D^b_c(X)$, then the following are equivalent:

\vskip .1in

\noindent 1)\hskip .2in $\Fdot\in \mathbf D^{{}^{\leqslant 0}}(X)$ (resp. $\mathbf D^{{}^{\geqslant 0}}(X)$);

\vskip .1in

\noindent 2)\hskip .2in the morphism $\tau_{\leqslant 0}\Fdot\rightarrow\Fdot$ is an isomorphism (resp. the morphism
$\Fdot\rightarrow\tau^{\geqslant 0}\Fdot$ is an 

\noindent\hbox{}\hskip .35in isomorphism);

\vskip .1in

\noindent 3)\hskip .2in $\tau^{\geqslant 1}\Fdot = 0$ (resp. $\tau_{\leqslant -1}\Fdot = 0)$);

\vskip .1in

\noindent 4)\hskip .2in $\tau^{\geqslant i}\Fdot = 0$ for all $i\geqslant 1$ (resp. $\tau_{\leqslant i}\Fdot = 0$ for all
$i\leqslant -1$);

\vskip .1in

\noindent 5)\hskip .2in there exists $a$ such that $\Fdot\in \mathbf D^{{}^{\leqslant a}}(X)$ and ${}^tH^i(\Fdot)=0$ for all
$i\geqslant 1$ (resp. there exists $a$ such that $\Fdot\in \mathbf D^{{}^{\geqslant a}}(X)$ and ${}^tH^i(\Fdot)=0$ for all
$i\leqslant -1$).

\vskip .3in

It follows that, if $\Fdot\in \mathbf D^b_c(X)$, then the following are equivalent:

\vskip .1in

\noindent 1)\hskip .2in $\Fdot\in \mathcal C$;

\vskip .1in

\noindent 2)\hskip .2in  ${}^tH^0(\Fdot)$ is isomorphic to $\Fdot$;

\vskip .1in

\noindent 3)\hskip .2in  there exist $a$ and $b$ such that $\Fdot\in \mathbf
D^{{}^{\geqslant a}}(X)$, $\Fdot\in \mathbf
D^{{}^{\leqslant b}}(X)$, and ${}^tH^n(\Fdot) = 0$ 

\noindent{}\hskip .34in for all $n \neq 0$.

\vskip .3in

As the heart is an Abelian category, we may talk about exact sequences in $\mathcal C$.  Any
distinguished triangle in $\mathbf D^b_c(X)$ determines a long exact sequence of objects in the heart
of the $t$-structure;  if  

\vskip .2in

\vbox{
$\hskip 2.5in \mathbf E^\bullet \longrightarrow \ \Fdot$ \vskip .03in
$\hskip 2.35in \ _{[1]} \ \nwarrow\hskip .15in \swarrow$  
\vskip .03in $\hskip
2.75in \mathbf G^\bullet$
}

\vskip .1in

\noindent is a distinguished triangle in $\mathbf D^b_c(X)$, then the associated long exact sequence in
$\mathcal C$ is
$$\dots\rightarrow {}^tH^{-1}(\mathbf G^\bullet)\rightarrow {}^tH^0(\mathbf E^\bullet)
\rightarrow {}^tH^0(\Fdot)\rightarrow {}^tH^0(\mathbf G^\bullet)
\rightarrow {}^tH^1(\mathbf E^\bullet)\rightarrow\dots .$$

\vskip .2in

We are finished now with our generalities on $t$-structures and wish to, at last, give our two primary
examples.

\vskip .3in

\noindent{\bf The ``ordinary'' $t$-structure}

\vskip .2in

The ``ordinary'' $t$-structure on $\mathbf D^b_c(X)$ is given by
$$\mathbf D^{{}^{\leqslant 0}}(X) = \{ \Fdot\in\mathbf D^b_c(X)\ |\ \mathbf H^i(\Fdot) =
0 \text{ for all } i>0 \}$$
and
$$\mathbf D^{{}^{\geqslant 0}}(X) = \{ \Fdot\in\mathbf D^b_c(X)\ |\ \mathbf H^i(\Fdot) =
0 \text{ for all } i<0 \} .$$

\vskip .1in

The associated truncation functors are the ordinary ones described in \cite{G-M3}.  If $\Fdot\in\mathbf
D^b_c(X)$, then

\vskip .1in

$$\left(\tau_{\leqslant p}\Fdot\right)^n = \begin{cases} 
\mathbf F^n &\text{ if } n<p\\
\text{\bf ker }d^p &\text{ if }  n=p\\0 &  \text{ if }  n>p\end{cases}$$

and

$$\left(\tau^{\geqslant p}\Fdot\right)^n = \begin{cases}
0 &\text{ if } n<p\\
\text{\bf coker }d^{p-1} &\text{ if }  n=p\\\mathbf F^n &  \text{ if }  n>p .\end{cases}$$

\vskip .1in

These truncated complexes are naturally quasi-isomorphic to the complexes 

\vskip .1in

$$\left(\tilde\tau_{\leqslant p}\Fdot\right)^n =\begin{cases}
\mathbf F^n &\text{ if } n\leqslant p\\
\text{\bf Im }d^p &\text{ if }  n=p+1\\0 &  \text{ if }  n>p+1\end{cases}$$

and

$$\left(\tilde\tau^{\geqslant p}\Fdot\right)^n = \begin{cases}
0 &\text{ if } n<p-1\\
\text{\bf Im }d^{p-1} &\text{ if }  n=p-1\\\mathbf F^n &  \text{ if }  n\geqslant p .\end{cases}$$

\vskip .3in

If $\Adot, \Bdot\in\mathbf D^b_c(X)$, then

\vskip .1in

\noindent 1. $\left(\tau_{\leqslant p}\Adot\right)_x = \tau_{\leqslant p}\left(\Adot_x\right)$;

\vskip .2in

\noindent 2. $\mathbf H^k\left(\tau_{\leqslant p}\Adot\right)_x = \begin{cases}\mathbf H^k\left(\Adot\right)_x &\text{ if } k\leqslant p\\0 &\text{ for }k>p .\end{cases}$

\vskip .2in

\noindent 3. If $\phi:\Adot\rightarrow \Bdot$ is a morphism of complexes of sheaves which induces
isomorphisms on the associated cohomology sheaves 
$$\phi^*: \mathbf H^n(\Adot)\cong \mathbf H^n(\Bdot) \text{ for all } n\leqslant p ,$$
then $\tau_{\leqslant p}\phi: \tau_{\leqslant p}\Adot\rightarrow \tau_{\leqslant p}\Bdot$ is a
quasi-isomorphism.

\vskip .2in

\noindent 4.  If $f: X\rightarrow Y$ is a continuous map and $\mathbf C^\bullet$ is a complex of sheaves on $Y$, then 
$$
\tau_{\leqslant p}f^*(\mathbf C^\bullet) \cong f^*\tau_{\leqslant p}(\mathbf C^\bullet) .$$

\vskip .2in

\noindent 5.  If $R$ is a field and $\Adot$ is a complex of sheaves of $R$-modules on $X$ with locally constant
cohomology sheaves, then there are natural quasi-isomorphisms 
$$
\tau^{\geqslant -p}R\mathbf{Hom}^\bullet(\Adot, \mathbf R^\bullet_X)\rightarrow \tau^{\geqslant
-p}R\mathbf{Hom}^\bullet(\tau_{\leqslant p}\Adot, \mathbf R^\bullet_X)\leftarrow R\mathbf{Hom}^\bullet(\tau_{\leqslant p}\Adot, \mathbf R^\bullet_X) .
$$

\bigskip

The heart of
this $t$-structure consists of those complexes which have non-zero cohomology sheaves only in degree
$0$; such complexes are quasi-isomorphic to complexes which are non-zero only in degree $0$.  

The $t$-structure cohomology of a complex $\Fdot$ is essentially the sheaf cohomology of
$\Fdot$; ${}^tH^n(\Fdot)$ is quasi-isomorphic to a complex which has $\mathbf
H^n(\Fdot)$ in degree $0$ and is zero in all other degrees.  With this identification, the
$t$-structure long exact sequence associated to a distinguished triangle is merely the usual long
exact sequence on sheaf cohomology.

\bigskip

We are now going to give the construction of the intersection cohomology complexes as it is presented in \cite{G-M3}.  Our
indexing will look different from that of \cite{G-M3} for several reasons. 

First, we are dealing only with complex analytic
spaces, $X$, and we are using only middle perversity; this accounts for some of the indexing differences. In
addition, in this setting, the intersection cohomology complex defined in \cite{G-M3} would have possibly non-zero cohomology
only in degrees between
$-2\operatorname{dim}_{\mathbb C}X$ and
$-(\operatorname{dim}_{\mathbb C}X)-1$, inclusive. The definition below is shifted by $-\operatorname{dim}_{\mathbb C}X$ from the \cite{G-M3}
definition, and yields a perverse sheaf which has possibly non-zero stalk cohomology only in degrees between
$-\operatorname{dim}_{\mathbb C}X$ and $-1$, inclusive (assuming the space has no isolated points).

\medskip

Let $X$ be a purely $n$-dimensional complex analytic space with a complex analytic Whitney stratification $\mathcal S=
\{S_\alpha\}$. 

  For all
$k$, let
$X^k$ denote the union of the strata of dimension less than or equal to
$k$.  By convention, we set $X^{-1} =
\emptyset$.  Hence, we have a filtration
$$
\emptyset = X^{-1}\subseteq X^0\subseteq X^1\subseteq\dots \subseteq X^{n-1}\subseteq X^n = X .
$$
For all $k$, let $\mathcal U_k := X - X^{n-k}$, and let $i_k$ denote the inclusion $\mathcal U_k\hookrightarrow \mathcal U_{k+1}$.  Let $\mathcal
L^\bullet_{{}_{\mathcal U_1}}$ be a local system on the top-dimensional strata.

Then, the intersection cohomology complex on $X$ with coefficients in $\mathcal L^\bullet_{{}_{\mathcal U_1}}$, as described in section
2, is given by
$$
\mathbf{IC}^\bullet_X(\mathcal L^\bullet_{{}_{\mathcal U_1}}) := \tau_{{}_{\leqslant -1}}R{i_n}_*\dots \tau_{{}_{\leqslant 1-n}}R{i_2}_*
 \tau_{{}_{\leqslant -n}}R{i_1}_*\big(\mathcal L^\bullet_{{}_{\mathcal U_1}}[n]\big) .
$$

\vskip .1in

Up to quasi-isomorphism, this complex is independent of the stratification. Note that the cohomology sheaves of 
$\mathbf{IC}^\bullet_X(\mathcal L^\bullet_{{}_{\mathcal U_1}})$ are supported only in degrees $k$ for which $-n\leqslant k\leqslant -1$ 
(unless $X$ is $0$-dimensional, and then $\mathbf{IC}^\bullet_X(\mathcal L^\bullet_{{}_{\mathcal U_1}})\cong \mathcal L^\bullet_{{}_{\mathcal
U_1}}$).  

\vskip .1in

Also note
that it follows from the construction that there is always a canonical map from the shifted constant sheaf $\mathbf R^\bullet_X[n]$
to 
$\mathbf{IC}^\bullet_X(\mathbf R^\bullet_{{}_{\mathcal U_1}})$ which induces an isomorphism when restricted to $\mathcal U_1$.  

To see
this, consider the canonical morphism $\mathbf R^\bullet_{{}_{\mathcal
U_{k+1}}}[n]\rightarrow  R{i_k}_*i_k^*\mathbf R^\bullet_{{}_{\mathcal U_{k+1}}}[n]$ for each 
$k\geqslant 1$.  As $i_k^*\mathbf R^\bullet_{{}_{\mathcal U_{k+1}}}[n]
\cong
\mathbf R^\bullet_{{}_{\mathcal U_k}}[n]$, we have a canonical map $\mathbf R^\bullet_{{}_{\mathcal U_{k+1}}}[n]\rightarrow 
R{i_k}_*\mathbf R^\bullet_{{}_{\mathcal U_k}}[n]$ and, hence, a canonical map between the truncations 
$\tau_{{}_{\leqslant k-n-1}}\big(\mathbf R^\bullet_{{}_{\mathcal U_{k+1}}}[n]\big)\rightarrow 
\tau_{{}_{\leqslant k-n-1}}R{i_k}_*\big(\mathbf R^\bullet_{{}_{\mathcal U_k}}[n]\big)$.  But, $$\tau_{{}_{\leqslant
k-n-1}}\big(\mathbf R^\bullet_{{}_{\mathcal U_{k+1}}}[n]\big) \cong \mathbf R^\bullet_{{}_{\mathcal U_{k+1}}}[n]$$ and so we have a
canonical map 
$ \mathbf R^\bullet_{{}_{\mathcal U_{k+1}}}[n]\rightarrow \tau_{{}_{\leqslant k-n-1}}R{i_k}_*\big(\mathbf R^\bullet_{{}_{\mathcal
U_k}}[n]\big)$.  By piecing all of these maps together, one obtains the desired morphism.

\vskip .3in

\noindent{\bf The perverse $t$-structure}

\vskip .2in

The {\it perverse} $t$-structure (with middle perversity $\mu$) on $\mathbf D^b_c(X)$ is given by
$${}^{{}^\mu}\mathbf D^{{}^{\leqslant 0}}(X) = \{ \Fdot\in\mathbf D^b_c(X)\ |\
\operatorname{dim}\operatorname{supp^{-j}}\Fdot\leqslant j \text{ for all } j \}$$
and
$${}^{{}^\mu}\mathbf D^{{}^{\geqslant 0}}(X) = \{ \Fdot\in\mathbf D^b_c(X)\ |\
\operatorname{dim}\operatorname{cosupp^j}\Fdot\leqslant j \text{ for all } j \}  .$$
Note that the heart of this $t$-structure is precisely ${\operatorname{Perv}}(X)$, the category of perverse sheaves.   Thus, every distinguished 
triangle in $\mathbf D^b_c(X)$ determines a long exact sequence in the Abelian category ${\operatorname{Perv}}(X)$.

\vskip .1in

We naturally call the $t$-structure cohomology associated to the perverse $t$-structure
 the {\it perverse cohomology} or {\it perverse projection} and denote it in degree $n$ by ${}^{\mu}\hskip -0.02in H^n(\Fdot)$.

\vskip .2in

Let $d$ be an integer, and let $f:Y\rightarrow X$ be a morphism of complex spaces such that, for all $\mathbf x\in X$, 
$\operatorname{dim}f^{-1}(\mathbf x)\leqslant d$. Then,

\vskip .1in

\noindent 1) $f^*$ sends ${}^{{}^\mu}\mathbf D^{{}^{\leqslant 0}}(X)$ to ${}^{{}^\mu}\mathbf D^{{}^{\leqslant d}}(Y)$;

\vskip .1in

\noindent 2) $f^!$ sends ${}^{{}^\mu}\mathbf D^{{}^{\geqslant 0}}(X)$ to ${}^{{}^\mu}\mathbf D^{{}^{\geqslant -d}}(Y)$;

\vskip .1in

\noindent 3) if $\Fdot\in {}^{{}^\mu}\mathbf D^{{}^{\leqslant 0}}(Y)$ and  $Rf_!\Fdot\in \mathbf D^b_c(X)$, then $Rf_!\Fdot\in
{}^{{}^\mu}\mathbf D^{{}^{\leqslant d}}(X)$;

\vskip .1in

\noindent 4) if $\Fdot\in {}^{{}^\mu}\mathbf D^{{}^{\geqslant 0}}(Y)$ and  $Rf_*\Fdot\in \mathbf D^b_c(X)$, then $Rf_*\Fdot\in
{}^{{}^\mu}\mathbf D^{{}^{\geqslant -d}}(X)$.

\vskip 0.1in

\noindent Furthermore, if $f$ is a smooth map of constant relative dimension $d$, then,

\vskip 0.1in

\noindent 5) $f^*$ sends ${}^{{}^\mu}\mathbf D^{{}^{\geqslant 0}}(X)$ to ${}^{{}^\mu}\mathbf D^{{}^{\geqslant d}}(Y)$;

\vskip 0.1in

\noindent 6) $f^!$ sends ${}^{{}^\mu}\mathbf D^{{}^{\leqslant 0}}(X)$ to ${}^{{}^\mu}\mathbf D^{{}^{\leqslant -d}}(Y)$.

\vskip 0.3in

For closed embeddings, we also have the following:

\vskip 0.1in

\noindent Let $g_1, \dots, g_e$ be complex analytic functions on $X$. Let $m$ denote the inclusion of $Y:=V(g_1, \dots, g_e)$ into $X$. Then,

\vskip 0.1in

\noindent i) $m^*$ sends ${}^{{}^\mu}\mathbf D^{{}^{\geqslant 0}}(X)$ to ${}^{{}^\mu}\mathbf D^{{}^{\geqslant -e}}(Y)$, i.e., $m^*[-e]$ sends ${}^{{}^\mu}\mathbf D^{{}^{\geqslant 0}}(X)$ to ${}^{{}^\mu}\mathbf D^{{}^{\geqslant 0}}(Y)$;

\vskip 0.1in

\noindent ii) $m^!$ sends ${}^{{}^\mu}\mathbf D^{{}^{\leqslant 0}}(X)$ to ${}^{{}^\mu}\mathbf D^{{}^{\leqslant e}}(Y)$, i.e., $m^![e]$ sends ${}^{{}^\mu}\mathbf D^{{}^{\leqslant 0}}(X)$ to ${}^{{}^\mu}\mathbf D^{{}^{\leqslant 0}}(Y)$.

\vskip .3in

Consider, in particular, the case in which $m$ is the inclusion of $Y=V(f)$ into $X$, where $V(f)$ contains no irreducible component of $X$, and suppose that the shifted constant sheaf $\Pdot:=\mathbf R^\bullet_X[\operatorname{dim} X]$ is perverse (so $X$ is necessarily pure-dimensional). Then, it is trivial that $m^*[-1]\Pdot\cong \mathbf R^\bullet_Y[\operatorname{dim} Y]$ is in ${}^{{}^\mu}\mathbf D^{{}^{\leqslant 0}}(Y)$. Moreover, i) above implies that $m^*[-1]\Pdot$ is in ${}^{{}^\mu}\mathbf D^{{}^{\geqslant 0}}(Y)$. Thus, $\mathbf R^\bullet_Y[\operatorname{dim} Y]$ is perverse. Hence, by inducting, we recover L\^e's result that the shifted constant sheaf on a local complete intersection is perverse.

\vskip .2in

Suppose that $j:Y\hookrightarrow X$ is inclusion, and $\Adot\in  \mathbf D^b_c(X)$. Later, when we discuss the {\it intermediate extension}, it will be important that, if  $\Adot\in{}^{{}^\mu}\mathbf D^{{}^{\leqslant 0}}$ and $j^*\Adot=0$, then $j^*\pcoh(\Adot)\in {}^{{}^\mu}\mathbf D^{{}^{\leqslant -2}}$. To see this, one uses the perverse truncation triangle to conclude that, since $j^*\Adot=0$, we have an isomorphism $j^*\tau^{\geq 0}\Adot\cong (j^*\tau_{\leq -1}\Adot)[1]$. As $\Adot\in{}^{{}^\mu}\mathbf D^{{}^{\leqslant 0}}$, we have an isomorphism $\tau^{\geq 0}\Adot\cong \pcoh(\Adot)$. Hence, $j^*\pcoh(\Adot)\cong j^*[1]\tau_{\leq -1}\Adot$ and, as $\tau_{\leq -1}\Adot\in {}^{{}^\mu}\mathbf D^{{}^{\leqslant -1}}$ and by our previous property that $j^*$ for an inclusion takes $\mathbf D^{{}^{\leqslant 0}}$ to $\mathbf D^{{}^{\leqslant 0}}$, it follows that 
$$
j^*\pcoh(\Adot)\cong j^*[1]\tau_{\leq -1}\Adot\in  {}^{{}^\mu}\mathbf D^{{}^{\leqslant -2}}.
$$

Dually, if  $\Adot\in{}^{{}^\mu}\mathbf D^{{}^{\geqslant 0}}$ and $j^!\Adot=0$, then $j^!\pcoh(\Adot)\in {}^{{}^\mu}\mathbf D^{{}^{\geqslant 2}}$.

 In particular, if $j$ is the inclusion of a closed analytic subset, $i:X-Y\hookrightarrow X$ is the inclusion of the open complement, and $\Pdot$ is a perverse sheaf on $X-Y$, then
$$
j^*\pcoh(i_!\Pdot)\in {}^{{}^\mu}\mathbf D^{{}^{\leqslant -2}}\hskip .2in\text{and}\hskip .2in j^!\pcoh(Ri_*\Pdot)\in {}^{{}^\mu}\mathbf D^{{}^{\geqslant 2}},
$$
provided that $i_!\Pdot$ and $Ri_*\Pdot$ are constructible (which would be the case, for instance, if $\Pdot$ is isomorphic to the restriction of a constructible complex on all of $X$).

\vskip .2in

Let $f:Y\rightarrow X$ be a morphism of complex spaces such that each point in $X$ has an open neighborhood $\mathcal U$ such that
$f^{-1}(\mathcal U)$ is a Stein space (e.g., an affine map between algebraic varieties). Then,

\vskip .1in

\noindent a) if $\Fdot\in {}^{{}^\mu}\mathbf D^{{}^{\leqslant 0}}(Y)$ and  $Rf_*\Fdot\in \mathbf D^b_c(X)$, then $Rf_*\Fdot\in
{}^{{}^\mu}\mathbf D^{{}^{\leqslant 0}}(X)$;

\vskip .1in

\noindent b) if $\Fdot\in {}^{{}^\mu}\mathbf D^{{}^{\geqslant 0}}(Y)$ and  $Rf_!\Fdot\in \mathbf D^b_c(X)$, then $Rf_!\Fdot\in
{}^{{}^\mu}\mathbf D^{{}^{\geqslant 0}}(X)$.

\vskip .3in

If $f: X\rightarrow\mathbb C$ is an analytic map, then the functors $\psi_f[-1]$ and $\phi_f[-1]$ are
{\it $t$-exact} with respect to the perverse $t$-structures; this means that if $\mathbf E^\bullet\in
{}^{{}^\mu}\mathbf D^{{}^{\leqslant 0}}(X)$ and $\Fdot\in
{}^{{}^\mu}\mathbf D^{{}^{\geqslant 0}}(X)$, then $\psi_f\mathbf E^\bullet[-1]$ and 
$\phi_f\mathbf E^\bullet[-1]$ are in ${}^{{}^\mu}\mathbf D^{{}^{\leqslant 0}}(f^{-1}(0))$, and 
 $\psi_f\Fdot[-1]$ and 
$\phi_f\Fdot[-1]$ are in ${}^{{}^\mu}\mathbf D^{{}^{\geqslant 0}}(f^{-1}(0))$.

In particular, $\psi_f[-1]$ and $\phi_f[-1]$ take perverse sheaves to perverse sheaves and, for any
$\Fdot\in \mathbf D^b_c(X)$,
 $${}^{\mu}\hskip -0.02in H^n(\psi_f\Fdot[-1]) \cong \psi_f{}^{\mu}\hskip -0.02in H^n(\Fdot)[-1] \text{ \ and \ }
{}^{\mu}\hskip -0.02in H^n(\phi_f\Fdot[-1]) \cong
\phi_f{}^{\mu}\hskip -0.02in H^n(\Fdot)[-1] .$$

While there are several proofs in the literature of the fact  that $\psi_f[-1]$ is a perverse functor, i.e.,  takes perverse sheaves to perverse sheaves, it is not so easy to find proofs that $\phi_f[-1]$ is a perverse functor. Let $j$ denote the inclusion of $f^{-1}(0)$ into $X$. Let $\Pdot$ be a perverse sheaf on $X$. By 1) and 2) above, $j^*\Pdot\in {}^{{}^\mu}\mathbf D^{{}^{\leqslant 0}}(f^{-1}(0))$ and $j^!\Pdot\in{}^{{}^\mu}\mathbf D^{{}^{\geqslant 0}}(f^{-1}(0))$. Knowing that $\psi_f[-1]\Pdot$ is perverse, and applying perverse cohomology to the distinguished triangles
$$\begin{CD}
\phi_f[-1]\Pdot@>\ \operatorname{var}\ >>\psi_f[-1]\Pdot\longrightarrow j^![1]\Pdot\rightarrow\phi_f\Pdot\end{CD}
$$
and
$$\begin{CD}
 j^*[-1]\Pdot\rightarrow\psi_f[-1]\Pdot@>\ \ \operatorname{can}\ \ >>\phi_f[-1]\Pdot\rightarrow  j^*\Pdot,\end{CD}$$
 
 \medskip
 
\noindent one finds, respectively, that ${}^{\mu}\hskip -0.02in H^i(\phi_f[-1]\Pdot)$ is zero if $i\leqslant -1$ and if $i\geqslant 1$, i.e., $\phi_f[-1]\Pdot$ is perverse.

\vskip .3in

If the base ring is a field, then the functor 
$\pcoh$ also commutes with  Verdier
dualizing; that is, there is a  natural isomorphism $$\mathcal D\circ \pcoh\cong \pcoh\circ\mathcal
D.$$

\medskip

Intersection cohomology has a functorial definition related to $\pcoh$. Let $X$ be a  purely $n$-dimensional complex analytic set, and let $\mathcal L$ be a local system (in degree $0$) on a smooth, open dense subset, $\mathcal U$, of the smooth part of $X$. Let $i$ denote the inclusion of $\mathcal U$ into $X$. Then, by applying the functor $\pcoh$ to the canonical map $i_!\mathcal L^\bullet[n]\rightarrow Ri_*\mathcal L^\bullet[n]$, one obtains a morphism $\omega:\pcoh\big(i_!\mathcal L^\bullet[n]\big)\rightarrow \pcoh\big(Ri_*\mathcal L^\bullet[n]\big)$ in the category of perverse sheaves on $X$. The image of $\omega$ in this Abelian category is naturally isomorphic to the intersection cohomology complex $\mathbf {IC}_{{}_X}^\bullet(\mathcal L)$.

If $X$ is not pure-dimensional, one can still use essentially the construction above, except that $\mathcal L^\bullet[n]$ needs to be replaced by a ``shifted piecewise local system'', i.e., one takes the image of $\omega:\pcoh\big(i_!\Pdot\big)\rightarrow \pcoh\big(Ri_*\Pdot\big)$, where $\Pdot$ is a complex on $\mathcal U$, whose restriction to each connected component $C$ of $\mathcal U$ is a local system on $C$, shifted by $\dim C$ (i.e., shifted into degree $-\dim C$); note that such a shifted piecewise local system $\Pdot$ is a perverse sheaf on $\mathcal U$.

More generally, suppose that $X$ is a complex analytic space, that $Y$ is a closed complex analytic proper subspace of $X$, and that $j:Y\rightarrow X$ and $i:X-Y\rightarrow X$ are the inclusions. Let ${\operatorname{Perv}}^X(X-Y)$ denote the full subcategory of ${\operatorname{Perv}}(X-Y)$ consisting of perverse sheaves which are the restrictions of constructible complexes on $X$ (in the algebraic setting, this would be all of ${\operatorname{Perv}}(X-Y)$). For $\Pdot$ in ${\operatorname{Perv}}^X(X-Y)$, $i_!\Pdot$ and $Ri_*\Pdot$ are in $\mathbf D^b_c(X)$, and the functor from ${\operatorname{Perv}}^X(X-Y)$ to ${\operatorname{Perv}}(X)$ which takes $\Pdot$ to the image (in ${\operatorname{Perv}}(X)$) of the canonical map ${}^{\mu}\hskip -0.02in H^0(i_!\Pdot)\rightarrow {}^{\mu}\hskip -0.02in H^0(Ri_*\Pdot)$ is denoted by $i_{!*}$, and is called the {\bf intermediate extension} functor. 

The perverse sheaf $i_{!*}\Pdot$ is the unique (up to isomorphism) perverse extension $\Idot$ of $\Pdot$ from $X-Y$ to $X$ satisfying any/all of the three equivalent conditions:

\medskip

\noindent 1) ${}^{\mu}\hskip -0.02in H^0(j^*\Idot)=0$ and ${}^{\mu}\hskip -0.02in H^0(j^!\Idot)=0$;

\medskip

\noindent 2) $j^*[-1]\Idot\in {}^{{}^\mu}\mathbf D^{{}^{\leqslant 0}}(Y)$ and $j^![1]\Idot\in {}^{{}^\mu}\mathbf D^{{}^{\geqslant 0}}(Y)$;

\medskip

\noindent 3) $\Idot$ has no non-trivial subobjects or quotient objects in ${\operatorname{Perv}}(X)$ whose support is contained in $Y$.

\medskip

\noindent Furthermore, if $Y=V(f)$ where $f$ is a complex analytic function on $X$ (which is not identically zero), then these are equivalent to:

\medskip

\noindent 4) $j^*[-1]\Idot$ and $j^![1]\Idot$ are perverse.

\bigskip

If $X-Y$ is dense in $X$, then the intermediate extension takes intersection cohomology, with possibly local system coefficients, on $X-Y$ to intersection cohomology on $X$, with the same local system coefficients.

\bigskip

Suppose now that we have $f:X\rightarrow\mathbb C$, $i$ denotes the inclusion of $X-f^{-1}(0)$ into $X$, and $j$ denotes the inclusion of $f^{-1}(0)$ into $X$. Further suppose that $\Idot=i_{!*}\Pdot$, where $\Pdot\in {\operatorname{Perv}}^X(X-f^{-1}(0))$. Then we obtain that the two canonical distinguished triangles yield short exact sequences in ${\operatorname{Perv}}(f^{-1}(0))$:
$$\begin{CD}
0\rightarrow\phi_f[-1]\Idot@>\ \operatorname{var}\ >>\psi_f[-1]\Idot\longrightarrow j^![1]\Idot\longrightarrow 0\end{CD}
$$
and
$$\begin{CD}
0\rightarrow j^*[-1]\Idot\rightarrow\psi_f[-1]\Idot@>\ \  \operatorname{can}\ \ >>\phi_f[-1]\Idot\rightarrow 0.\end{CD}$$

\medskip

\noindent Thus, the image, in ${\operatorname{Perv}}(f^{-1}(0))$, of the endomorphism $\operatorname{var}\circ \operatorname{can} = \operatorname{id}-T_f$ on $\psi_f[-1]\Idot$ is isomorphic to $\phi_f[-1]\Idot$.

\bigskip

Let $\mathbf F^\bullet$  be a bounded  complex of sheaves on $X$ which is constructible with respect to a connected Whitney
stratification $\{S_\alpha\}$ of $X$, and let $d_\alpha:=\operatorname{dim}S_\alpha$. Then,  ${}^{\mu}\hskip -0.02in H^0(\Fdot)$ is
also constructible with respect to
$\mathcal S$, and $\big({}^{\mu}\hskip -0.02in H^0(\Fdot)\big)_{|_{\mathbb N_\alpha}}[-d_\alpha]$ is naturally isomorphic to
${}^{\mu}\hskip -0.02in H^0(\Fdot_{|_{\mathbb N_\alpha}}[-d_\alpha])$, where $\mathbb N_\alpha$ denotes a normal slice to $S_\alpha$.

Let
$S_{\operatorname{max}}$ be a maximal stratum contained in the support of $\Fdot$, and let  $m=\dm S_{\operatorname{max}}$. Then,
$\left(\pcoh(\Fdot)\right)_{|_{S_{\operatorname{max}}}}$ is isomorphic (in the derived category) to the complex
which has $\left(\mathbf H^{-m}(\Fdot)\right)_{|_{S_{\operatorname{max}}}}$ in degree
$-m$ and zero in all other degrees. 

In particular, $\operatorname{supp}\Fdot = \bigcup_i \operatorname{supp}{}^{\mu}\hskip -0.02in H^i(\Fdot)$, and if $\Fdot$
is supported on an isolated point, $\mathbf q$, then
$H^0(\pcoh (\mathbf F^\bullet))_{\mathbf q}\cong H^0(\mathbf F^\bullet)_{\mathbf q}.$ From this, and the fact that perverse cohomology commutes with nearby 
and vanishing cycles shifted by $-1$, one easily concludes that, at all points $\mathbf x\in X$,

$$
\chi(\Fdot)_\mathbf x = \sum_k(-1)^k\chi\big({}^{\mu}\hskip -0.02in H^k(\Fdot)\big)_\mathbf x.
$$

\bigskip

\noindent{\bf Switching Coefficients}

\medskip

Suppose that the base ring $R$ is a PID. For each prime ideal $\mathfrak p$ of $R$, let $k_{\mathfrak p}$ denote the field of
fractions of $R/\mathfrak p$, i.e., $k_0$ is the field of fractions of $R$, and for $\mathfrak p\neq 0$, $k_{\mathfrak p} = R/\mathfrak p$.
There are the obvious functors $\delta_{\mathfrak p}: \mathbf D^b_c(R_{{}_X})\rightarrow \mathbf D^b_c({(k_{\mathfrak
p})}_{{}_X})$, which sends $\Fdot$ to $\Fdot\lotimes (k_{\mathfrak p})^\bullet_{{}_X}$, and $\epsilon_{\mathfrak p}: \mathbf
D^b_c({(k_{\mathfrak p})}_{{}_X})\rightarrow\mathbf D^b_c(R_{{}_X})$, which considers $k_{\mathfrak p}$-vector spaces as $R$-modules. 

If $\Adot$ is a
complex of $k_{\mathfrak p}$-vector spaces, we may consider the perverse cohomology of $\Adot$,
${}^{\mu}\hskip -0.02in  H^i_{{}_{k_{\mathfrak p}}}(\Adot)$, or the perverse cohomology of $\epsilon_{\mathfrak p}(\Adot)$, which we denote by
${}^{\mu}\hskip -0.02in  H^i_{{}_R}(\Adot)$. If
$\Adot\in\mathbf D^b_c({(k_{\mathfrak p})}_{{}_X})$ and 
$S_{\operatorname{max}}$ is a maximal stratum contained in the support of $\Adot$, then there is a canonical isomorphism
$$\epsilon_{\mathfrak p}\big(({}^{\mu}\hskip -0.02in  H^i_{{}_{k_{\mathfrak p}}}(\Adot))_{|_{S_{\operatorname{max}}}}\big)\cong ({}^{\mu}\hskip -0.02in 
H^i_{{}_{R}}(\Adot))_{|_{S_{\operatorname{max}}}};$$
in particular, $\operatorname{supp}{}^{\mu}\hskip -0.02in  H^i_{{}_{k_{\mathfrak p}}}(\Adot) = \operatorname{supp}{}^{\mu}\hskip -0.02in 
H^i_{{}_{R}}(\Adot)$.

If
$\Fdot\in\mathbf D^b_c({R}_{{}_X})$, 
$S_{\operatorname{max}}$ is a maximal stratum contained in the support of $\Fdot$, and $\mathbf x\in S_{\operatorname{max}}$, 
then for some prime ideal
$\mathfrak p\subset R$ and for some integer $i$, $H^i(\Fdot)_\mathbf x\otimes k_{\mathfrak p}\neq 0$; it follows that
$S_{\operatorname{max}}$ is also a maximal stratum in the support of $\Fdot\lotimes {(k_{\mathfrak p})}^\bullet_{{}_X}$.
 Thus, 
$$
\operatorname{supp}\Fdot = \bigcup_{\mathfrak p} \operatorname{supp}(\Fdot\lotimes
{(k_{\mathfrak p})}^\bullet_{{}_X})
$$
and so
$$
\operatorname{supp}\Fdot=\bigcup_{i, \mathfrak p}
\operatorname{supp}{}^{\mu}\hskip -0.02in  H^i_{{}_{k_{\mathfrak p}}}(\Fdot\lotimes {(k_{\mathfrak p})}^\bullet_{{}_X}),
$$
where the boundedness and constructibility of $\Fdot$ imply that this union is locally finite.


\begin{thebibliography}{G-M3}


\bibitem{A'C}   N. A'Campo,  {\it Le nombre de Lefschetz d'une
monodromie}, Proc. Kon. Ned. Akad. Wet., Series A, vol. 76,
1973, 113--118  \vskip .1in

\bibitem{BBD} A. Beilinson, J. Bernstein, and P. Deligne, {\it Faisceaux Pervers}, Ast\'erisque vol. 100, Soc. Math de France, 1983
 \vskip .1in

\bibitem{Bo}  A. Borel et al., 
{\it Intersection Cohomology}, Progress in Math. 50, Birkh\"auser, 1984.
 \vskip .1in

\bibitem{BMM}   J. Brian\c con, P. Maisonobe, and M. Merle, {\it  Localisation de
syst\`emes diff\'erentiels, stratifications de Whitney et condition de Thom}, Invent. Math., vol. 117, 1994, 531--550  \vskip .1in


\bibitem{Br}  J. Brylinski,    {\it Transformations canoniques,
Dualit\'e projective, Th\'eorie de Lefschetz, Transformations de Four\-ier et
sommes trigonom\'etriques}, Ast\'erisque, vol. 140,
Soc. Math de France, 1986  \vskip .1in


\bibitem{De}  P. Deligne, {\it Comparaison avec la th\'eorie
transcendante}, S\'eminaire de  g\'eom\'etrie alg\'ebrique du
Bois-Marie, SGA 7 II, Springer Lect. Notes, vol. 340, 1973   \vskip .1in

\bibitem{Di} A. Dimca,
{\it Sheaves in Topology},  Universitext, Springer-Verlag, 2004
 \vskip .1in

\bibitem{Gi} V. Ginsburg, {\it Characteristic Varieties and
Vanishing Cycles}, Invent. Math., vol. 84, 1986, 327--403   
\vskip .1in


\bibitem{G-M1}  M. Goresky and R. MacPherson,  {\it Morse Theory and Intersection Homology}, Analyse et Topologie sur les Espaces Singuliers, Ast\'erisque, vol. 101,  Soc.
Math. France, 1983, 135--192     \vskip .1in


\bibitem{G-M2} M. Goresky and R. MacPherson, {\it Stratified Morse Theory},
Springer-Verlag, Berlin, Ergebnisse der Math. vol. 14, 1988
\vskip .1in

\bibitem{G-M3}  M. Goresky and R. MacPherson, {\it Intersection homology II}, Invent. Math., vol. 71, 1983,  77-129\vskip .1in


\bibitem{G-M4}  M. Goresky and R. MacPherson,  {\it Intersection homology
theory}, Topology, vol. 19, 1980, 135--162  \vskip .1in

\bibitem{Gr} P. P. Grivel, {\it Les foncteurs de la cat\'egorie des faisceaux associ\'es \`a une application
continue}, Intersection Cohomology, Prog. in Math., vol. 50, 1984, 183--207
\vskip .1in


\bibitem{H}  R. Hartshorne,  {\it Residues and Duality}, Springer Lecture Notes, vol. 20, Springer-Verlag, 1966
 \vskip .1in

\bibitem {I}  B. Iversen, {\it Cohomology of Sheaves}, Springer-Verlag, 1986  \vskip .1in

\bibitem {K-S}  M. Kashiwara and P. Schapira, {\it Sheaves on
Manifolds}, Grund. der math. Wiss., vol. 292, Springer-Verlag, 1990
 \vskip .1in

\bibitem{Le1}  L\^e D. T.,    {\it Morsification of $\mathbf D$-Modules}, Bol. Soc. Mat. Mexicana
(3) vol. 4, 229--248, 1998

\vskip .1in

\bibitem{Le2}  L\^e D. T.,    {\it Sur les cycles
\'evanouissants des espaces analytiques}, C.R. Acad.
Sci. Paris, Ser. A, vol. 288, 283--285, 1979   \vskip .1in


\bibitem{Lo}  E. Looijenga, {\it Isolated Singular Points on Complete Intersections}, Cambridge Univ. Press, 1984
 \vskip .1in

\bibitem{Mac1}   R. MacPherson, {\it Global Questions in the Topology of Singular Spaces},
Proc. Internat. Congress of Math., Warsaw, 213--235, 1983   \vskip .1in

\bibitem{Mac2}  R. MacPherson, {\it Intersection Homology and Perverse Sheaves}, unpublished AMS notes, 1--161, 1990   \vskip .1in

\bibitem{M-V}   R. MacPherson and K. Vilonen, {\it Elementary construction of perverse sheaves}, Invent. Math., vol. 84, 403--435, 1986   \vskip .1in

\bibitem{Nee}   A. Neeman, {\it Triangulated Categories}, Annals of Math. Studies, vol. 148, 2001   \vskip .1in

\bibitem{P}   A. Parusi\'nski,  {\it Limits of Tangent Spaces to fibers and the $w_f$ Condition}, Duke Math. J., vol. 72, 99--108, 1993     \vskip .1in


\bibitem{Sa}   M. Saito, {\it Modules de Hodge polarisables}, Publ. RIMS, Kyoto Univ., vol. 24, 849--995, 1988   \vskip .1in

\bibitem{Sch}  J. Sch\"urmann, {\it Topology of Singular Spaces and Constructible Sheaves}, Monografie Matematyczne, vol. 63, Birkh\"auser, 2003
    \vskip .1in

\bibitem{V}   J. L. Verdier, {\it Cat\'egories d\'eriv\'ees}, \'Etat 0, SGA $4 \frac12$,
 Lecture Notes in Math., vol. 569, 262--311, 1977   

\end{thebibliography}
\end{document}